\documentclass[11pt,reqno]{amsart}

\usepackage{geometry}

\title{Hochschild Homology of analytic rings and $B_{\on{dR}}^+$-cohomology}
\author{Zhiyou Wu}

\usepackage{amsmath,amssymb,fancyhdr,url,amsthm,mathtools,mathrsfs,setspace,stmaryrd}
\usepackage[pdftex]{graphicx}
\usepackage{tikz-cd}
\usetikzlibrary{cd}
\usepackage[mathscr]{euscript}
\usepackage{tikz}
\usepackage{dsfont}
\usetikzlibrary{matrix}

\makeatletter
\newcommand*{\rom}[1]{\expandafter\@slowromancap\romannumeral #1@}
\makeatother

\usepackage{relsize} 
\usepackage[bbgreekl]{mathbbol} 
\usepackage{amsfonts} 
\DeclareSymbolFontAlphabet{\mathbb}{AMSb} 
\DeclareSymbolFontAlphabet{\mathbbl}{bbold} 
\newcommand{\prism}{{\mathlarger{\mathbbl{\Delta}}}}

\usepackage[toc,page]{appendix}



\tikzset{
    labl/.style={anchor=north, rotate=90, inner sep=.8mm}
}

\newtheorem{proposition}{Proposition}[section]
\newtheorem{theorem}[proposition]{Theorem}
\newtheorem{example}[proposition]{Example}
\newtheorem{corollary}[proposition]{Corollary}
\newtheorem{definition}[proposition]{Definition}
\newtheorem{remark}[proposition]{Remark}
\newtheorem{lemma}[proposition]{Lemma}

\newcommand{\nc}{\newcommand}
\nc{\on}{\operatorname}

\begin{document}

\maketitle

\begin{abstract}
    We use Hochschild homology to construct $\on{B}_{\on{dR}}^+$-cohomology theories  for Clausen-Scholze's analytic rings over perfectoid rings, and Breuil-Kisin type cohomology theories for analytic rings over discrete valued non-archimedean fields.  We establish basic results on HKR filtration of Hochschild homology of  analytic rings, and the Hodge completed derived de Rham complex of analytic rings, following \cite{Raksit2026}. 
\end{abstract}

\tableofcontents

\section{Introduction}

The first main result of this paper is to give a Hochschild homology construction of the $\on{B}_{\on{dR}}^+$-cohomology  for smooth rigid analytic spaces over a perfectoid field. There are actually two  closely related $\on{B}_{\on{dR}}^+$-cohomology theories, first constructed in  \cite{SCHOLZE2013} and \cite[section 13]{Bhatt2018} respectively, and we are able to reconstruct both of them directly from the Hochschild homology. Both of the two constructions have been generalized by Haoyang Guo in \cite{Guo2025Crystalline} and  \cite{Guo2026Prismatic} respectively, and we adopt this perspective in this work.

We work with analytic rings in the framework of Clausen-Scholze, and the construction generalizes easily to an arbitrary analytic ring defined over a perfectoid field, so we can in particular generalize the $\on{B}_{\on{dR}}^+$-cohomology theories to derived adic spaces (in the sense of \cite{RodriguezCamargo2024DeRham}) defined over perfectoid fields. 

On the other hand, with similar constructions, we can construct  Breuil-Kisin type cohomology theories for analytic rings defined over discrete valued non-archimedean fields, and in particular for rigid analytic varieties defined over discrete valued non-archimedean fields. 

The construction parallels perfectly with (and is clearly motivated by) the work of \cite{BMS2}. The main difference is that we are in the characteristic 0 case, so there is no difference between algebra and topology, namely the sphere spectrum becomes rational numbers after localization at all primes, which means that we do not need to consider topological Hochschild homology as in \cite{BMS2}, and the more classical algebraic Hochschild homology suffices. The price to pay in the characteristic 0 situation is that the $p$-adic completion does not work naively, and the  completion procedure is done, as is becoming standard nowadays, in the solid setting of Clausen-Scholze. Another difference is that in the algebraic setting, we can use the more classical Hochschild-Kostant-Rosenberg filtration instead of the even filtration, and one of the main technical result is to establish the basic results of HKR filtration for analytic rings. We note that it is observed in \cite{Wagner2026Habiro} that even filtration (more precisely, the version developed in \cite{E1even}) is subtle to generalize to analytic rings, whereas we prove that the HKR filtration generalizes to complete generality.

As a preliminary input for the construction of cohomology theories, we construct  HKR filtered Hochschild homology and  Hodge completed derived de Rham cohomology for analytic rings. 
We follow the strategy of \cite{Raksit2026}
closely here. It is possible to construct the HKR filtration directly by animating the polynomial case as is done in the classical approach, however, the filtered $\mathbb{S}^1$-action and the resulting filtration for homotopy fixed points and Tate construction are not easy to derive from the naive approach.

Let us explain our work on the generalization of \cite{Raksit2026} to analytic rings. There are two steps here, the first one is to generalize it to condensed rings, i.e. sheaves of rings on the light profinite sets, and the second step is to take care of the analytic structures. For the first step, we can not cite directly the general construction of \cite{Raksit2026} 
since condensed rings does not come from a derived algebraic context as defined in $loc.cit.$.  The problem is the lack of compact projective generators in condensed abelian groups. However, this is also not hard to do since we can just consider the condensed derived rings (together with the filtered and graded ones), i.e. sheaves on light profinite sets with values in derived rings,  and bootstrap the constructions in $loc.cit.$. We need to prove that it has the right monadic properties, but this is all straightforward. 

For the second step on analytic structures, the problem is more subtle. We do not define a theory of analytic derived rings, since there are multiple choice of $t$-structures on non-connective rings and it is not clear to the author which one is natural to choose. Instead, we define a minimal framework of filtered (and graded) analytic derived rings that contains both Hochschild homology (Hodge completed derived de Rham complex) of analytic animated rings and the filtered circle $\mathbb{T}_{\on{Fil}}^{\vee}$, which is non-connective and is the main reason we need to consider derived rings even if we only care about Hochschild homology of animated rings. The framework is not optimal, it is the smallest one where we can carry the  usual constructions of analytic animated rings to the filtered (graded) derived algebra context. In particular, we need to prove that it has completed and noncompleted version, and there is a left adjoint completion functor from uncompleted ones to completed ones. 

Let us be more precise on what is proved and why we need to proceed the way we present.  Let $\mathcal{A} \longrightarrow \mathcal{B}$ be a map of analytic animated rings, and we want to construct the HKR filtration on Hochschild homology $\on{HH}(\mathcal{B}/\mathcal{A}) $, which should have graded pieces $\on{Sym}^{\bullet}_{\mathcal{B}} \Sigma \mathbb{L}_{\mathcal{B}/\mathcal{A}}$. A natural approach is to first construct the filtration for Hochschild homology $\on{HH}(\mathcal{B}^{\triangleright}/\mathcal{A}^{\triangleright})$ of underlying condensed animated rings, and then take the completion functor from $\mathcal{B}^{\triangleright}$-modules to  $\mathcal{B}$-modules. The HKR filtration on $\on{HH}(\mathcal{B}^{\triangleright}/\mathcal{A}^{\triangleright})$ is a direct algebraic construction which has graded pieces $\on{Sym}^{\bullet}_{\mathcal{B}^{\triangleright}} \Sigma \mathbb{L}_{\mathcal{B}^{\triangleright}/\mathcal{A}^{\triangleright}}$, and we need to take $\mathcal{B}$-completion to obtain the desired $\on{Sym}^{\bullet}_{\mathcal{B}} \Sigma \mathbb{L}_{\mathcal{B}/\mathcal{A}}$. This can be done and does give the correct answer, but it is not clear that the filtration is $\mathbb{S}^1$-equivariant, let alone it has a filtered circle action.  The problem is that the circle action is only $\mathcal{A}^{\triangleright}$-linear, and we are taking $\mathcal{B}$-completion, there is a priori no reason why the circle action should pass to the $\mathcal{B}$-completion. 

In order to prove this, we set up a framework of (completed and uncompleted)  filtered analytic derived rings, and view  
$\on{HH}(\mathcal{B}^{\triangleright}/\mathcal{A}^{\triangleright})$
as an uncompleted analytic filtered derived ring with a filtered circle action. The $\mathcal{B}$-completion is replaced by the completion functor from uncompleted filtered  analytic derived rings to complete filtered analytic derived rings, which preserves the filtered circle action and its effect on graded pieces is simply the $\mathcal{B}$-completion. This works as expected and the computation of the induced filtration on homotopy fixed points and Tate construction follows directly. 

The precise theorem on $\on{B}_{\on{dR}}^+$-cohomology theories is as follows.

\begin{theorem}
    Let $(K, K^+)$ be a perfectoid ring over $(\mathbb{Q}_p, \mathbb{Z}_p)$, and $(K,K^+)_{\blacksquare} \longrightarrow \mathcal{R}$
    a map of analytic animated rings,
    then we have two cohomology theories 
    $ R\Gamma_{\on{inf}}(\mathcal{R}/\mathbb{B}_{\on{dR}}^+(K) ) $ 
    and 
    $ \widehat{\prism}_{\mathcal{R}/\mathbb{B}_{\on{dR}}^+(K)   } $
    of $\mathcal{R}$,  with coefficient in the de Rham period ring $\mathbb{B}_{\on{dR}}^+(K)$ satisfying the following properties.  

    (1)  We have the de Rham comparison with the derived de Rham complex of $\mathcal{R} $ over $(K,K^+)_{\blacksquare}$, 
    \[
     R\Gamma_{\on{inf}}(\mathcal{R}/\mathbb{B}_{\on{dR}}^+(K) ) \otimes_{(\mathbb{B}_{\on{dR}}^+(K),\mathbb{Z}_p)_{\blacksquare}  } K 
     \cong 
     \widehat{\on{dR}}_{\mathcal{R}/ (K, K^+)_{\blacksquare}}.
    \]

   (2) There exists a Hodge-Tate filtration on
   \[
   \overline{\prism}_{\mathcal{R}/\mathbb{B}_{\on{dR}}^+(K)   } 
    := \widehat{\prism}_{\mathcal{R}/\mathbb{B}_{\on{dR}}^+(K)   } \otimes_{ (\mathbb{B}_{\on{dR}}^+(K), \mathbb{Z}_p )_{\blacksquare}} K
   \]
    with graded pieces
    \[
    \on{gr}_{\on{HT}}^{\bullet} \overline{\prism}_{\mathcal{R}/\mathbb{B}_{\on{dR}}^+(K)   }  
    \cong 
    \Sigma^{-2*}\on{Sym}_{\mathcal{R}}^{\bullet} (\Sigma \mathbb{L}_{\mathcal{R}/(K,K^+)_{\blacksquare}}[1]).
    \]

(3) There exists a natural filtration on 
$ R\Gamma_{\on{inf}}(\mathcal{R}/\mathbb{B}_{\on{dR}}^+(K) ) $
with graded pieces
\[
\on{gr}^i R\Gamma_{\on{inf}}(\mathcal{R}/\mathbb{B}_{\on{dR}}^+(K) )
\cong 
F^{i}_{\on{HT}}\overline{\prism}_{\mathcal{R}/\mathbb{B}_{\on{dR}}^+(K)  }.
\]

(4) There exists a 
natural comparison map
\[
     R\Gamma_{\on{inf}}(\mathcal{R}/\mathbb{B}_{\on{dR}}^+(K) ) 
     \longrightarrow 
     \widehat{\prism}_{\mathcal{R}/\mathbb{B}_{\on{dR}}^+(K)   }
   \]
   that induces an equivalence 
   \[
     R\Gamma_{\on{inf}}(\mathcal{R}/\mathbb{B}_{\on{dR}}^+(K) ) [\xi^{-1}]
     \cong 
     \widehat{\prism}_{\mathcal{R}/\mathbb{B}_{\on{dR}}^+(K)   }[\xi^{-1}].
   \]

   (5) When $\mathcal{R} = (R, R^+)_{\blacksquare}$ for a smooth affinoid $(R,R^+)$ ring over a perfectoid field $(K,K^+)$, then 
   $ R\Gamma_{\on{inf}}(\mathcal{R}/\mathbb{B}_{\on{dR}}^+(K) ) $ 
    and 
    $ \widehat{\prism}_{\mathcal{R}/\mathbb{B}_{\on{dR}}^+(K)   } $
    are the same as the crystalline and prismatic $\on{B}_{\on{dR}}^+$-cohomology theory defined in \cite{Guo2025Crystalline} and \cite{Guo2026Prismatic} respectively, which also coincides with the one constructed in \cite[section 13]{Bhatt2018} and the proetale $\on{B}_{\on{dR}}^+$-cohomology in \cite{SCHOLZE2013} respectively. 
\end{theorem}

We can construct these cohomology theories using filtration on Hochschild homology following the strategy of Bhatt-Morrow-Scholze (\cite{BMS2}). 

   \begin{theorem} 
   Let $(K, K^+)$ be a perfectoid ring over $(\mathbb{Q}_p, \mathbb{Z}_p)$, and $(K,K^+)_{\blacksquare} \longrightarrow \mathcal{R}$
    a map of analytic animated rings,
   there exists filtration on $\on{HH}(\mathcal{R}/(\mathbb{Q}_p, \mathbb{Z}_p)_{\blacksquare})^{t\mathbb{S}^1}$
and $\on{HH}(\mathcal{R}/(\mathbb{Q}_p, \mathbb{Z}_p)_{\blacksquare})^{h\mathbb{S}^1}$
with graded pieces
\[
\on{gr}^{i}
\on{HH}(\mathcal{R}/(\mathbb{Q}_p, \mathbb{Z}_p)_{\blacksquare})^{t\mathbb{S}^1}
= 
\Sigma^{2i}  R\Gamma_{\on{inf}}(\mathcal{R}/\mathbb{B}_{\on{dR}}^+(K) ) 
\]
\[
\on{gr}^{i}
\on{HH}(\mathcal{R}/(\mathbb{Q}_p, \mathbb{Z}_p)_{\blacksquare})^{h\mathbb{S}^1}
= \Sigma^{2i} R\Gamma_{\on{inf}}(\mathcal{R}/\mathbb{B}_{\on{dR}}^+(K) )^{\geq i}
\]
such that the canonical map from homotopy fixed points to Tate construction induces the filtration maps on graded pieces.

 Moreover, there exists natural HKR filtration on Hochschild homology with graded pieces 
\[
\on{gr}_{\on{HKR}}^{i}\on{HH}(\mathcal{R}/(\mathbb{Q}_p, \mathbb{Z}_p)_{\blacksquare})
\cong
\Sigma^{2i}
F^{i}_{\on{HT}}\overline{\prism}_{\mathcal{R}/\mathbb{B}_{\on{dR}}^+(K)  }.
\] 
\end{theorem}

The construction of $\widehat{\prism}_{\mathcal{R}/\mathbb{B}_{\on{dR}}^+(K)   } $
is an analouge of the Frebenius descent in \cite[section 11]{BMS2}. 

\begin{theorem}
 Let $(K, K^+)$ be a perfectoid ring over $(\mathbb{Q}_p, \mathbb{Z}_p)$, and $(K,K^+)_{\blacksquare} \longrightarrow \mathcal{R}$
    a map of analytic animated rings, then we have
  \[
    \on{gr}^{\bullet}
\on{HH}(\mathcal{R}/(\mathbb{Q}_p, \mathbb{Z}_p)_{\blacksquare})^{h\mathbb{S}^1}[\Tilde{u}^{-1}]
\cong 
\widehat{\prism}_{\mathcal{R}/\mathbb{B}_{\on{dR}}^+(K)   } [ \Tilde{u}, \Tilde{u}^{-1}]
    \]
where  $\Tilde{u}$ of homological degree 2 and graded degree 1
is a choice of generator 
\[
\on{gr}^{\bullet}
\on{HH}((K, K^+)_{\blacksquare}/(\mathbb{Q}_p, \mathbb{Z}_p)_{\blacksquare})^{h\mathbb{S}^1}
\cong \mathbb{B}_{\on{dR}}^+(K) [\Tilde{u},v]/(\Tilde{u}v-\xi).
\]

   Moreover, the localization map 
   \[
   \on{HH}(\mathcal{R}/(\mathbb{Q}_p, \mathbb{Z}_p)_{\blacksquare})^{h\mathbb{S}^1}
   \longrightarrow
   \on{HH}(\mathcal{R}/(\mathbb{Q}_p, \mathbb{Z}_p)_{\blacksquare})^{h\mathbb{S}^1}[\Tilde{u}^{-1}]
   \]
   induces a comparison map 
   \[
     R\Gamma_{\on{inf}}(\mathcal{R}/\mathbb{B}_{\on{dR}}^+(K) ) 
     \longrightarrow 
     \widehat{\prism}_{\mathcal{R}/\mathbb{B}_{\on{dR}}^+(K)   }.
   \]
  
\end{theorem}

We can also  construct a variant of the $\on{B}_{\on{dR}}^+ $-cohomology for
 analytic animated rings defined over a discrete valued non-archimedean field, similar to  Breuil-Kisin cohomology theory. 

\begin{theorem}
    Let $k$ be a perfect field of characteristic $p$, and $K$ a finite ramified extension of $W(k)[p^{-1}]$ with uniformizer $\omega$. Let $(K,\mathbb{Z}_p)_{\blacksquare} \longrightarrow \mathcal{R}$
    a map of analytic animated rings, 
    and 
    \[
    \mathfrak{S} := \underline{W(k)[p^{-1}]}[T]^{\wedge}_E
    \]
    where $E$ is a minimal polynomial of $\omega$ over $W(k)[p^{-1}]$. 
    We have two cohomology theories 
    $ R\Gamma_{\on{inf}}(\mathcal{R}/\mathfrak{S} ) $ 
    and 
    $ \widehat{\prism}_{\mathcal{R}/\mathfrak{S}  } $
    of $\mathcal{R}$,  with coefficient in the de Rham period ring $\mathfrak{S}$ satisfying the following properties.  

    (1)  We have the de Rham comparison with the derived de Rham complex of $\mathcal{R} $ over $(K,\mathbb{Z}_p)_{\blacksquare}$, 
    \[
     R\Gamma_{\on{inf}}(\mathcal{R}/\mathfrak{S} ) \otimes_{(\mathfrak{S},\mathbb{Z}_p)_{\blacksquare}  } K 
     \cong 
     \widehat{\on{dR}}_{\mathcal{R}/ (K, \mathbb{Z}_p)_{\blacksquare}}.
    \]

   (2) There exists a Hodge-Tate filtration on
   \[
   \overline{\prism}_{\mathcal{R}/\mathfrak{S}  } 
    := \widehat{\prism}_{\mathcal{R}/\mathfrak{S}  } \otimes_{ (\mathfrak{S}, \mathbb{Z}_p )_{\blacksquare}} K
   \]
    with graded pieces
    \[
    \on{gr}_{\on{HT}}^{\bullet} \overline{\prism}_{\mathcal{R}/\mathfrak{S}  }  
    \cong 
    \Sigma^{-2*}\on{Sym}_{\mathcal{R}}^{\bullet} (\Sigma \mathbb{L}_{\mathcal{R}/(K,\mathbb{Z}_p)_{\blacksquare}}[1]).
    \]

(3) There exists a natural filtration $R\Gamma_{\on{inf}}(\mathcal{R}/\mathfrak{S} )^{\geq \bullet}$ on 
$ R\Gamma_{\on{inf}}(\mathcal{R}/\mathfrak{S} ) $
with graded pieces
\[
\on{gr}^i R\Gamma_{\on{inf}}(\mathcal{R}/\mathfrak{S} )
\cong 
F^{i}_{\on{HT}}\overline{\prism}_{\mathcal{R}/\mathfrak{S}  }.
\]

(4) There exists a 
natural comparison map
\[
     R\Gamma_{\on{inf}}(\mathcal{R}/\mathfrak{S} ) 
     \longrightarrow 
     \widehat{\prism}_{\mathcal{R}/\mathfrak{S}  }
   \]
   that induces an equivalence 
   \[
     R\Gamma_{\on{inf}}(\mathcal{R}/\mathfrak{S} ) [E^{-1}]
     \cong 
     \widehat{\prism}_{\mathcal{R}/\mathfrak{S}  }[E^{-1}].
   \]
\end{theorem}

The construction is also similar to the construction of Breuil-Kisin cohomology in \cite{BMS2}, where we take relative  Hochschild homology over a polynomial base. 

   \begin{theorem} 
   Let $K$ be a finite ramified  extension of $W(k)[p^{-1}]$ as in the above theorem,   $(K,\mathbb{Z}_p)_{\blacksquare} \longrightarrow \mathcal{R}$ 
    a map of analytic animated rings,
   there exists filtration on $\on{HH}(\mathcal{R}/(\mathbb{Q}_p[T], \mathbb{Z}_p)_{\blacksquare})^{t\mathbb{S}^1}$
and $\on{HH}(\mathcal{R}/(\mathbb{Q}_p[T], \mathbb{Z}_p)_{\blacksquare})^{h\mathbb{S}^1}$
with graded pieces
\[
\on{gr}^{i}
\on{HH}(\mathcal{R}/(\mathbb{Q}_p[T], \mathbb{Z}_p)_{\blacksquare})^{t\mathbb{S}^1}
= 
\Sigma^{2i}  R\Gamma_{\on{inf}}(\mathcal{R}/\mathfrak{S} ) 
\]
\[
\on{gr}^{i}
\on{HH}(\mathcal{R}/(\mathbb{Q}_p[T], \mathbb{Z}_p)_{\blacksquare})^{h\mathbb{S}^1}
= \Sigma^{2i} R\Gamma_{\on{inf}}(\mathcal{R}/\mathfrak{S} )^{\geq i}
\]
such that the canonical map from homotopy fixed points to Tate construction induces the filtration maps on graded pieces.

 Moreover, there exists natural HKR filtration on Hochschild homology with graded pieces 
\[
\on{gr}_{\on{HKR}}^{i}\on{HH}(\mathcal{R}/(\mathbb{Q}_p[T], \mathbb{Z}_p)_{\blacksquare})
\cong
\Sigma^{2i}
F^{i}_{\on{HT}}\overline{\prism}_{\mathcal{R}/\mathfrak{S}  },
\] 
and we have 
  \[
    \on{gr}^{\bullet}
\on{HH}(\mathcal{R}/(\mathbb{Q}_p[T], \mathbb{Z}_p)_{\blacksquare})^{h\mathbb{S}^1}[\Tilde{u}^{-1}]
\cong 
\widehat{\prism}_{\mathcal{R}/\mathfrak{S}  } [ \Tilde{u}, \Tilde{u}^{-1}]
    \]
where  $\Tilde{u}$ of homological degree 2 and graded degree 1
is a choice of generator 
\[
\on{gr}^{\bullet}
\on{HH}((K, \mathbb{Z}_p)_{\blacksquare}/(\mathbb{Q}_p[T], \mathbb{Z}_p)_{\blacksquare})^{h\mathbb{S}^1}
\cong \mathfrak{S} [\Tilde{u},v]/(\Tilde{u}v-\xi).
\]

  Lastly, the localization map 
   \[
   \on{HH}(\mathcal{R}/(\mathbb{Q}_p[T], \mathbb{Z}_p)_{\blacksquare})^{h\mathbb{S}^1}
   \longrightarrow
   \on{HH}(\mathcal{R}/(\mathbb{Q}_p[T], \mathbb{Z}_p)_{\blacksquare})^{h\mathbb{S}^1}[\Tilde{u}^{-1}]
   \]
   induces a comparison map 
   \[
     R\Gamma_{\on{inf}}(\mathcal{R}/\mathfrak{S} ) 
     \longrightarrow 
     \widehat{\prism}_{\mathcal{R}/\mathfrak{S}  }.
   \]
  
\end{theorem}

\begin{remark}
    We state only the affinoid version, but it is straightforward to glue the construction to obtain cohomology theories for global analytic spaces along Zariski topology for example. 
\end{remark}

\subsection*{Conventions} 
We use the framework of (light) condensed mathematics as developed by Clausen-Scholze, and in particular we follow the treatment of \cite{Notesonsolidgeometry}. We use the notation
\[
\on{Cond} D(\mathbb{Z}) := \on{Shv}(\on{Prof}^{\on{light}}, D(\mathbb{Z})) 
\]
the hypersheaves on the site of light profinite sets with values in $D(\mathbb{Z})$, and similar for condensed rings and etc.  To distinguish between condensed rings and analytic rings, we will generally use $\mathcal{A}$ to denote
analytic rings, and $\mathcal{A}^{\triangleright}$ its underlying condensed ring, whereas $A$ is used for condensed rings without analytic structures. 

\subsection*{Acknowledgment} 
We would like to thank the support from National Key Research and Development Program of China 2025YFA1018000. We would like to also thank the support from BR Talent program E4290101. 

\subsection*{Declaration of AI assistance}
The paper is completely written by human, and the author is responsible for all the content. AI is used mostly for searching literature, but there is one exception, the proof of lemma \ref{fhowufew0u9f0ew}, which is found by ChatGPT. The author has checked the details and present it with my own understanding.

\section{Filtered analytic derived rings}

We define a version of filtered and graded analytic derived rings,  which only serves as technical purpose here. We follow the treatment of \cite{Raksit2026} closely.  Let 
\[
\on{Fil}\on{Cond} D (\mathbb{Z}) := \on{Fun}(\mathbb{Z}^{\on{op}}, \on{Cond} D_{} (\mathbb{Z})) 
\]
\[
\on{Gr}\on{Cond} D (\mathbb{Z}) := \on{Fun}(\mathbb{Z}^{\on{disc}}, \on{Cond} D_{} (\mathbb{Z})) 
\]
be the filtered/graded condensed connective $\mathbb{Z}$-modules. We see easily that  
\[
\on{Fil}\on{Cond} D_{} (\mathbb{Z}) 
\cong 
\on{Cond} \on{Fun}(\mathbb{Z}^{\on{op}},D_{} (\mathbb{Z})) = \on{Cond} \on{Fil}D_{} (\mathbb{Z})
\cong \on{Cond}\on{Ani} \otimes \on{Fil}D_{} (\mathbb{Z})
\]
\[
\on{Gr}\on{Cond} D_{} (\mathbb{Z}) 
\cong 
\on{Cond} \on{Fun}(\mathbb{Z}^{\on{disc}},D_{} (\mathbb{Z})) = \on{Cond} \on{Gr}D_{} (\mathbb{Z})
\cong \on{Cond}\on{Ani} \otimes \on{Gr}D_{} (\mathbb{Z})
\]
where the tensor products are taken in $\on{Pr}^L$. 

Recall (\cite[Construction 4.3.4]{Raksit2026}) that Bhatt-Matthew's theory of derived rings produces  monads $\on{LSym}^{\bullet}= \oplus \on{LSym}^{i}$ on $\on{Fil}D_{} (\mathbb{Z})$  and $\on{Gr}D_{} (\mathbb{Z})$  
that refines the commutative algebra monad $X \mapsto \oplus_n (X^{\otimes n})_{\Sigma_n}$ on $\on{Fil}D_{} (\mathbb{Z})$ and $\on{Gr}D_{} (\mathbb{Z})$, where $\on{Fil}D_{} (\mathbb{Z})$ and $\on{Gr}D_{} (\mathbb{Z})$ are equipped with the Day convolution product. Let 
\[
\on{DAlg}\on{Fil}D_{} (\mathbb{Z}) := \on{LMod_{\on{LSym}^{\bullet}}} \on{Fil}D_{} (\mathbb{Z})
\]
\[
\on{DAlg}\on{Gr}D_{} (\mathbb{Z}) := \on{LMod_{\on{LSym}^{\bullet}}} \on{Gr}D_{} (\mathbb{Z})
\]
which by definition comes with  monadic adjoint functors
\[
\on{Fil}D_{} (\mathbb{Z}) \rightleftarrows
\on{DAlg}\on{Fil}D_{} (\mathbb{Z})
\]
\[
\on{Gr}D_{} (\mathbb{Z}) \rightleftarrows
\on{DAlg}\on{Gr}D_{} (\mathbb{Z})
\]
with the right adjoints the forgetful functors. Further, the forgetful functors factorize through 
\[
U:
\on{DAlg}\on{Fil}D_{} (\mathbb{Z}) \overset{\Theta}{\longrightarrow} \on{CAlg}\on{Fil}D_{} (\mathbb{Z}) \longrightarrow \on{Fil}D_{} (\mathbb{Z})
\]
\[
U:
\on{DAlg}\on{Gr}D_{} (\mathbb{Z}) \overset{\Theta}{\longrightarrow} \on{CAlg}\on{Gr}D_{} (\mathbb{Z}) \longrightarrow \on{Gr}D_{} (\mathbb{Z})
\]
as right adjoint functors, and the functors $\Theta$ preserve both limits and colimits  by \cite[Proposition 4.2.27]{Raksit2026}.

Let 
\[
\on{DAlg}\on{Cond}\on{Fil}D_{} (\mathbb{Z}) := \on{Cond} \on{DAlg}\on{Fil}D_{} (\mathbb{Z}) \cong \on{Cond} \on{Ani} \otimes \on{DAlg} \on{Fil} D(\mathbb{Z}),
\]
\[
\on{DAlg}\on{Cond}\on{Gr}D_{} (\mathbb{Z}) := \on{Cond} \on{DAlg}\on{Gr}D_{} (\mathbb{Z}) \cong \on{Cond} \on{Ani} \otimes \on{DAlg} \on{Gr} D(\mathbb{Z}),
\]
which are hypersheaves on light profinite sets with values in $\on{DAlg}\on{Fil}D_{} (\mathbb{Z})$ and $\on{DAlg}\on{Gr}D_{} (\mathbb{Z})$.
We see that it comes with a monadic adjoint pair
\[
\on{Sym}^{\bullet}:
\on{Cond}\on{Fil}D_{} (\mathbb{Z}) \rightleftarrows
\on{DAlg}\on{Cond}\on{Fil}D_{} (\mathbb{Z}): U
\]
\[
\on{Sym}^{\bullet}:
\on{Cond}\on{Gr}D_{} (\mathbb{Z}) \rightleftarrows
\on{DAlg}\on{Cond}\on{Gr}D_{} (\mathbb{Z}): U
\]
since tensoring with a presentable  category preserves monadic left adjoints in $\on{Pr}^L$, 
\footnote{Given a monadic left adjoint functor $F: C \longrightarrow D$ in $\on{Pr}^L$ with right adjoint $U: D \longrightarrow C$, and let $A \in \on{Pr}^L$, then $\on{id}_A \otimes F: \on{Fun}^R(A^{\on{op}}, C)\cong A\otimes C \longrightarrow A \otimes D\cong \on{Fun}^R(A^{\on{op}}, D)$ has right adjoint $R: \on{Fun}^R(A^{\on{op}}, D) \longrightarrow \on{Fun}^R(A^{\on{op}}, C)$ given by composition with $U$, see \cite[Remark 9.7.5.12]{Kerodon}. $R$ is conservative since $U$ is conservative, and by Barr-Beck-Lurie, we need only to prove that $R$ preserves colimits of diagrams $\gamma: \Delta^{op} \longrightarrow \on{Fun}^R(A^{\on{op}}, D) $ 
such that $R \circ \gamma$ are split simplicial diagrams.

Let $x= \on{colim}_{\Delta^{\on{op}}} R\circ \gamma \in \on{Fun}^R(A^{\on{op}}, C) $.
Since $R \circ \gamma$ is split and split diagrams are universal,  for any object  $a \in A^{\on{op}}$,  $\on{ev}_a \circ R\circ \gamma = U \circ \on{ev}_a \circ \gamma$ 
is also split. Since $U$ is monadic and $U \circ \on{ev}_a \circ \gamma$ is split, $U$ preserves the colimit of $\on{ev}_a \circ \gamma$.   
Moreover, $\on{ev}_a$ preserves the colimit of $R \circ \gamma$ by the universality of split diagrams, so 
\[
\on{ev}_a(x) = \on{colim}_{\Delta^{\on{op}}} \on{ev}_a \circ R \circ \gamma = \on{colim}_{\Delta^{\on{op}}} U \circ \on{ev}_a \circ \gamma
= U (\on{colim}_{\Delta^{\on{op}}}\on{ev}_a \circ \gamma).
\]
Let $d \in \on{Fun}(A^{\on{op}},D)$ be the colimit of $\gamma$ composed with the inclusion $ \on{Fun}^R(A^{\on{op}},D) \subset  \on{Fun}(A^{\on{op}},D)$, then the colimit is computed pointwise, so we have that $\on{ev}_a(d) = \on{colim}_{\Delta^{\on{op}}}\on{ev}_a \circ \gamma $, and 
\begin{equation} \label{xmnfs}
\on{ev}_a(x)= U (\on{colim}_{\Delta^{\on{op}}}\on{ev}_a \circ \gamma)  = U \circ \on{ev}_a (d).
\end{equation}
If $d \in \on{Fun}^R(A^{\on{op}},D)$,  then it is the colimit of $\gamma$, and we have that 
\[
\on{ev}_a(x) = U \circ \on{ev}_a (d) =\on{ev}_a \circ R (d)
\]
for every $a \in A^{\on{op}}$, proving that $x=R(d)$, i.e. $R$ preserves the colimit of $\gamma$. 

Thus we are left to show that $d \in \on{Fun}^R(A^{\on{op}},D)$. 
Given a diagram $a_i \in A^{\on{op}}$ with limit $a$, we have 
\[
U(d(a)) = U \circ \on{ev}_a (d) = \on{ev}_a(x) = x(a) = x(\on{lim}_i a_i) = \on{lim}_i x(a_i) = \on{lim}_i \on{ev}_{a_i}(x)
\]
\[
= \on{lim}_i U \circ \on{ev}_{a_i} (d) =  \on{lim}_i U(d(a_i))
=U (\on{lim}_i d(a_i))
\]
where we use that $x \in \on{Fun}^R(A^{\on{op}}, C)$ and (\ref{xmnfs}). Since $U$ is conservative, we have $d(a)=\on{lim}_i d(a_i)$, proving $d \in \on{Fun}^R(A^{\on{op}},D)$. 
}
where the right adjoint is the forgetful functor.

\begin{lemma} \label{areyurwo}
There exists right adjoint forgetful functors
\[
\on{DAlg}\on{Cond}\on{Fil}D_{} (\mathbb{Z}) 
\longrightarrow 
\on{CAlg}\on{Cond}\on{Fil}D_{} (\mathbb{Z}),
\]
\[
\on{DAlg}\on{Cond}\on{Gr}D_{} (\mathbb{Z}) 
\longrightarrow 
\on{CAlg}\on{Cond}\on{Gr}D_{} (\mathbb{Z}),
\]
that preserve both limits and colimits.
\end{lemma}

\begin{proof}
Both sides are identified with $\on{Cond}\on{Ani}$ tensored with $\on{DAlg}\on{Fil}D_{} (\mathbb{Z})$ (resp. $\on{DAlg}\on{Gr}D_{} (\mathbb{Z})$) and $\on{CAlg}\on{Fil}D_{} (\mathbb{Z})$ (resp. $\on{CAlg}\on{Gr}D_{} (\mathbb{Z})$) respectively, where the left hand sides are by definition, and the right hand sides hold by \cite[Proposition 7.7]{Gepner2015}. Indeed, (7.8) in $loc.cit.$ proves that for any Cartesian closed symmetric monoidal presentable category $\mathcal{C}$ and symmetric monoidal presentable category $D$, we have $\on{CAlg}(\mathcal{C} \otimes D) \cong \mathcal{C} \otimes \on{CAlg}(D)$.

The functors in the lemma are $\Theta$'s tensored with $\on{Cond}\on{Ani}$, since $\Theta$'s preserve both colimits and limits, they have both left and right adjoints, i.e. they are internal right adjoints in $\on{Pr}^L$,  and tensor product with $\on{Cond}\on{Ani}$, being an $(\infty,2)$-functor from $\on{Pr}^L$ to itself, preserves internal right adjoints. 
\end{proof}

Let $A \in \on{DAlg}\on{Cond}\on{Fil}D_{} (\mathbb{Z})$, then $A$ has a commutative algebra structure in $\on{Cond}\on{Fil}D_{} (\mathbb{Z})$, so we have the module category  $\on{Mod}_{A} \on{Cond}\on{Fil}D_{} (\mathbb{Z})$. Similarly, for $A'\in \on{DAlg}\on{Cond}\on{Gr}D_{} (\mathbb{Z})$,
we have the module category 
$\on{Mod}_{A'} \on{Cond}\on{Gr}D_{} (\mathbb{Z})$. 

\begin{lemma} \label{mxcnaaadfs}
The forgetful functors
\[
\on{Mod}_{A} \on{Cond}\on{Fil}D_{} (\mathbb{Z})
\longleftarrow
\on{DAlg}\on{Cond}\on{Fil}D_{} (\mathbb{Z})_{A/}:U
\]
\[
\on{Mod}_{A'} \on{Cond}\on{Gr}D_{} (\mathbb{Z})
\longleftarrow
\on{DAlg}\on{Cond}\on{Gr}D_{} (\mathbb{Z})_{A'/}: U
\]
from the coslice categories
to 
module categories
are monadic right adjoint functors, so we have monadic adjunctions
\[
\on{Sym}_A^{\bullet}:
\on{Mod}_{A} \on{Cond}\on{Fil}D_{} (\mathbb{Z}) 
\leftrightarrows
\on{DAlg}\on{Cond}\on{Fil}D_{} (\mathbb{Z})_{A/}: U
\]
\[
\on{Sym}_{A'}^{\bullet}:
\on{Mod}_{A'} \on{Cond}\on{Gr}D_{} (\mathbb{Z}) 
\leftrightarrows
\on{DAlg}\on{Cond}\on{Gr}D_{} (\mathbb{Z})_{A'/}: U.
\]
Moreover, we have  commutative diagrams
\[
\begin{tikzcd}
 \on{Mod}_{A} \on{Cond}\on{Fil}D_{} (\mathbb{Z}) 
 \arrow[r,"\on{Sym}_A^{\bullet}"]
 \arrow[d,"\on{gr}^{\bullet}"]
 & 
 \on{DAlg}\on{Cond}\on{Fil}D_{} (\mathbb{Z})_{A/}
 \arrow[d,"\on{gr}^{\bullet}"]
    \\
 \on{Mod}_{\on{gr}^{\bullet}A} \on{Cond}\on{Gr}D_{} (\mathbb{Z}) 
 \arrow[r,"\on{Sym}_{\on{gr}^{\bullet}A}^{\bullet}"]
 &
\on{DAlg}\on{Cond}\on{Gr}D_{} (\mathbb{Z})_{\on{gr}^{\bullet}A/}
\end{tikzcd}
\]
\[
\begin{tikzcd}
 \on{Mod}_{A} \on{Cond}\on{Fil}D_{} (\mathbb{Z}) 
 \arrow[d,"\on{gr}^{\bullet}"]
 & 
 \on{DAlg}\on{Cond}\on{Fil}D_{} (\mathbb{Z})_{A/}
 \arrow[d,"\on{gr}^{\bullet}"]
 \arrow[l,"U"]
    \\
 \on{Mod}_{\on{gr}^{\bullet}A} \on{Cond}\on{Gr}D_{} (\mathbb{Z}) 
 &
\on{DAlg}\on{Cond}\on{Gr}D_{} (\mathbb{Z})_{\on{gr}^{\bullet}A/}.
\arrow[l,"U"]
\end{tikzcd}
\]
\end{lemma}

\begin{proof}
Both functors factorize through
\[
\on{DAlg}\on{Cond}\on{Fil}D_{} (\mathbb{Z})_{A/} 
\longrightarrow 
\on{CAlg}\on{Cond}\on{Fil}D_{} (\mathbb{Z})_{A/}
\longrightarrow 
\on{Mod}_{A} \on{Cond}\on{Fil}D_{} (\mathbb{Z})
\]
\[
\on{DAlg}\on{Cond}\on{Gr}D_{} (\mathbb{Z})_{A/} 
\longrightarrow 
\on{CAlg}\on{Cond}\on{Gr}D_{} (\mathbb{Z})_{A/}
\longrightarrow 
\on{Mod}_{A} \on{Cond}\on{Gr}D_{} (\mathbb{Z}),
\]
where the right functors are the forgetful functors sending commutative algebras  to its underlying modules under the equivalence \[
\on{CAlg}\on{Cond}\on{Fil}D_{} (\mathbb{Z})_{A/} \cong \on{CAlg} \on{Mod}_{A} \on{Cond}\on{Fil}D_{} (\mathbb{Z}) 
\]
\[
\on{CAlg}\on{Cond}\on{Gr}D_{} (\mathbb{Z})_{A/} \cong \on{CAlg} \on{Mod}_{A} \on{Cond}\on{Gr}D_{} (\mathbb{Z}) 
\]
induced by \cite[Corollary 3.4.1.7]{Higheralgebra}, which  preserve sifted colimits by \cite[Corollary 3.2.3.2]{Higheralgebra}. The left functors preserve colimits since we have seen it is so before passing to coslice categories, so the composition functors preserve sifted colimits. The composites are conservative  since their composition with the  forgetful functors to $\on{Cond}\on{Fil}D_{} (\mathbb{Z})$ and $\on{Cond}\on{Gr}D_{} (\mathbb{Z})$
are conservative, thus they are monadic by Barr-Beck-Lurie. 

To prove the commutative diagrams, we observe that when $A=\mathbb{Z}$, this is \cite[remark 4.2.25.]{Raksit2026} and passing to condensed objects. It then induces the commutative diagrams we want
by passing to coslice/module category. 
\end{proof}

\begin{remark}
    There is a non-negatively filtered/graded version of condensed derived rings, defined to be 
    \[
\on{DAlg}\on{Cond}\on{Fil}_{\geq 0}D_{} (\mathbb{Z}) 
:= 
\on{Cond} \on{DAlg}\on{Fil}_{\geq 0}D_{} (\mathbb{Z}) 
\cong 
\on{Cond} \on{Ani} \otimes \on{DAlg} \on{Fil}_{\geq 0} D(\mathbb{Z}),
\]
\[
\on{DAlg}\on{Cond}\on{Gr}_{\geq 0}D_{} (\mathbb{Z}) 
:= 
\on{Cond} \on{DAlg}\on{Gr}_{\geq 0}D_{} (\mathbb{Z}) 
\cong 
\on{Cond} \on{Ani} \otimes \on{DAlg} \on{Gr}_{\geq 0} D(\mathbb{Z}),
\]
  and the claims about filtered/graded rings generalizes directly in this setting. In particular, we have monadic adjunctions  
  \[
\on{Sym}_A^{\bullet}:
\on{Mod}_{A} \on{Cond}\on{Fil}_{\geq 0}D_{} (\mathbb{Z}) 
\leftrightarrows
\on{DAlg}\on{Cond}\on{Fil}_{\geq 0}D_{} (\mathbb{Z})_{A/}: U
\]
\[
\on{Sym}_A^{\bullet}:
\on{Mod}_{A} \on{Cond}\on{Gr}_{\geq 0}D_{} (\mathbb{Z}) 
\leftrightarrows
\on{DAlg}\on{Cond}\on{Gr}_{\geq 0}D_{} (\mathbb{Z})_{A/}: U.
\]
for $A \in \on{DAlg}\on{Cond}\on{Fil}_{\geq 0}D_{} (\mathbb{Z})$ or $\on{DAlg}\on{Cond}\on{Gr}_{\geq 0}D_{} (\mathbb{Z})$. 
\end{remark}

There is a grading functor from filtered to graded objects, and it is compatible with the derived algebra structure, i.e. there exists a (left adjoint) functor 
\[
\on{gr}^{\bullet}: \on{DAlg} \on{Fil} D(\mathbb{Z}) \longrightarrow 
\on{DAlg} \on{Gr} D(\mathbb{Z}) 
\]
that refines the graded functor from $\on{Fil} D(\mathbb{Z}) $ to $\on{Gr} D(\mathbb{Z}) $, see \cite[Remark 4.3.6]{Raksit2026}. Tensoring with $\on{Cond}\on{Ani}$, we obtain a left adjoint functor
\[
\on{gr}^{\bullet}: \on{DAlg}\on{Cond}\on{Fil}D_{} (\mathbb{Z}) \longrightarrow 
\on{DAlg}\on{Cond}\on{Gr}D_{} (\mathbb{Z}). 
\]
Moreover, we have from $loc.cit.$ a right adjoint functor 
\[
\on{ev}^0: \on{DAlg} \on{Gr} D(\mathbb{Z}) \longrightarrow 
\on{DAlg}  D(\mathbb{Z}) 
\]
which is the evaluation at degree 0. Tensoring with $\on{Cond}\on{Ani}$, we have a right adjoint
\[
\on{ev}^{0}: \on{DAlg}\on{Cond}\on{Gr}D_{} (\mathbb{Z}) \longrightarrow 
\on{DAlg}\on{Cond}D_{} (\mathbb{Z}), 
\]
and we define 
\[
\on{gr}^0: \on{DAlg}\on{Cond}\on{Fil}D_{} (\mathbb{Z}) \longrightarrow 
\on{DAlg}\on{Cond}D_{} (\mathbb{Z})
\]
to be the composition of $\on{ev}^0$
with $gr^{\bullet}$. 

Note that $\on{gr}^0$ is a composition of a right adjoint with a left adjoint, hence does not necessarily preserve either limits or colimits. The situation is better if we restrict to the case of non-negatively filtered objects. The key difference is that $\on{ev}^0 : \on{Gr}_{\geq 0} D(\mathbb{Z}) \longrightarrow D(\mathbb{Z})$ is symmetric monoidal, which is easily seen defining a morphism of derived algebra context in the sense of \cite[Definition 4.2.1]{Raksit2026}, giving rise to a left adjoint 
\[
\on{ev}^{0}: \on{DAlg}\on{Cond}\on{Gr}_{\geq 0}D_{} (\mathbb{Z}) \longrightarrow 
\on{DAlg}\on{Cond}D_{} (\mathbb{Z}).
\] 
Then we have a left adjoint functor
\[
\on{gr}^0: \on{DAlg}\on{Cond}\on{Fil}_{\geq0}D_{} (\mathbb{Z}) \longrightarrow 
\on{DAlg}\on{Cond}D_{} (\mathbb{Z}) 
\]
given by the composition of $\on{ev}^0$ and $\on{gr}^{\bullet}$. 

The derived algebra structure has a connective cover, i.e. we have a right adjoint functor
\[
\tau_{\geq 0}:
\on{DAlg}D_{} (\mathbb{Z}) 
\longrightarrow \on{DAlg}D_{\geq 0} (\mathbb{Z})
 \]
 that refines the $\tau_{\geq 0}$ on the underlying filtered category, see \cite[Proposition 2.16]{BENBASSAT2025110530}. We know from \cite[Remark 4.3.3.]{Raksit2026} that the right handside is the category of animated rings.  Tensoring with $\on{Cond}\on{Ani}$ gives us right adjoint functor
\[
\tau_{\geq 0}:
\on{DAlg}\on{Cond}D_{} (\mathbb{Z}) 
\longrightarrow \on{DAlg}\on{Cond}D_{\geq 0} (\mathbb{Z}) = \on{Cond} \on{AniRing}. 
 \]

We also need to recall another standard notion of completely filtered objects. 
We have a full subcategory of completely filtered objects, namely
 \[
 \on{Fil}^{\wedge}D_{} (\mathbb{Z}) \subset \on{Fil}D_{} (\mathbb{Z})
 \]
 is the full subcategory such that $X^{\bullet} \in \on{Fil}^{\wedge}D_{} (\mathbb{Z})$
 if and only if 
 $\lim_{i \rightarrow \infty} X^i = 0$. It comes with a left adjoint whose kernel is a tensor ideal, and it equips $\on{Fil}^{\wedge}D_{} (\mathbb{Z})$ 
 with a symmetric monoidal structure such that the left adjoint is symmetric monoidal. Moreover, the graded functor
 \[
  \on{Fil}^{\wedge}D_{} (\mathbb{Z})
  \longrightarrow 
   \on{Gr}D_{} (\mathbb{Z})
 \]
 is conservative, and by \cite[theorem 3.2.14]{Raksit2026}, it induces an equivalence
 \begin{equation} \label{nxcsfd}
  \on{Fil}^{\wedge}D_{} (\mathbb{Z})
  \cong 
  \on{Mod}_{\mathbb{D}_-} \on{Gr}D_{} (\mathbb{Z})
 \end{equation}
 of symmetric monoidal categories
 for some cocommutative coalgebra object of derived algebra in graded categories 
 \[
 \mathbb{D}_- = \mathbb{Z}[\epsilon
 ]/\epsilon^2  \in \on{cCAlg}( \on{DAlg} (\on{Gr}D_{} (\mathbb{Z}))) \footnote{To see that it has a derived algebra structure, we know from \cite[proposition 5.1.7.]{Raksit2026} that $\mathbb{D}_+^{\vee} \in \on{cCAlg}( \on{DAlg} (\on{Gr}D_{} (\mathbb{Z})))$, we simply observe that $\mathbb{D}_- = I( \mathbb{D}_+^{\vee})$, where $I: \on{Gr}D_{} (\mathbb{Z})) \longrightarrow \on{Gr}D_{} (\mathbb{Z}))$ 
 is the functor that reverse the grading, i.e. $I(X^{\bullet}) = X^{-\bullet}$, which is easily seen to be a equivalence of derived algebra contexts.  It follows that the identification $\mathbb{D}_- = I( \mathbb{D}_+^{\vee})$ induces the coalgebra of a derived algebra structure on $\mathbb{D}_-$. }
 \]
 as defined in \cite[notation 3.2.11, construction 5.1.1.]{Raksit2026}, 
 where $\epsilon$ is of  graded degree $1$ and homotopy degree $-1$. 
 We give the symmetric monoidal structure on $\on{Mod}_{\mathbb{D}_-} \on{Gr}D_{} (\mathbb{Z})$ 
as in \cite[variant 2.2.3.]{Raksit2026} so that the tensor product is computed in
$\on{Gr}D_{} (\mathbb{Z})$, in other words, the forgetful functor 
$\on{Mod}_{\mathbb{D}_-} \on{Gr}D_{} (\mathbb{Z}) \longrightarrow \on{Gr}D_{} (\mathbb{Z})$
is symmetric monoidal. Note that this is not the relative tensor product over $\mathbb{D}_-$.

 We now define 
 \[
 \on{Cond} \on{Fil}^{\wedge}D_{} (\mathbb{Z})
 \cong 
 \on{Cond}\on{Ani} \otimes
  \on{Fil}^{\wedge}D_{} (\mathbb{Z})
 \]
which is then a full subcategory of $ \on{Cond} \on{Fil}D_{} (\mathbb{Z})$
with  symmetric monoidal left adjoint
\[
(-)^{\wedge}: 
\on{Cond} \on{Fil}D_{} (\mathbb{Z}) 
\longrightarrow
\on{Cond} \on{Fil}^{\wedge}D_{} (\mathbb{Z}).
\]
In particular, $\on{Cond} \on{Fil}^{\wedge}D_{} (\mathbb{Z})$ is tensored over $\on{Cond} \on{Fil}D_{} (\mathbb{Z})$ 
through the left adjoint, so $\on{Mod}_A \on{Cond} \on{Fil}^{\wedge}D_{} (\mathbb{Z})$
is defined for $A \in \on{DAlg}\on{Cond} \on{Fil}D_{} (\mathbb{Z})$, and we have an adjunction
\[
(-)^{\wedge}:
\on{Mod}_A \on{Cond} \on{Fil}D_{} (\mathbb{Z}) 
\leftrightarrows 
\on{Mod}_A \on{Cond} \on{Fil}^{\wedge}D_{} (\mathbb{Z}) : U
\]
which is induced by the adjunction pair on the underlying category.  Moreover, (\ref{nxcsfd}) induces an equivalence
\begin{equation} \label{xncvuyedi}
\on{gr}^{\bullet}:  \on{Mod}_A\on{Cond} \on{Fil}^{\wedge}D_{} (\mathbb{Z})
  \cong 
  \on{Mod}_{\mathbb{D}_-} 
  \on{Mod}_{\on{gr}^{\bullet} A}\on{Gr}\on{Cond}D_{} (\mathbb{Z})
 \end{equation}
 of symmetric monoidal categories, where again the symmetric monoidal structure on the right hand side is given by  \cite[variant 2.2.3.]{Raksit2026} such that forgetting the $\mathbb{D}_-$-module structure is a symmetric monoidal functor.

\begin{lemma} \label{xqrupdn}
    The monad $U \circ \on{Sym}_A^{\bullet}$ on $\on{Mod}_A \on{Cond} \on{Fil}D_{} (\mathbb{Z})$ 
    in lemma \ref{mxcnaaadfs} induces a monad 
    \begin{equation} \label{wemvcnw}
    T_A := (-)^{\wedge} \circ U \circ \on{Sym}_A^{\bullet}
    \end{equation}
    on 
    $\on{Mod}_A \on{Cond} \on{Fil}^{\wedge}D_{} (\mathbb{Z})$, and we denote the corresponding algebra by 
    \[
    \on{DAlg}\on{Cond}\on{Fil}^{\wedge}D_{} (\mathbb{Z}) := \on{LMod}_{T_{\mathbb{Z}}} \on{Cond}\on{Fil}^{\wedge}D_{} (\mathbb{Z})
    \]
    when $A=\mathbb{Z}$. 
\end{lemma}

\begin{proof}
    We need to prove that for $M \longrightarrow N$ which becomes an isomorphism after applying $(-)^{\wedge}$, then  $T(M) \longrightarrow T(N)$ is also an isomorphism. We know $M\longrightarrow N$ becomes an isomorphism after applying $(-)^{\wedge}$ if and only if the associated graded is an equivalence. Applying the same criterion to $T= (-)^{\wedge} \circ U \circ \on{Sym}_A^{\bullet}$, we see that it is enough to prove that 
    $\on{gr}^{\bullet} \circ U \circ \on{Sym}_A^{\bullet} $
    preserves maps which induces an isomorphism on graded pieces, but this follows immediately from the computation 
     $\on{gr}^{\bullet} \circ U \circ \on{Sym}_A^{\bullet}= 
    U \circ \on{Sym}_{\on{gr}^{\bullet}A}^{\bullet} \circ  \on{gr}^{\bullet} $,
   which is the commutative diagrams in lemma \ref{mxcnaaadfs}. 
\end{proof}

\begin{remark}
    We see easily that lemma \ref{mxcnaaadfs} can be generalized to the completely filtered setting, i.e. the forgetful functors
\[
\on{Mod}_{A^{\wedge}} \on{Cond}\on{Fil}^{\wedge}D_{} (\mathbb{Z})
\longleftarrow
\on{DAlg}\on{Cond}\on{Fil}^{\wedge}D_{} (\mathbb{Z})_{A^{\wedge}/}:U
\]
from the coslice categories
to 
module categories
are monadic right adjoint functors, so we have monadic adjunctions
\[
\on{Sym}_{A^{\wedge}}^{\bullet}:
\on{Mod}_{A^{\wedge}} \on{Cond}\on{Fil}^{\wedge}D_{} (\mathbb{Z}) 
\leftrightarrows
\on{DAlg}\on{Cond}\on{Fil}^{\wedge}D_{} (\mathbb{Z})_{A^{\wedge}/}: U.
\]
Moreover, $\on{Sym}_{A^{\wedge}}^{\bullet} = (-)^{\wedge} \circ \on{Sym}_{A}^{\bullet}$.  We have a natural equivalence
    \[
    \on{DAlg}\on{Cond}\on{Fil}^{\wedge}D_{} (\mathbb{Z})_{A^{\wedge}/} \cong \on{LMod}_{T_{A}} 
    \on{Mod}_A
    \on{Cond}\on{Fil}^{\wedge}D_{} (\mathbb{Z}).
    \]
\end{remark}

\begin{definition}
    A filtered uncompleted analytic derived ring is a pair $(A, D(A))$, where $A \in \on{DAlg}\on{Cond}\on{Fil}D_{} (\mathbb{Z}) $, and 
    \[
  U_A:  D(A) \subset \on{Mod}_{A} \on{Cond}\on{Fil}^{\wedge}D_{} (\mathbb{Z})
    \]
    is a full subcategory such that 

    (1) $D(A)$ is closed under limits and colimits.
    
    (2) The kernel of the left adjoint $L_A$ to $U_A$ is  a tensor ideal, providing  $D(A)$ with a canonical symmetric monoidal structure such that $L_A$
    is symmetric monoidal. In particular, $D(A)$ is tensored over $ \on{Mod}_{A} \on{Cond}\on{Fil}^{\wedge}D_{} (\mathbb{Z})$.

    (3) $D(A)$ is closed under the shift of degrees functor $[1]$.

    (4) The monad $T_A$  in (\ref{wemvcnw})  preserves the class of maps which becomes isomorphism after applying $L_A$, so that $
    \hat{T}_A :=
    L_A \circ T_A \circ  U_A$ is a monad on $D(A)$ such that $\hat{T}_A \circ L_A = L_A \circ T_A$. 
\end{definition}

\begin{definition}
    A filtered analytic derived ring is a filtered uncompleted analytic derived ring $(A, D(A))$ such that $A\in D(A)$. 
\end{definition}

\begin{remark}
    Compared to the definition of analytic rings, we do not require that $L_A$ preserves connective objects, this is because that we are in the possibly nonconnective derived rings situation, and the t-structure on $\mathbb{Z}$-modules does not induce a t-structure on $A$-modules. It is possible to define a t-structure on $A$-modules, but there are multiple choices, and the t-structure can be degenerate  in general, i.e. either the connective or coconnective subcategory is zero. 
\end{remark}

\begin{remark}
    We do not know if condition (4) is always satisfied as in the case of analytic rings, see \cite[proposition 4.2.6.]{Notesonsolidgeometry}. Note that this is one of the very nice features of the light condensed setting. In the classical setting of condensed math where one uses all profinite sets, we need to put a condition on Frobenius to make sure that (4) holds. 
\end{remark}

We also need the definition of graded analytic derived rings, which is a simple variant of the filtered version above.

\begin{definition}
    A graded uncompleted analytic derived ring is a pair $(A, D(A))$, where $A \in \on{DAlg}\on{Cond}\on{Gr}D_{} (\mathbb{Z}) $, and 
    \[
  U_A:  D(A) \subset \on{Mod}_{A} \on{Cond}\on{Gr}D_{} (\mathbb{Z})
    \]
    is a full subcategory such that 

    (1) $D(A)$ is closed under limits and colimits.
    
    (2) The kernel of the left adjoint $L_A$ to $U_A$ is  a tensor ideal, providing  $D(A)$ with a canonical symmetric monoidal structure such that $L_A$
    is symmetric monoidal. In particular, $D(A)$ is tensored over $ \on{Mod}_{A} \on{Cond}\on{Gr}D_{} (\mathbb{Z})$.

    (3) $D(A)$ is closed under the shift of degrees functor $[1]$.

    (4) The monad $T_A := U_A \circ \on{Sym}_A $ in lemma \ref{mxcnaaadfs}  preserves the class of maps which becomes isomorphism after applying $L_A$, so that $
    \hat{T}_A :=
    L_A \circ T_A \circ  U_A$ is a monad on $D(A)$ such that $\hat{T}_A \circ L_A = L_A \circ T_A$. 
\end{definition}

\begin{definition}
    A graded analytic derived ring is a filtered uncompleted analytic derived ring $(A, D(A))$ such that $A\in D(A)$. 
\end{definition}

We first prove that the filtered uncompleted analytic ring structure is uniquely characterized by  the corresponding subcategory of the graded category. 

 Let $(A, D(A))$ be a filtered uncompleted analytic derived ring, and we have the corresponding graded category
    \[
 \on{Mod}_{C\tau} D(A), 
\]
where 
\[
C\tau := \on{cofib}(\mathbb{S}[1] \overset{\tau}{\longrightarrow} \mathbb{S}) \in \on{Fil} \on{Sp}
\]
is a commutative algebra object in  $\on{Fil} \on{Sp}$ such that $\on{Gr} \on{Sp} \cong \on{Mod}_{C\tau}  \on{Fil} \on{Sp}$ as in \cite[section 3.2.2.]{Lecturesynthetic}, $\mathbb{S} \in \on{Fil} \on{Sp}$ is the unit object, $\tau$ is defined in 
\cite[definition 3.1.3.]{Lecturesynthetic} that corresponds to transition maps of filtered objects, and we use that $D(A)$ is tensored over $\on{Fil} \on{Sp}$.  
More precisely, the  $\on{Fil} \on{Sp}$-linear structure of $D(A)$
is given by the composite of symmetric monoidal left adjoint functors
\[
 \on{Fil} \on{Sp} \longrightarrow 
 \on{Cond}\on{Fil}D_{} (\mathbb{Z}) 
 \overset{A\otimes -}{\longrightarrow}
 \on{Mod}_{A} \on{Cond}\on{Fil}D_{} (\mathbb{Z}) 
 \longrightarrow
 \on{Mod}_{A} \on{Cond}\on{Fil}^{\wedge}D_{} (\mathbb{Z})
 \overset{L_1}{\longrightarrow}
 D(A)
\]
where the first functor is the $\on{Fil}$ of the initial symmetric monoidal functor from $\on{Sp}$ to the stable symmetric monoidal category $\on{Cond}D_{} (\mathbb{Z})$, the third functor is the left adjoint to the inclusion of completely filtered to filtered objects, and the last functor is symmetric monoidal thanks to condition (2) of $D(A)$. 

\begin{proposition}
    Let $(A, D(A))$ be a filtered uncompleted analytic derived ring, then $(\on{gr}^{\bullet} A,  \on{Mod}_{C\tau} D(A)  )$ 
    is a graded analytic uncompleted analytic derived ring. 
    
    More precisely,
 there exists a fully faithful embedding 
    \[
    U_{\on{gr}^{\bullet}A}:
 \on{Mod}_{C\tau} D(A) \hookrightarrow 
 \on{Mod}_{\on{gr}^{\bullet}A} \on{Cond}\on{Gr}D_{} (\mathbb{Z}).
\]
such that

  (1)  $U_{\on{gr}^{\bullet}A}$ preserves both limits and colimits.

(2) $\on{Mod}_{C\tau} D(A)$ is equipped with a symmetric monoidal structure such that left adjoint $L_{\on{gr}^{\bullet}A}$ of $U_{\on{gr}^{\bullet}A}$ is symmetric monoidal.

(3) $\on{Mod}_{C\tau} D(A)$ is closed under shifts of degree 1 functor $[1]$. 

(4) 
The class of maps which becomes isomorphism after applying $L_{\on{gr}^{\bullet}A}$
is preserved by the monad 
\[
T_{\on{gr}^{\bullet}A}:= U \circ \on{Sym}^{\bullet}_{\on{gr}^{\bullet}A}
\]
of the monadic adjunction 
\[
\on{Sym}_{\on{gr}^{\bullet}A}^{\bullet}:
\on{Mod}_{\on{gr}^{\bullet}A} \on{Cond}\on{Gr}D_{} (\mathbb{Z}) 
\leftrightarrows
\on{DAlg}\on{Cond}\on{Gr}D_{} (\mathbb{Z})_{\on{gr}^{\bullet}A/}: U.
\]
of lemma \ref{mxcnaaadfs}, 
so that $
    \hat{T}_{\on{gr}^{\bullet}A} :=
    L_{\on{gr}^{\bullet}A} \circ T_{\on{gr}^{\bullet}A} \circ  U_{\on{gr}^{\bullet}A}$ is a monad on $\on{Mod}_{C\tau} D(A)$ such that $\hat{T}_{\on{gr}^{\bullet}A} \circ L_{\on{gr}^{\bullet}A} = L_{\on{gr}^{\bullet}A} \circ T_{\on{gr}^{\bullet}A}$. 

    Moreover, we have that 
    \begin{equation} \label{808098080}
    D(A) = \on{Mod}_{\mathbb{D}_-} \on{Mod}_{C\tau} D(A). 
    \end{equation}
\end{proposition}

\begin{proof}
    Given  $(A, D(A))$ a filtered uncompleted analytic derived ring,  let 
    \[
G: \on{Mod}_{\on{gr}^{\bullet}A} \on{Cond}\on{Gr}D_{} (\mathbb{Z})
\longrightarrow 
\on{Mod}_{A} \on{Cond}\on{Fil}^{\wedge}D_{} (\mathbb{Z}) 
\]
be the right adjoint of $\on{gr}^{\bullet}$. 
Then we define a subcategory $C $ of  $\on{Mod}_{\on{gr}^{\bullet}A} \on{Cond}\on{Gr}D_{} (\mathbb{Z}) $
by the pullback
\begin{equation} \label{werererewdfd}
\begin{tikzcd}
\on{Mod}_{\on{gr}^{\bullet}A} \on{Cond}\on{Gr}D_{} (\mathbb{Z}) 
\arrow[d,"G"]
\arrow[r,"U_{\on{gr}^{\bullet}A}",hookleftarrow] 
&
C
\arrow[d,"G"]
\arrow[dl, phantom, "\urcorner", very near start]
\\
\on{Mod}_{A} \on{Cond}\on{Fil}^{\wedge}D_{} (\mathbb{Z})
\arrow[r,hookleftarrow,"U_A"]
&
D(A).
&
\end{tikzcd}
\end{equation}
We observe that $C$ can be identified as
\begin{equation}  \label{5765586}
C \cong \on{Mod}_{C\tau} D(A) 
\end{equation}

To prove (\ref{5765586}), we first observe that 
\[
\on{Mod}_{\on{gr}^{\bullet}A} \on{Cond}\on{Gr}D_{} (\mathbb{Z})  
\cong 
\on{Mod}_{C\tau} \on{Mod}_{A} \on{Cond}\on{Fil}^{\wedge}D_{} (\mathbb{Z})
\]
by \cite[proposition C.2.]{Burklund2026Galois}, and the functor $G$ corresponds to forgetting the $C\tau$-module structure under this equivalence, since  $\on{gr}^{\bullet}$ corresponds to the free $C\tau$-module functor $\otimes C\tau$.
We then have a commutative diagram
\begin{equation} \label{weyruewir}
\begin{tikzcd}
\on{Mod}_{C\tau} \on{Mod}_{A} \on{Cond}\on{Fil}^{\wedge}D_{} (\mathbb{Z})
\arrow[d,"G"]
\arrow[r,"U_{\on{gr}^{\bullet}A}",hookleftarrow] 
&
\on{Mod}_{C\tau} D(A)
\arrow[d,"G"]
\\
\on{Mod}_{A} \on{Cond}\on{Fil}^{\wedge}D_{} (\mathbb{Z})
\arrow[r,hookleftarrow,"U_{A}"]
&
D(A)
\end{tikzcd}
\end{equation}
which is a pullback diagram by \cite[proposition 4.1.9.]{Raksit2026} applied to the monad of tensoring with $C\tau$ on
$ \on{Mod}_{A} \on{Cond}\on{Fil}^{\wedge}D_{} (\mathbb{Z})$. 

\begin{lemma} \label{797384378}
    The pullback diagram (\ref{weyruewir}) is left adjointable, i.e. we have a commutative diagram
    \[
    \begin{tikzcd}
\on{Mod}_{C\tau} \on{Mod}_{A} \on{Cond}\on{Fil}^{\wedge}D_{} (\mathbb{Z})
\arrow[d,"G"]
\arrow[r,"L_{\on{gr}^{\bullet}A}"] 
&
\on{Mod}_{C\tau} D(A)
\arrow[d,"G"]
\\
\on{Mod}_{A} \on{Cond}\on{Fil}^{\wedge}D_{} (\mathbb{Z})
\arrow[r,"L_{A}"]
&
D(A)
\end{tikzcd}
    \]
\end{lemma}

\begin{proof}
    By lemma \ref{urowewqpepwp}, it is enough to prove that the right adjoint 
    \[
    R: \on{Mod}_{A} \on{Cond}\on{Fil}^{\wedge}D_{} (\mathbb{Z})
    \longrightarrow 
    \on{Mod}_{C\tau} \on{Mod}_{A} \on{Cond}\on{Fil}^{\wedge}D_{} (\mathbb{Z})
    \]
    of $G$ has the property that 
    $R \circ U_A$ has image in the full subcategory $\on{Mod}_{C\tau} D(A) $, and this is equivalent to $G\circ R $ preserves the subcategory $D(A)$. We can compute directly that $R = \Sigma^{-1}[1] \circ \on{gr}^{\bullet}$, so $G\circ R = \Sigma^{-1}[1] \circ G \circ \on{gr}^{\bullet} $, but 
    $G \circ \on{gr}^{\bullet} (X^{\bullet}) = \on{cofib} (X^{\bullet}[1] \overset{\tau}{ \longrightarrow} X^{\bullet})$ 
    preserves the subcategory $D(A)$ by axiom (1) and (3), so does $G\circ R$ as well.
\end{proof}

(1)   
We observe that $G$ preserves both limits and colimits as it sends graded objects to filtered objects with 0 transition maps, so $G$ has right adjoint as well. Since $U_{A}$ has both left and right adjoint by definition of filtered analytic derived rings, so does  $U_{\on{gr}^{\bullet}A}$ as the pullback can be taken both in $\on{Pr}^L$ and $\on{Pr}^R$. 

(2) We observe that (\ref{5765586}) equips $C$ with  a symmetric monoidal structure and the left adjoint $L_{\on{gr}^{\bullet}A}$ of $U_{\on{gr}^{\bullet}A}$ is symmetric monoidal since it is induced by $L_A$, which is symmetric monoidal as $D(A)$ is a filtered analytic derived ring.

(3) $C$ is stable under $[1]$ since $D(A)$ is so. 

(4) Since $\on{gr}^{\bullet}$ is a morphism of derived algebra context as in \cite[remark 4.3.6.]{Raksit2026}, passing to condensed objects we have that 
$T_{\on{gr}^{\bullet}A} \circ \on{gr}^{\bullet} =\on{gr}^{\bullet} \circ T_A$. Taking left adjoint of the commutative diagram (\ref{weyruewir}), we have that 
$\on{gr}^{\bullet} \circ L_{A} = L_{\on{gr}^{\bullet} A} \circ \on{gr}^{\bullet}$, 
and similarly taking the right adjoint of the commutative diagram in lemma \ref{797384378}, we have 
$\on{gr}^{\bullet} \circ U_{A} = U_{\on{gr}^{\bullet} A} \circ \on{gr}^{\bullet}$.
It follows that 
$T_{\on{gr}^{\bullet}A} \circ  U_{\on{gr}^{\bullet} A} \circ
L_{\on{gr}^{\bullet} A} = U_{\on{gr}^{\bullet} A} \circ L_{\on{gr}^{\bullet} A} \circ T_{\on{gr}^{\bullet}A}$, which is $\on{gr}^{\bullet}$ of the formula 
$T_A \circ U_A \circ L_A = U_A \circ L_A \circ T_A$ that holds thanks to the last axiom of filtered uncompleted analytic derived rings. 

To prove the last claim (\ref{808098080}), 
we apply \cite[proposition 4.1.9.]{Raksit2026}  to the subcategory $C$ of $\on{Mod}_{\on{gr}^{\bullet}A} \on{Cond}\on{Gr}D_{} (\mathbb{Z}) $ and the monad of tensoring with $\mathbb{D}_-$ on $\on{Mod}_{\on{gr}^{\bullet}A} \on{Cond}\on{Gr}D_{} (\mathbb{Z}) $ to obtain the cartesian diagram 
\[
\begin{tikzcd}
\on{Mod}_{\mathbb{D}_-}
\on{Mod}_{\on{gr}^{\bullet}A} \on{Cond}\on{Gr}D_{} (\mathbb{Z})
\arrow[r,hookleftarrow,"U_A"]
\arrow[d]
& 
\on{Mod}_{\mathbb{D}_-} C
\arrow[d]
\arrow[dl, phantom, "\urcorner", very near start]
\\
\on{Mod}_{\on{gr}^{\bullet}A} \on{Cond}\on{Gr}D_{} (\mathbb{Z}) 
\arrow[r,hookleftarrow,"U_{\on{gr}^{\bullet} A}"]
&
C
\end{tikzcd}
\] 
with vertical arrows forgetful functors. Composing with the cartesian diagram (\ref{werererewdfd}), we have the cartesian diagram
\[
\begin{tikzcd}
\on{Mod}_{\mathbb{D}_-}
\on{Mod}_{\on{gr}^{\bullet}A} \on{Cond}\on{Gr}D_{} (\mathbb{Z})
\arrow[r,hookleftarrow,"U_A"]
\arrow[d]
& 
\on{Mod}_{\mathbb{D}_-} C
\arrow[d]
\arrow[dl, phantom, "\urcorner", very near start]
\\
\on{Mod}_{\on{gr}^{\bullet}A} \on{Cond}\on{Gr}D_{} (\mathbb{Z}) 
\arrow[d,"G"]
\arrow[r,"U_{\on{gr}^{\bullet}A}",hookleftarrow] 
&
C
\arrow[d,"G"]
\arrow[dl, phantom, "\urcorner", very near start]
\\
\on{Mod}_{A} \on{Cond}\on{Fil}^{\wedge}D_{} (\mathbb{Z})
\arrow[r,hookleftarrow,"U_A"]
&
D(A).
&
\end{tikzcd}
\]

Now using the equivalence
 (\ref{xncvuyedi}), and the fact that  the forgetting $\mathbb{D}_-$-module functor corresponds to  $\on{gr}^{\bullet}$ under this equivalence, we see that the composite of left vertical maps is $G\circ \on{gr}^{\bullet} $, which as we have seen in the proof of lemma  \ref{797384378} has the form 
 $ G \circ \on{gr}^{\bullet} (X^{\bullet}) = \on{cofib} (X^{\bullet}[1] \overset{\tau}{ \longrightarrow} X^{\bullet}) $, and we need to prove that 
 $X^{\bullet} \in  D(A)$ if and only if 
 $G \circ \on{gr}^{\bullet} (X^{\bullet}) \in D(A)$, from which we can conclude that $D(A) = \on{Mod}_{\mathbb{D}_-} C$ since the diagram is cartesian. 
 
 Now  assume that $X^{\bullet} \in  D(A)$, then $ G \circ \on{gr}^{\bullet} (X^{\bullet}) = \on{cofib} (X^{\bullet}[1] \overset{\tau}{ \longrightarrow} X^{\bullet}) \in D(A)$
 by axiom (1) and (3) of filtered analytic derived rings. To see the other direction, assume that $\on{cofib} (X^{\bullet}[1] \overset{\tau}{ \longrightarrow} X^{\bullet}) \in D(A)$, then 
 $\on{cofib} (X^{\bullet}[n] \overset{\tau^n}{ \longrightarrow}X^{\bullet}) \in D(A) $
 since $\on{cofib} (X^{\bullet}[n] \overset{\tau^n}{ \longrightarrow}X^{\bullet})  $ 
 is an extension of $\on{cofib} (X^{\bullet}[n] \overset{\tau^{n-1}}{ \longrightarrow}X^{\bullet}) $ 
 by 
 $\on{cofib} (X^{\bullet}[n] \overset{\tau}{\longrightarrow}X^{\bullet}) $, 
 so we can conclude by induction and uses axiom (1) that $D(A)$ is closed under limits and colimits.  Then we use that $X^{\bullet} = \on{lim}_{n}  \on{cofib} (X^{\bullet}[n] \overset{\tau^n}{ \longrightarrow} X^{\bullet}) $, which holds since $X^{\bullet}$ is complete, proving $X^{\bullet} \in D(A)$ since $D(A)$ is closed under limits by axiom (1) again. 
\end{proof}

\begin{lemma} \label{urowewqpepwp}
    Let $U: \mathcal{C} \hookrightarrow \mathcal{D}$ 
    be a full subcategory of a presentable category $\mathcal{D}$  with left adjoint $L: \mathcal{D} \longrightarrow \mathcal{C}$. Let $F: \mathcal{D}' \longrightarrow \mathcal{D} $ be a left adjoint functor between presentable categories with right adjoint $G: \mathcal{D}  \longrightarrow 
    \mathcal{D}'$, we assume that $F$ also has a left adjoint.  Let $\mathcal{C}' \subset \mathcal{D}'$ 
    \[
    \begin{tikzcd}
       \mathcal{D}'
   \arrow[d,"F"]
   \arrow[r,hookleftarrow, "U'"]
   & 
        \mathcal{C}'
        \arrow[d,"F_{\mathcal{C}}"]
        \\ 
     \mathcal{D} 
     \arrow[r,hookleftarrow, "U"]
     & 
        \mathcal{C}
    \end{tikzcd}
    \]
be    the pullback of $\mathcal{C}$ along $F$ with left adjoint $L': \mathcal{D}' \longrightarrow \mathcal{C}'$. Then there exists a commutative diagram
    \[
    \begin{tikzcd}
   \mathcal{D}'
   \arrow[d,"F"]
   \arrow[r,"L'"]
   & 
        \mathcal{C}'
        \arrow[d,"F_{\mathcal{C}}"]
        \\ 
     \mathcal{D} 
     \arrow[r,"L"]
     & 
        \mathcal{C}
    \end{tikzcd}
    \]
    if  the functor $ G \circ U: \mathcal{D} \longrightarrow \mathcal{D}'$ has image in the full subcategory $\mathcal{C}'$. 
\end{lemma}

\begin{proof} 
Since $F$ also has left adjoint, $\mathcal{C}'$ can be viewed as pullback in $\on{Pr}^R$, so $U'$ has left adjoint $L'$ as well. 
    There exists a natural transform
    $L \circ F \Rightarrow F_{\mathcal{C}} \circ L' $, which by adjunction is given by 
    $F \Rightarrow F \circ U'\circ L' =  U \circ  F_{\mathcal{C}} \circ L'$
    with the first natural transform the unit of the adjunction $L'$ and $U' $.

    Let $d' \in \mathcal{D}'$, and $c \in \mathcal{C}$, then 
    \begin{align*}
      \on{Map}(L\circ F (d'), c)
        & =  \on{Map}(d', G\circ U(c))
        \\
         & =  \on{Map}(U'\circ L' (d'),  G\circ U(c))
          \\
  & = \on{Map}(F \circ U' \circ L'(d'), U(c))
  \\
  & = \on{Map}(U\circ F_{\mathcal{C}} \circ L'(d'), U(c))
  \\
 &=  \on{Map}(F_{\mathcal{C}} \circ L'(d'), c) 
    \end{align*}
    where in the second equivalence we use our assumption that $G \circ U $ lands in the full subcategory $\mathcal{C}'$.
\end{proof}

We now prove that the category $\on{Mod}_{C\tau} D(A) $ also characterize $D(A)$. 

\begin{proposition} \label{2378937583}
Let 
$A \in \on{DAlg}\on{Cond}\on{Fil}D_{} (\mathbb{Z})$, and $(\on{gr}^{\bullet}A, C)$ is a graded uncompleted analytic derived ring. In other words, 
\[
U_{\on{gr}^{\bullet}A}: C \hookrightarrow  
\on{Mod}_{\on{gr}^{\bullet}A} \on{Cond}\on{Gr}D_{} (\mathbb{Z})
\]
is
a full subcategory such that 

 (1) $C$ is closed under both limits and colimits.

(2) $C$ is equipped with a symmetric monoidal structure such that left adjoint $L_{\on{gr}^{\bullet}A}$ of $U_{\on{gr}^{\bullet}A}$ is symmetric monoidal.

(3) $C$ is closed under shifts of degree 1 functor $[1]$.

(4)The class of maps which becomes isomorphism after applying $L_{\on{gr}^{\bullet}A}$
is preserved by the monad 
$
T_{\on{gr}^{\bullet}A}$, so that it induces a monad $\hat{T}_{\on{gr}^{\bullet} A}$ on $C$ such that $\hat{T}_{\on{gr}^{\bullet} A} \circ L_{\on{gr}^{\bullet} A} =L_{\on{gr}^{\bullet} A} \circ T_{\on{gr}^{\bullet} A}$. 

Then 
\[
D(A) := \on{Mod}_{\mathbb{D}_-} C 
\]
is a full subcategory of 
$\on{Mod}_{A} \on{Cond}\on{Fil}^{\wedge}D_{} (\mathbb{Z})$ such that $(A, D(A))$ is a filtered uncompleted analytic derived ring, and 
the left adjoint $L_A$ to the inclusion $D(A) \subset \on{Mod}_{A} \on{Cond}\on{Fil}^{\wedge}D_{} (\mathbb{Z})$ satisfies 
    \[
    \on{gr}^{\bullet} \circ L_A = 
    L_{\on{gr}^{\bullet}A}
    \circ \on{gr}^{\bullet}. 
    \]
Moreover, we have that 
\[
\on{Mod}_{C\tau}D(A) = C. 
\]
\end{proposition}

\begin{proof}
Using (\ref{xncvuyedi}) and \cite[proposition 4.1.9.]{Raksit2026} applied to the monad of tensoring with $\mathbb{D}_-$, we have a cartesian diagram
\begin{equation} \label{2463764736723}
\begin{tikzcd}
\on{Mod}_{A} \on{Cond}\on{Fil}^{\wedge}D_{} (\mathbb{Z})
\arrow[d,"\on{gr}^{\bullet}"]
\arrow[r,hookleftarrow,"U_A"]
&
D(A) 
\arrow[d,"\on{gr}^{\bullet}"]
\arrow[dl, phantom, "\urcorner", very near start]
\\
\on{Mod}_{\on{gr}^{\bullet}A} \on{Cond}\on{Gr}D_{} (\mathbb{Z})
\arrow[r,hookleftarrow,"U_{\on{gr}^{\bullet} A}"]
&
C.
\end{tikzcd}
\end{equation}  
The axioms (1) and (3) of filtererd uncompleted analytic derived rings follows immidiately from our hypothesis (1) and (3) on $C$. 

\begin{lemma} \label{erererere}
The cartesian diagram is left adjointable, i.e. there is a commutative diagram
\[
\begin{tikzcd}
\on{Mod}_{A} \on{Cond}\on{Fil}^{\wedge}D_{} (\mathbb{Z})
\arrow[d,"\on{gr}^{\bullet}"]
\arrow[r,"L_A"]
&
D(A) 
\arrow[d,"\on{gr}^{\bullet}"]
\\
\on{Mod}_{\on{gr}^{\bullet}A} \on{Cond}\on{Gr}D_{} (\mathbb{Z})
\arrow[r,"L_{\on{gr}^{\bullet} A}"]
&
C
\end{tikzcd}
\]
where the horizontal arrows are the left adjoint to the inclusion.
\end{lemma}

\begin{proof}
By lemma \ref{urowewqpepwp}, it is enough to prove that 
the right adjoint $G$ of $\on{gr}^{\bullet}$, sending a graded object $X^{\bullet}$ to the filtered object $\cdots \rightarrow X^{i+1} \overset{0}{\rightarrow} X^i\overset{0}{\rightarrow} \cdots$
with $0$ transition maps (on underlying objects in $\on{Gr} \on{Cond} D(\mathbb{Z})$), has the property that the image of $G\circ U_{\on{gr}^{\bullet} A}$ lies in the full subcategory $D(A)$. Indeed, since $D(A)$ is defined by pullback, this is equivalent to that the image of $\on{gr}^{\bullet} \circ G \circ U_{\on{gr}^{\bullet} A} $ lies in the image of $U_{\on{gr}^{\bullet} A}$, and we compute directly that 
$\on{gr}^{\bullet} \circ G(X^{\bullet}) = X^{\bullet} \oplus \Sigma X^{\bullet}[1]$, and  $C$
is closed under $\Sigma$ and $\oplus$ by the condition (1) of $C$. It is closed under $[1]$ by condition (3) of $C$.  
\end{proof}

We can now prove (2), namely the kernel of $L_A$ is a tensor ideal. Given $x \in \on{Ker}(L_A)$, and $y \in \on{Mod}_{A} \on{Cond}\on{Fil}^{\wedge}D_{} (\mathbb{Z})$, we compute that
\[
\on{gr}^{\bullet} \circ L_A(x\otimes y) \cong L_{\on{gr}^{\bullet} A}\circ \on{gr}^{\bullet} (x\otimes y) 
\cong L_{\on{gr}^{\bullet} A} \circ \on{gr}^{\bullet} (x)\otimes  L_{\on{gr}^{\bullet} A} \circ\on{gr}^{\bullet} (y) 
\cong  \on{gr}^{\bullet} \circ L_A(x) \otimes  L_{\on{gr}^{\bullet} A} \circ\on{gr}^{\bullet} (y) 
\cong 0.
\]
Since $\on{gr}^{\bullet}$ is conservative on completely filtered objects, this proves that $x\otimes y \in \on{Ker}(L_A)$.  

To prove axiom (4), note that  lemma \ref{56326} applied to $D(A) := \on{Mod}_{\mathbb{D}_-} C$ tells us that the monad $\hat{T}_{\on{gr}^{\bullet} A}$ induces a monad $\hat{T}_A$ on $D(A)$ such that $ \hat{T}_A \circ L_A=L_A \circ T_A$. 

 To prove last claim,  we apply \cite[proposition 4.1.9.]{Raksit2026}  to the subcategory $D(A)$ of $\on{Mod}_{A} \on{Cond}\on{Fil}^{\wedge}D_{} (\mathbb{Z})$ and the monad of tensoring with $C\tau$ on $\on{Mod}_{A} \on{Cond}\on{Fil}^{\wedge}D_{} (\mathbb{Z})$ to obtain the cartesian diagram 
\[
\begin{tikzcd}
\on{Mod}_{C\tau}\on{Mod}_{A} \on{Cond}\on{Fil}^{\wedge}D_{} (\mathbb{Z})
\arrow[r,hookleftarrow,"U_A"]
\arrow[d]
& 
\on{Mod}_{C\tau} D(A)
\arrow[d]
\arrow[dl, phantom, "\urcorner", very near start]
\\
\on{Mod}_{A} \on{Cond}\on{Fil}^{\wedge}D_{} (\mathbb{Z})
\arrow[r,hookleftarrow,"U_A"]
&
D(A) 
\end{tikzcd}
\] 
with vertical arrows forgetful functors. Composing with the cartesian diagram (\ref{2463764736723}), we have cartesian diagram
\[
\begin{tikzcd}
\on{Mod}_{C\tau}\on{Mod}_{A} \on{Cond}\on{Fil}^{\wedge}D_{} (\mathbb{Z})
\arrow[r,hookleftarrow,"U_A"]
\arrow[d,"G"]
& 
\on{Mod}_{C\tau} D(A)
\arrow[d]
\arrow[dl, phantom, "\urcorner", very near start]
\\
\on{Mod}_{A} \on{Cond}\on{Fil}^{\wedge}D_{} (\mathbb{Z})
\arrow[r,hookleftarrow,"U_A"]
\arrow[d,"\on{gr}^{\bullet}"]
&
D(A) 
\arrow[dl, phantom, "\urcorner", very near start]
\arrow[d,""]
\\
\on{Mod}_{\on{gr}^{\bullet}A} \on{Cond}\on{Gr}D_{} (\mathbb{Z})
\arrow[r,hookleftarrow,"U_{\on{gr}^{\bullet} A}"]
&
C.
\end{tikzcd}
\]
Now using the equivalence
 \[
\on{Mod}_{\on{gr}^{\bullet}A} \on{Cond}\on{Gr}D_{} (\mathbb{Z})  
\cong 
\on{Mod}_{C\tau} \on{Mod}_{A} \on{Cond}\on{Fil}^{\wedge}D_{} (\mathbb{Z})
\]
of \cite[proposition C.2.]{Burklund2026Galois}, and the fact that  the forgetting $C\tau$-module functor corresponds to the right adjoint $G$ of $\on{gr}^{\bullet}$ under this equivalence, we compute directly that 
\[
\on{gr}^{\bullet} \circ G (-) \cong \on{id}(-) \oplus \Sigma (-) [1]
\]
using the explicit description that $G$ sends a graded object $X^{\bullet}$ to the filtered object $\cdots \rightarrow X^{i+1} \overset{0}{\rightarrow} X^i\overset{0}{\rightarrow} \cdots$
with $0$ transition maps. It follows that $\on{Mod}_{C\tau} D(A)$ is the full subcategory of 
$\on{Mod}_{\on{gr}^{\bullet}A} \on{Cond}\on{Gr}D_{} (\mathbb{Z})$ 
such that $X^{\bullet} \in \on{Mod}_{C\tau} D(A)$ 
if and only if
$\on{gr}^{\bullet} \circ G (X^{\bullet})=X^{\bullet}\oplus \Sigma X^{\bullet}[1] \in C$, but $C$ is closed under colimits and $[1]$ by condition (1) and (3),  so  $\on{Mod}_{C\tau} D(A)$ as a full subcategory is exactly $C$. 
\end{proof}

We observe that the induced analytic structure construction also works in our the filtered analytic derived algebra setting.

\begin{lemma} \label{56326}
Let $C$ be a presentable symmetric monoidal stable category, and $T: C \rightarrow C$ be a monad on $C$ together with a map of monad $\on{Sym}^{\otimes}_{C} \rightarrow T$,
where $\on{Sym}^{\otimes}_C$ 
is the commutative algebra monad
$X \mapsto \oplus_n (X^{\otimes n})_{\Sigma_n}$.  It gives rise to a functor 
\[
\Theta_C: 
\on{LMod}_{T} C \longrightarrow \on{CAlg}(C)
\]
and for any $A \in \on{LMod}_{T} C$, and we define $\on{Mod}_A C$ by viewing $A$ as a commutative algebra through the functor. We assume also that $T$ preserves sifted colimits, so that $\Theta_C$ preserves both limits and sifted colimits as well. 

Let $U: D \subset C$ be a full  subcategory of $C$ with  left adjoint $L: C \rightarrow D$ whose kernel is a tensor ideal so that $D$ is equipped with a symmetric monoidal structure such that $L$ is a  symmetric monoidal functor. We assume also that $D$ is closed under both limits and colimits.  

Assume that $U \circ L \circ T = T \circ U \circ L$, equivalently, $T$ preserves the class of arrows which becomes equivalence after applying $L$, then $\hat{T} :=L \circ T \circ U $ is a monad on $D$. Let $A \in \on{LMod}_{T} C$, then $T$  (resp. $\hat{T}$) induces a monad $T_A$ (resp. $\hat{T}_A$) on $\on{Mod}_AC$ (resp. $\on{Mod}_A D = \on{Mod}_{L(A)} D$) such that 
$\hat{T}_A \circ L_A = L_A \circ T_A$. 
\end{lemma}

\begin{proof}
By \cite[proposition 4.1.9.]{Raksit2026},
    $\hat{T}$ is a monad and we have that 
    \[
    U:
    \on{LMod}_{\hat{T}} D \subset \on{LMod}_{T} C 
    \]
    is the full subcategory $\on{LMod}_{T_A} C \times_{C} D$ with left adjoint induced by $L$. In other words, we have a commutative diagrams
    \begin{equation} \label{238938493}
    \begin{tikzcd}
    \on{LMod}_{T} C 
    \arrow[d]
    &
          \on{LMod}_{\hat{T}} D
          \arrow[l,hook, "U"] 
          \arrow[d]
        \\
       C 
       & 
       D
       \arrow[l,hook,"U"]
    \end{tikzcd}
    \end{equation}
    \begin{equation} \label{23437439fe}
    \begin{tikzcd}
    \on{LMod}_{T} C 
    \arrow[r,"L"]
    \arrow[d]
    &
          \on{LMod}_{\hat{T}} D 
          \arrow[d]
        \\
       C 
       \arrow[r,"L"]
       & 
       D
    \end{tikzcd}
    \end{equation}
    where the vertical maps are forgetful functors. 

    We observe that the map of monads $\on{Sym}^{\otimes}_{C} \rightarrow T$
    induces a map of monads 
    \[
    \on{Sym}^{\otimes}_{D} \longrightarrow \hat{T}
    \]
    since $L$ is by our assumption symmetric monoidal. Moreover, $\hat{T}$ preserves sifted colimits as well, since $T$ does so and both $L$ and $U$ preserves all colimits by our assumption. It implies that the functor 
    \[
    \Theta_D: \on{LMod}_{\hat{T}} D
    \longrightarrow 
    \on{CAlg}(D)
    \]
    preserves both limits and sifted colimits. Indeed, we know that the forgetful fucntor  to $D$ is conservative and  creates limits  (see \cite[corollary 4.2.3.3.]{Higheralgebra}), and it also creates sifted colimits since both the monads $\hat{T}$ and $\on{Sym}^{\otimes}_{D}$
    preserves sifted colimits, see \cite[corollary 4.2.3.5.]{Higheralgebra}. 

    The forgetful functor 
    \[
    D \longleftarrow \on{LMod}_{\hat{T}} D: U
    \]
    induces a monadic adjunction 
    \[
    \hat{\on{Sym}}_{L(A)}: \on{Mod}_{L(A)} D \leftrightarrows 
    (\on{LMod}_{\hat{T}} D)_{L(A)/} :U
    \]
    by passing to coslice and module categories. More precisely, the right adjoint is the composite 
    \[
    \on{Mod}_{L(A)} D
    \longleftarrow
    \on{CAlg}(D)_{L(A)/} 
    \overset{\Theta_D}{\longleftarrow}
    (\on{LMod}_{\hat{T}} D)_{L(A)/} :U
    \]
    and it has a monadic left adjoint $ \hat{\on{Sym}}_{L(A)}$ since $U$ is conservative and preserves all limits and sifted colimits. 
Similarly, we have monadic adjunctions
\[
    \on{Sym}_A: \on{Mod}_{A} C\leftrightarrows 
    (\on{LMod}_{T} C)_{A/} :U. 
    \]

    Let 
    \[
    \hat{T}_A:= U \circ \hat{\on{Sym}}_{L(A)}
    \]
    and 
    \[
    T_A:= U \circ \on{Sym}_A
    \]
    and we want to prove that $\hat{T}_A \circ L_A = L_A \circ T_A$, where 
    \[
    L_A:  \on{Mod}_{L(A)} C \longrightarrow  \on{Mod}_{L(A)} D
    \]
    is the left adjoint to the inclusion $ \on{Mod}_{L(A)} D \subset  \on{Mod}_{L(A)} C$ induced by $L$. 

    The commutative diagrams (\ref{238938493}) and (\ref{23437439fe})
    induces commutative diagrams
    \begin{equation} \label{238933}
    \begin{tikzcd}
    (\on{LMod}_{T} C )_{A/}
    \arrow[d,"U"]
    &
          (\on{LMod}_{\hat{T}} D)_{L(A)/}
          \arrow[l,hook, "U_A"] 
          \arrow[d,"U"]
        \\
     \on{Mod}_A  C 
       & 
     \on{Mod}_{L(A)}  D
       \arrow[l,hook,"U_A"]
    \end{tikzcd}
    \end{equation}
    \begin{equation} \label{27439fe}
    \begin{tikzcd}
    (\on{LMod}_{T} C )_{A/}
    \arrow[r,"L_A"]
    \arrow[d,"U"]
    &
          (\on{LMod}_{\hat{T}} D )_{L(A)/}
          \arrow[d,"U"]
        \\
      \on{Mod}_A C 
       \arrow[r,"L_A"]
       & 
       \on{Mod}_{L(A)}  D
    \end{tikzcd}
    \end{equation}
    by passing to the coslice/module categories. 
Taking left adjoint of the diagram (\ref{238933}), we have a commutative diagram
\begin{equation} \label{238933324793749}
    \begin{tikzcd}
    (\on{LMod}_{T} C )_{A/}
    \arrow[r,"L_A"]
    &
          (\on{LMod}_{\hat{T}} D)_{L(A)/}
        \\
     \on{Mod}_A  C 
     \arrow[u,"\on{Sym}_A"]
     \arrow[r,"L_A"]
       & 
     \on{Mod}_{L(A)}  D
       \arrow[u,"\hat{\on{Sym}}_{L(A)}"]
    \end{tikzcd}
    \end{equation}
    Composing (\ref{238933324793749}) and (\ref{27439fe}) 
    we have that $\hat{T}_A \circ L_A = L_A \circ T_A$. 
\end{proof}

\begin{corollary} \label{74329874389}
    Let $(A, D(A))$ be an uncompleted filtered/graded analytic derived ring, and $B \in \on{DAlg}\on{Cond}\on{Fil}D_{} (\mathbb{Z})/\on{DAlg}\on{Cond}\on{Gr}D_{} (\mathbb{Z}) $ be an $A$-algebra, then $(B, \on{Mod}_BD(A))$ is also an uncompleted filtered/graded analytic derived ring. 
\end{corollary}

\begin{proof}
    The axiom (1), (2) and (3) for $\on{Mod}_BD(A)$ is clear. (4) follows from the lemma applied to $C = \on{Mod}_A \on{Cond}\on{Fil}D_{} (\mathbb{Z})/\on{Mod}_A \on{Cond}\on{Gr}D_{} (\mathbb{Z}) $ 
    and $D= D(A)$. 
\end{proof}

\begin{definition}
    A map of (uncompleted) filtered/graded analytic derived rings 
    $(A, D(A)) \longrightarrow (B, D(B))$
    is a map 
    $A \longrightarrow B$ in $\on{DAlg}\on{Cond}\on{Fil}D_{} (\mathbb{Z})/\on{DAlg}\on{Cond}\on{Gr}D_{} (\mathbb{Z})$ such that the restriction functor sends $D(B)$ to $D(A)$. 

    The ($\infty$-)category of uncompleted filtered/graded analytic derived rings is denoted $\on{FilAnDAG}^{\on{un}}/\on{GrAnDAG}^{\on{un}}$, which has a full subcategory $\on{FilAnDAG}/\on{GrAnDAG}$ 
    of filtered/graded analytic derived rings. 
\end{definition}

\begin{lemma} \label{weuryeiwyrue}
The category
    $\on{FilAnDAG}^{\on{un}}/\on{GrAnDAG}^{\on{un}}$ has all colimits, given a diagram $(A_i, D(A_i))$ of uncompleted filtered/graded analytic derived rings, the colimit is $(A, D(A))$, where $A= \on{colim}_i A_i$ computed in $\on{DAlg}\on{Cond}\on{Fil}D_{} (\mathbb{Z})/\on{DAlg}\on{Cond}\on{Gr}D_{} (\mathbb{Z})$, and $D(A)= \cap_i \on{Mod}_A D(A_i) $ where the intersection is limits of the  full subcategories $\on{Mod}_A D(A_i) \subset \on{Mod}_A \on{Mod}_{A_i} \on{Cond}\on{Fil}^{\wedge}D_{} (\mathbb{Z}) = \on{Mod}_A  \on{Cond}\on{Fil}^{\wedge}D_{} (\mathbb{Z})$, resp. $\on{Mod}_A D(A_i) \subset \on{Mod}_A \on{Mod}_{A_i} \on{Cond}\on{Gr}D_{} (\mathbb{Z}) = \on{Mod}_A  \on{Cond}\on{Gr}D_{} (\mathbb{Z})$. 
\end{lemma}

\begin{proof}
    This is identical to the proof of \cite[proposition 4.1.11.]{Notesonsolidgeometry}. We first check that $(A, D(A))$ is a filtered/graded uncompleted analytic derived ring, the conditions (1), (2) and (3) are clear. To see (4), we observe that an arrow  $f$ in 
    $\on{Mod}_A  \on{Cond}\on{Fil}^{\wedge}D_{} (\mathbb{Z})/\on{Mod}_A  \on{Cond}\on{Gr}D_{} (\mathbb{Z})$
    becomes an equivalence in  the intersection $D(A)$ if and only if it becomes an equivalence in $\on{Mod}_A D(A_i)$ for every $i$, and $T_A$ preserves them by applying lemma \ref{56326} to $\on{Mod}_A D(A_i)$. 

    We can now see that $(A, D(A))$ is indeed a colimit. Recall from the definition of  $\on{FilAnDAG}^{\on{un}}$ that 
    \[
    \on{Map}_{\on{FilAnDAG}^{\on{un}}}((A,D(A)), (B,D(B)) \subset 
    \on{Map}_{\on{DAlg}\on{Cond}\on{Fil}D_{} (\mathbb{Z})} (A, B)
    \]
    is the sub-anima consisting of the connected components of $\on{Map}_{\on{DAlg}\on{Cond}\on{Fil}D_{} (\mathbb{Z})} (A, B)$ which has $\pi_0$ those maps $A \longrightarrow B$  that sends $D(B)$ to $D(A)$. We see  from the definition of $(A,D(A))$ that the left hand side corresponds to the sub-anima $\{f: A \longrightarrow B\}$ of the right hand side such that $D(B)$ is sent to $D(A_i)$ under $A_i \rightarrow A \overset{f}{\rightarrow} B$ for every $i$, and this is exactly $\on{lim}_i \on{Map}_{\on{FilAnDAG}^{\on{un}}}((A_i,D(A_i)), (B,D(B))$. Similar for the graded case. 
\end{proof}

\begin{proposition} \label{wyrw3874}
    The full subcategory $\on{FilAnDAG} \subset \on{FilAnDAG}^{\on{un}}$ of filtered analytic derived rings has a left adjoint which sends $(A, D(A))$ to $(L_A(A), D(A))$. In particular, $\on{FilAnDAG} $ has all colimits. 

    Similarly, the full subcategory $\on{GrAnDAG} \subset \on{GrAnDAG}^{\on{un}}$ of filtered analytic derived rings has a left adjoint which sends $(A, D(A))$ to $(L_A(A), D(A))$, and $\on{GrAnDAG} $ has all colimits.
\end{proposition}

\begin{proof}
    This is the same as the proof of \cite[proposition 2.3.12.]{Mann}, one observes directly that $(L_A(A), D(A))$ is the initial object in $\on{FilAnDAG}$ lies over $(A, D(A))$. The only non-trivial thing is that $L_A(A)$ actually lives in $\on{DAlg}\on{Cond}\on{Fil}D_{} (\mathbb{Z})$, but this follows from our axiom (4) of filtered uncompleted analytic derived rings. More precisely, it follows from the commutative diagram (\ref{27439fe}) when $C = \on{Mod}_A \on{Cond}\on{Fil}D_{} (\mathbb{Z}) $ 
    and $D= D(A)$. The graded case is identical. 
\end{proof}

Let us describe a class of examples of filtered analytic derived rings which motivates the definition given here.

\begin{proposition} \label{xnvmnvd}
Let $A \in \on{DAlg}\on{Cond}\on{Fil}D_{} (\mathbb{Z})$ 
   and we are given an uncompleted analytic ring structure on $\tau_{\geq 0}\on{gr}^0A$, i.e. a full subcategory  
    \[
D(\tau_{\geq 0}\on{gr}^0A)  
\subset
\on{Mod}_{\tau_{\geq0}\on{gr}^0A} \on{Cond}D_{} (\mathbb{Z}) 
    \]
    that is closed under limits and colimits, the kernel of the left adjoint to the inclusion is a tensor ideal, and the left adjoint preserves connective objects. Then 
    \[
(\on{gr}^{\bullet} A,  \on{Mod}_{\on{gr}^{\bullet} A} \on{Gr}D_{} (\tau_{\geq 0}\on{gr}^0A))
    \]
    is a graded uncompleted analytic derived ring. 
    
    Let $D(A)$ be the full subcategory
    \[
    D(A) \subset \on{Mod}_{A} \on{Cond}\on{Fil}^{\wedge}D_{} (\mathbb{Z})
    \]
    such that $M \in D(A)$ if and only if $\on{gr}^{\bullet}M$ lies in the full subcategory 
    \[
    \on{Mod}_{\on{gr}^{\bullet} A} \on{Gr}D_{} (\tau_{\geq 0}\on{gr}^0A) 
\subset    
\on{Mod}_{\on{gr}^{\bullet}A} \on{Cond}\on{Gr}D_{} (\mathbb{Z}),
    \]
    then $(A, D(A))$ is a filtered analytic derived ring. Equivalently, 
    \[
    D(A) = \on{Mod}_{\mathbb{D}_-}
    \on{Mod}_{\on{gr}^{\bullet} A} \on{Gr}D_{} (\tau_{\geq 0}\on{gr}^0A). 
    \]

    Moreover, the left adjoint $L_A$ to the inclusion $D(A) \subset \on{Mod}_{A} \on{Cond}\on{Fil}^{\wedge}D_{} (\mathbb{Z})$ satisfies 
    \[
    \on{gr}^{i} \circ L_A = 
    L_{\tau_{\geq 0}\on{gr}^0A}
    \circ \on{gr}^{i}.
    \]
\end{proposition}

\begin{proof}
We first explain the embedding 
\[
\on{Mod}_{\on{gr}^{\bullet} A} \on{Gr}D_{} (\tau_{\geq 0}\on{gr}^0A)
\subset
\on{Mod}_{\on{gr}^{\bullet}A} \on{Cond}\on{Gr}D_{} (\mathbb{Z}).
\]
This is induced by the map of condensed graded derived rings 
\begin{equation} \label{qeiruhx}
\on{ins}^0 \tau_{\geq 0} \on{gr}^0A \longrightarrow \on{ins}^0\on{gr}^0A \longrightarrow \on{gr}^{\bullet}A
\end{equation}
where 
\[\on{ins}^0:  
\on{DAlg}\on{Cond}D_{} (\mathbb{Z})
\longrightarrow 
\on{DAlg}\on{Cond}\on{Gr}D_{} (\mathbb{Z}) 
\]
is the left adjoint of $\on{ev}^0$, which sends an ungraded object to a graded object concentrated in degree 0, see \cite[remark 4.3.6.]{Raksit2026}. The first map is the functoriality of $\on{ins}^0$, and the second is the counit of the adjunction. We have 
\[
\on{Mod}_{\on{ins}^0 \tau_{\geq 0} \on{gr}^0A} \on{Cond}\on{Gr}D_{} (\mathbb{Z}) 
\cong 
\on{Gr}
\on{Mod}_{\tau_{\geq 0} \on{gr}^0A} \on{Cond}D_{} (\mathbb{Z}).
\]
Indeed, $\on{ins}^{0}: \on{Cond}D_{} (\mathbb{Z})
\longrightarrow
\on{Cond}\on{Gr}D_{} (\mathbb{Z})$
is symmetric monoidal, and 
we have that that 
 \begin{align*}
 \on{Mod}_{\on{ins}^0 \tau_{\geq 0} \on{gr}^0A} \on{Cond}\on{Gr}D_{} (\mathbb{Z})
& \cong 
 \on{Mod}_{ \tau_{\geq 0} \on{gr}^0A} \on{Cond}D_{} (\mathbb{Z}) \otimes_{\on{Cond}D_{} (\mathbb{Z}), \on{ins}^0 } 
 \on{Cond}\on{Gr}D_{} (\mathbb{Z})
\\
&\cong \on{Mod}_{ \tau_{\geq 0} \on{gr}^0A} \on{Cond}D_{} (\mathbb{Z}) \otimes_{\on{Cond}D_{} (\mathbb{Z}), \on{ins}^0 } 
\on{Gr} \on{Cond}D_{} (\mathbb{Z})
\\
& \cong 
\on{Mod}_{ \tau_{\geq 0} \on{gr}^0A} \on{Cond}D_{} (\mathbb{Z}) \otimes_{\on{Cond}D_{} (\mathbb{Z}), \on{ins}^0 } 
( \on{Cond}D_{} (\mathbb{Z}) \otimes \on{Gr} \on{Sp} )
\\
&\cong \on{Gr}\on{Mod}_{ \tau_{\geq 0} \on{gr}^0A} \on{Cond}D_{} (\mathbb{Z}) 
\end{align*}
where we use \cite[Theorem 4.8.4.6.]{Higheralgebra} in the first equivalence, and \cite[remark 3.1.4]{Raksit2026} in the third and fourth equivalence.
We now have 
\[
\on{Mod}_{\on{gr}^{\bullet}A} \on{Cond}\on{Gr}D_{} (\mathbb{Z}) 
\cong 
\on{Mod}_{\on{gr}^{\bullet}A} \on{Mod}_{\on{ins}^0 \tau_{\geq 0} \on{gr}^0A}\on{Cond}\on{Gr}D_{} (\mathbb{Z})
\]
\[
\cong 
\on{Mod}_{\on{gr}^{\bullet}A} \on{Gr}
\on{Mod}_{\tau_{\geq 0} \on{gr}^0A} \on{Cond}D_{} (\mathbb{Z}).
\]
Since $\tau_{\geq 0}\on{gr}^0A$ is an analytic ring, 
$D_{} (\tau_{\geq 0}\on{gr}^0A) \subset \on{Mod}_{\tau_{\geq 0} \on{gr}^0A} \on{Cond}D_{} (\mathbb{Z})$ 
is a symmetric monoidal full subcategory stable under limits and colimits, with the localization functor being symmetric monoidal and preserving connective objects.  We see immediately that the same holds for the graded subcategory 
\[
R:
\on{Gr}D_{} (\tau_{\geq 0}\on{gr}^0A) 
\subset 
\on{Gr}
\on{Mod}_{\tau_{\geq 0} \on{gr}^0A} \on{Cond}D_{} (\mathbb{Z}),
\]
where both sides are equipped with the neutral t-structure as defined in \cite[construction 3.3.2.]{Raksit2026}, i.e. a graded object is connective/coconnective if and only if each graded piece is connective/coconnective.  
We have in particular 
$\on{Gr}D_{} (\tau_{\geq 0}\on{gr}^0A)$
is enriched over 
$\on{Gr}
\on{Mod}_{\tau_{\geq 0} \on{gr}^0A} \on{Cond}D_{} (\mathbb{Z})$
through the left adjoint to the inclusion, so $\on{Mod}_{\on{gr}^{\bullet} A} \on{Gr}D_{} (\tau_{\geq 0}\on{gr}^0A)$ makes sense. 

Now the fully faithful embedding being right adjoint to a symmetric monoidal functor is lax symmetric monoidal, so it induces a functor
\[
U_{\on{gr}^{\bullet} A}:
\on{Mod}_{\on{gr}^{\bullet} A} \on{Gr}D_{} (\tau_{\geq 0}\on{gr}^0A) 
\subset 
\on{Mod}_{\on{gr}^{\bullet}A} \on{Gr}
\on{Mod}_{\tau_{\geq 0} \on{gr}^0A} \on{Cond}D_{} (\mathbb{Z}).
\]

By proposition \ref{2378937583}, it is enough to prove that 
$\on{Mod}_{\on{gr}^{\bullet} A} \on{Gr}D_{} (\tau_{\geq 0}\on{gr}^0A) $
satisfies the conditions (1) to (4) there, from which (1), (2) and (3) are easy to see. 

\begin{lemma}
    $\on{Mod}_{\on{gr}^{\bullet} A} \on{Gr}D_{} (\tau_{\geq 0}\on{gr}^0A) $
    is  a full  subcategory  
stable under limits and colimits and $[1]$. Moreover, the left adjoint to $U_{\on{gr}^{\bullet} A}$ is symmetric monoidal. 
\end{lemma}

\begin{proof}
First, $\on{Gr}D_{} (\tau_{\geq 0}\on{gr}^0A) \subset \on{Gr}
\on{Mod}_{\tau_{\geq 0} \on{gr}^0A} \on{Cond}D_{} (\mathbb{Z})$
is closed under limits, colimits since this is true for $D_{} (\tau_{\geq 0}\on{gr}^0A) \subset \on{Mod}_{\tau_{\geq 0} \on{gr}^0A} \on{Cond}D_{} (\mathbb{Z})$ by definition of analytic animated rings, and limits and colimits in graded categories are computed pointwise. $\on{Gr}D_{} (\tau_{\geq 0}\on{gr}^0A) $ is closed under $[1]$ since it is a graded category itself. Passing to module categories preserves these properties, since they are computed underlying.

The symmetric monoidal Bousfield localization
\[
\on{Gr}
\on{Mod}_{\tau_{\geq 0} \on{gr}^0A} \on{Cond}D_{} (\mathbb{Z})
\longrightarrow 
\on{Gr}D_{} (\tau_{\geq 0}\on{gr}^0A) 
\]
induces a symmetric monoidal left adjoint
\[
L_{\on{gr}^{\bullet} A}:
\on{Mod}_{\on{gr}^{\bullet}A} \on{Gr}
\on{Mod}_{\tau_{\geq 0} \on{gr}^0A} \on{Cond}D_{} (\mathbb{Z})
\longrightarrow 
\on{Mod}_{\on{gr}^{\bullet} A} \on{Gr}D_{} (\tau_{\geq 0}\on{gr}^0A), 
\]
whose right adjoint is also computed in the underlying category by \cite[example 7.3.2.8.]{Higheralgebra}, i.e. 
we have a commutative diagram
\[
\begin{tikzcd}
 \on{Mod}_{\on{gr}^{\bullet}A} \on{Gr}
\on{Mod}_{\tau_{\geq 0} \on{gr}^0A} \on{Cond}D_{} (\mathbb{Z})   
\arrow[d]
\arrow[r,hookleftarrow, "U_{\on{gr}^{\bullet} A}"]
&
\on{Mod}_{\on{gr}^{\bullet} A} \on{Gr}D_{} (\tau_{\geq 0}\on{gr}^0A) 
\arrow[d]
    \\
 \on{Gr}
\on{Mod}_{\tau_{\geq 0} \on{gr}^0A} \on{Cond}D_{} (\mathbb{Z})
\arrow[r,hookleftarrow, ""]
&
\on{Gr}D_{} (\tau_{\geq 0}\on{gr}^0A)
\end{tikzcd}
\]
with vertical arrows forgetful functors. The claims on $U_{\on{gr}^{\bullet} A}$ follow since the forgetful functors are conservative. 
\end{proof}

To prove (4), we observe that the monad $U \circ \on{Sym}^{\bullet}_{\tau_{\geq0}\on{gr}^0A}$
on
$
\on{Mod}_{\tau_{\geq0}\on{gr}^0A} \on{Cond}D_{} (\mathbb{Z})$
induces a monad on 
$D(\tau_{\geq 0}\on{gr}^0A)  
$ 
by \cite[proposition 4.2.6.]{Notesonsolidgeometry}, by passing to graded categories we obtain a monad on $\on{Gr} D(\tau_{\geq 0}\on{gr}^0A) $ compatible with the  one on 
$ \on{Gr} \on{Mod}_{\tau_{\geq0}\on{gr}^0A} \on{Cond}D_{} (\mathbb{Z})$. 
Passing to module categories gives rise to a monad on 
$ \on{Mod}_{\on{gr}^{\bullet} A} \on{Gr}D_{} (\tau_{\geq 0}\on{gr}^0A) $ 
compatible with one on $\on{Mod}_{\on{gr}^{\bullet} A} \on{Gr}\on{Mod}_{\tau_{\geq0}\on{gr}^0A} \on{Cond}D_{} (\mathbb{Z})$ by
 lemma \ref{56326}.

 Finally, the last claim follows from lemma \ref{erererere} and $L_{\on{gr}^{\bullet}A}  $  is induced from the underlying left adjoint, i.e. we have a commutative diagram
 \[
 \begin{tikzcd}
     \on{Mod}_{\on{gr}^{\bullet}A} \on{Cond}\on{Gr}D_{} (\mathbb{Z})
\arrow[r,"L_{\on{gr}^{\bullet} A}"]
\arrow[d]
&
\on{Mod}_{\on{gr}^{\bullet} A} \on{Gr}D_{} (\tau_{\geq 0}\on{gr}^0A)
\arrow[d]
\\
\on{Gr}
\on{Cond}D_{} (\mathbb{Z})
\arrow[r,"\on{Gr}L_{\tau_{\geq 0}\on{gr}^0A}"]
&
\on{Gr}D_{} (\tau_{\geq 0}\on{gr}^0A)
 \end{tikzcd}
 \]
 with vertical arrows forgetful functors. 
\end{proof}

We observe that the proof above actually proves the following proposition. 

\begin{proposition} \label{wury9387r9}
Let $A \in \on{DAlg}\on{Cond}\on{Gr}D_{} (\mathbb{Z})$ 
   and we are given an uncompleted analytic ring structure on $\tau_{\geq 0}\on{ev}^0A$, i.e. a full subcategory  
    \[
D(\tau_{\geq 0}\on{ev}^0A)  
\subset
\on{Mod}_{\tau_{\geq0}\on{ev}^0A} \on{Cond}D_{} (\mathbb{Z}) 
    \]
    that is closed under limits and colimits, the kernel of the left adjoint to the inclusion is a tensor ideal, and the left adjoint preserves connective objects. Then 
    \[
    (A, D(A)) :=
( A,  \on{Mod}_{ A} \on{Gr}D_{} (\tau_{\geq 0}\on{gr}^0A))
    \]
    is a graded uncompleted analytic derived ring. 
  
    Moreover, the left adjoint $L_A$ to the inclusion $D(A) \subset \on{Mod}_{A} \on{Cond}\on{Gr}D_{} (\mathbb{Z})$ satisfies 
    \[
    \on{ev}^{i} \circ L_A = 
    L_{\tau_{\geq 0}\on{ev}^0A}
    \circ \on{ev}^{i}.
    \]
\end{proposition}

We record the functoriality of the construction. 

\begin{lemma} \label{79832403983}
    Let $A \in \on{DAlg}\on{Cond}\on{Fil}D_{} (\mathbb{Z})$,  
   $\tau_{\geq 0}\on{gr}^0A$ 
be given an analytic animated ring structure, and 
   $D(A) \subset \on{Mod}_{A} \on{Cond}\on{Fil}^{\wedge}D_{} (\mathbb{Z})$
    as in proposition \ref{xnvmnvd}, then for any $(B,D(B)) \in \on{FilAnDAG}^{\on{un}}$, a map 
    $
    A \longrightarrow B 
    $
    in 
    $\on{DAlg}\on{Cond}\on{Fil}D_{} (\mathbb{Z})$
    induces a map in $\on{FilAnDAG}^{\on{un}}$ if and only if for each $M \in D(B)$, 
    \[
    \on{gr}^iM \in D(\tau_{\geq 0}\on{gr}^0A)
    \]
    for all $i$ under the forgetful functor 
    \[
    \on{Mod}_{\on{gr}^{\bullet}B} \on{Cond}\on{Gr}D_{} (\mathbb{Z})
    \longrightarrow 
    \on{Gr}
    \on{Mod}_{\on{gr}^{0}B} \on{Cond}D_{} (\mathbb{Z})
    \longrightarrow 
    \on{Gr}
    \on{Mod}_{\tau_{\geq 0}\on{gr}^{0}A} \on{Cond}D_{} (\mathbb{Z})
    \]
    along the map of graded rings
    $ \on{ins}^0 \tau_{\geq 0}\on{gr}^{0}A 
    \longrightarrow
   \on{ins}^0 \on{gr}^{0}A 
    \longrightarrow
    \on{gr}^{\bullet}A 
    \longrightarrow 
     \on{gr}^{\bullet}B$.
     
     In particular, if the analytic structure on  $B$ is also induced from an analytic ring structure on $  \tau_{\geq 0}\on{gr}^{0}B$, then a map $A \longrightarrow B$ in 
     $\on{DAlg}\on{Cond}\on{Fil}D_{} (\mathbb{Z})$ 
     induces a map in $\on{FilAnDAG}^{\on{un}}$ if the induced map 
     $\tau_{\geq 0}\on{gr}^{0}A
     \longrightarrow \tau_{\geq 0}\on{gr}^{0}B$
     is a map of analytic animated rings.
\end{lemma}

\begin{lemma} \label{wityuoiewuroe}
       Let $A \in \on{DAlg}\on{Cond}\on{Gr}D_{} (\mathbb{Z})$,  
   $\tau_{\geq 0}\on{ev}^0A$ 
be given an analytic animated ring structure, and 
   $D(A) \subset \on{Mod}_{A} \on{Cond}\on{Gr}D_{} (\mathbb{Z})$
    as in proposition \ref{wury9387r9}, then for any $(B,D(B)) \in \on{GrAnDAG}^{\on{un}}$, a map 
    $
    A \longrightarrow B 
    $
    in 
    $\on{DAlg}\on{Cond}\on{Gr}D_{} (\mathbb{Z})$
    induces a map in $\on{GrAnDAG}^{\on{un}}$ if and only if for each $M \in D(B)$, 
    \[
    \on{ev}^iM \in D(\tau_{\geq 0}\on{ev}^0A)
    \]
    for all $i$ under the forgetful functor 
    \[
    \on{Mod}_{B} \on{Cond}\on{Gr}D_{} (\mathbb{Z})
    \longrightarrow 
    \on{Gr}
    \on{Mod}_{\on{ev}^{0}B} \on{Cond}D_{} (\mathbb{Z})
    \longrightarrow 
    \on{Gr}
    \on{Mod}_{\tau_{\geq 0}\on{ev}^{0}A} \on{Cond}D_{} (\mathbb{Z})
    \]
    along the map of graded rings
    $ \on{ins}^0 \tau_{\geq 0}\on{ev}^{0}A 
    \longrightarrow
   \on{ins}^0 \on{ev}^{0}A 
    \longrightarrow
   A 
    \longrightarrow 
     B$.
     
     In particular, if the analytic structure on  $B$ is also induced from an analytic ring structure on $  \tau_{\geq 0}\on{ev}^{0}B$, then a map $A \longrightarrow B$ in 
     $\on{DAlg}\on{Cond}\on{Gr}D_{} (\mathbb{Z})$ 
     induces a map in $\on{GrAnDAG}^{\on{un}}$ if the induced map 
     $\tau_{\geq 0}\on{ev}^{0}A
     \longrightarrow \tau_{\geq 0}\on{ev}^{0}B$
     is a map of analytic animated rings.
\end{lemma}

We give another example of filtered analytic derived rings that will be important for us. 

\begin{example} \label{xfower870}
    Let $\mathbb{T}_{\on{Fil}}^{\vee} := \tau_{\geq \bullet} (\mathbb{Z}^{\mathbb{S}^1}) \in \on{Fil} D(\mathbb{Z})$ 
    be the (dual) $\mathbb{Z}$-linear  filtered circle as defined in \cite[notation 6.1.4.]{Raksit2026}. In other words, it is the Postnikov filtration of the complex representing singular cohomology of the circle $\mathbb{S}^1$ with $\mathbb{Z}$-coefficients. By \cite[theorem 6.1.5.]{Raksit2026}, we have 
    \[
    \mathbb{T}_{\on{Fil}}^{\vee} \in \on{cCAlg}(\on{DAlg} \on{Fil} D(\mathbb{Z}))
    \]
    which means it is a cocommutative coalgebra object in the category $\on{DAlg} \on{Fil} D(\mathbb{Z})$ with cocartesian symmetric monoidal structure. We observe that 
    \[
     \on{gr}^{\bullet} \mathbb{T}_{\on{Fil}}^{\vee} 
     \cong 
     \mathbb{D}_+^{\vee}
     \in \on{cCAlg}(\on{DAlg} \on{Gr} D(\mathbb{Z})),
    \]
    where $ \mathbb{D}_+^{\vee} = \mathbb{Z} \oplus \Sigma^{-1} \mathbb{Z} [-1]$
    is defined as in \cite[proposition 5.1.7.]{Raksit2026}. 

    We view $  \mathbb{T}_{\on{Fil}}^{\vee} $
    as an element 
    \[
\underline{\mathbb{T}_{\on{Fil}}^{\vee}} \in
    \on{DAlg}\on{Cond}\on{Fil}D_{} (\mathbb{Z})
    \]
    through the canonical symmetric monoidal functor 
    \begin{equation} \label{xcnvksjh}
    \on{DAlg}\on{Fil}D_{} (\mathbb{Z})
    \longrightarrow 
    \on{DAlg}\on{Cond}\on{Fil}D_{} (\mathbb{Z}).
    \end{equation}
    In other words, the functor sends $A \in \on{DAlg}\on{Fil}D_{} (\mathbb{Z}) $ to the constant sheaf with value $A$ on light profinite sets. We give the trivial filtered analytic derived ring structure on $\underline{\mathbb{T}_{\on{Fil}}^{\vee}}$, namely 
    \[
    D(\underline{\mathbb{T}_{\on{Fil}}^{\vee}})= \on{Mod}_{\mathbb{T}_{\on{Fil}}^{\vee}} 
    \on{Cond}\on{Fil}^{\wedge}D_{} (\mathbb{Z}). 
    \]
    Moreover, we see that 
    \[
   (\underline{ \mathbb{T}_{\on{Fil}}^{\vee}}, D(\underline{\mathbb{T}_{\on{Fil}}^{\vee}})) \in \on{cCAlg}(\on{FilAnDAG}^{\on{un}}) 
    \]
    where we give $\on{FilAnDAG}^{\on{un}}$ the cocartesian symmetric monoidal structure. Indeed, the symmetric monoidal functor (\ref{xcnvksjh}) preserves cocommutative algebra structures, and so does the functor from $\on{DAlg}\on{Cond}\on{Fil}^{}D_{} (\mathbb{Z})$ to $\on{FilAnDAG}^{\on{un}}$ equipping $A \in \on{DAlg}\on{Cond}\on{Fil}^{}D_{} (\mathbb{Z})$ with the trivial analytic structures. 

    Taking graded functor and we have 
    \[
    (\underline{ \mathbb{D}^{\vee}_+}, D(\underline{ \mathbb{D}^{\vee}_+}) ) \in \on{cCAlg}(\on{GrAnDAG}^{\on{un}}),
    \]
    where 
    \[
    D(\underline{ \mathbb{D}^{\vee}_+}):=
    \on{Mod}_{\mathbb{D}^{\vee}_+} \on{Cond} \on{Gr} D(\mathbb{Z}). 
    \]
\end{example}

\section{Hodge completed derived de Rham complex of analytic animated rings} \label{3958203834}

We take a hands-on approach to define the Hodge completed derived de Rham complex of analytic rings. It is possible to define it using universal properties as in \cite{Raksit2026}. 

Let us first define the Hodge completed derived de Rham complex for condensed animated rings, without the analytic structure. 

\begin{definition}
    Let $A \longrightarrow B$ be a map of condensed derived rings, then we define the Hodge completed derived de Rham complex 
    \[
    \on{L} \Omega_{B/A}^{+, \bullet}  \in \on{cMod}_{\mathbb{D}_+^{\vee}}  (\on{DAlg}\on{Cond} \on{Gr}D(\mathbb{Z})_{\on{ins}^0 A/})
    \]
    to be the hypersheafification of the presheaf 
    \[
    S \mapsto \on{L} \Omega_{B(S)/A(S)}^{+, \bullet}  \in \on{cMod}_{\mathbb{D}_+^{\vee}}  (\on{DAlg} \on{Gr}D(\mathbb{Z})_{\on{ins}^0 A(S)/})
    \]
    on light profinite sets, where
    $\on{L} \Omega_{B(S)/A(S)}^{+, \bullet}$ 
    is the Hodge completed derived de Rham complex for derived rings as defined in \cite[definition 5.3.3.]{Raksit2026}.
\end{definition}

\begin{remark} \label{735082730854}
    More precisely, we define $  \on{L} \Omega_{B/A}^{+, \bullet} $
    as a sheaf on light profinite sets with values in 
    $\on{cMod}_{\mathbb{D}_+^{\vee}}  (\on{DAlg} \on{Gr}D(\mathbb{Z}))$, which is the hypersheafification of the presheaf 
    $ S \mapsto \on{L} \Omega_{B(S)/A(S)}^{+, \bullet}$, 
    where we forget the $\on{ins}^0 A(S)$-algebra structure of 
    $\on{L} \Omega_{B(S)/A(S)}^{+, \bullet}$. Moreover, if we forget the comodules structure and view $\on{L} \Omega_{B/A}^{+, \bullet} $ 
   as sheafification of  presheaves with values in $\on{DAlg} \on{Gr} D(\mathbb{Z})$, then it does not change the answer, since the forgetful functor preserves both limits and colimits as $\mathbb{D}_+^{\vee}$ is dualizable. 
   
    The $\on{ins}^0 A(S)$-algebra structure on $\on{L} \Omega_{B(S)/A(S)}^{+, \bullet}$ 
    defines a map of hypersheaves 
    $\on{ins}^0 A \longrightarrow   \on{L} \Omega_{B/A}^{+, \bullet}$, giving rise to 
    \[ \on{L} \Omega_{B/A}^{+, \bullet}  \in \on{cMod}_{\mathbb{D}_+^{\vee}}  (\on{DAlg}\on{Cond} \on{Gr}D(\mathbb{Z})_{\on{ins}^0 A/}).
    \]

    We have used implicitly here that 
    \[
    \on{cMod}_{\mathbb{D}_+^{\vee}}  (\on{DAlg}\on{Cond} \on{Gr}D(\mathbb{Z})_{\on{ins}^0 A/})
    \cong 
    (\on{cMod}_{\mathbb{D}_+^{\vee}}  \on{DAlg}\on{Cond} \on{Gr}D(\mathbb{Z}))_{\on{ins}^0 A/}
    \]
    and 
    \begin{equation} \label{80234803}
    \on{cMod}_{\mathbb{D}_+^{\vee}}  (\on{DAlg}\on{Gr}D(\mathbb{Z})_{\on{ins}^0 A(S)/})
    \cong 
    (\on{cMod}_{\mathbb{D}_+^{\vee}}  \on{DAlg}\on{Gr}D(\mathbb{Z}))_{\on{ins}^0 A(S)/}
    \end{equation}
    where on the right hand side we give the trivial $\mathbb{D}_+^{\vee}$-comodule structure on $\on{ins}^0 A$, namely the one induced by the unit map $\mathbb{Z} \longrightarrow \mathbb{D}^{\vee}_+$. This is simply a way of saying that $A$-linear $\mathbb{D}^{\vee}_+$-comodules on $A$-algebras are the same as $\mathbb{D}^{\vee}_+$-comodules of algebras over $A$ with trivial comodule structures. 
\end{remark}

\begin{lemma} \label{4702374903}
Let $A \longrightarrow B$ be a map of condensed derived rings, then after forgetting the $\mathbb{D}_+^{\vee}$-comodule structure 
 we have 
 \[
  \on{L} \Omega_{B/A}^{+, \bullet} \cong \on{Sym}^{\bullet}_{\on{ins}^0 B} (\Sigma \mathbb{L}_{B/A}[1])
 \]
 as condensed graded $B$-algebras. Moreover, the $\mathbb{D}_+^{\vee}$-comodule structure on the degree 0, 1 is given by the universal differential 
 $d : B \longrightarrow \mathbb{L}_{B/A}$. More precisely, under the equivalence between $\mathbb{D}_+^{\vee}$-comodules and $\mathbb{D}_+=\mathbb{Z}[\epsilon]/\epsilon^2 $-modules, the action of $\epsilon$ is given by $d$. 
\end{lemma}

\begin{proof}
    This is the hypersheafification of \cite[theorem 5.3.6.]{Raksit2026}. 
\end{proof}

We can also deduce the universal property of 
$\on{L} \Omega_{B/A}^{+, \bullet}$
from that of \cite[definition 5.3.3.]{Raksit2026}. Note that the main subtlety in the proof below is when $A$ is not constant sheaf, so that relevant adjunction is not about sheaf categories, but slice categories of sheaf categories, and we cannot waive our hands by just reducing to section of sheaves. We proceed by reducing to the case when $A$ is the absolute base, so that the relevant categories are functor categories, where we can apply the pointwise adjunction established in \cite[definition 5.3.3.]{Raksit2026}. 

\begin{proposition} \label{423087403284}
    Let $A \longrightarrow B$ be a map of condensed derived rings, then we have a canonical equivalence
    \[
    \on{Map}_{\on{cMod}_{\mathbb{D}_+^{\vee}}  (\on{DAlg}\on{Cond} \on{Gr}D(\mathbb{Z})_{\on{ins}^0 A/})} (\on{L} \Omega_{B/A}^{+, \bullet}, X^{\bullet} )
    \cong 
    \on{Map}_{\on{DAlg}\on{Cond} D(\mathbb{Z})_{ A/}} 
    (B, \on{ev}^0 X^{\bullet})
    \]
    for any $X^{\bullet} \in \on{cMod}_{\mathbb{D}_+^{\vee}}  (\on{DAlg}\on{Cond} \on{Gr}D(\mathbb{Z})_{\on{ins}^0 A/})$. 
\end{proposition}

\begin{proof}
    Let 
    $\on{L} \Omega_{B/A}^{+, \bullet, \on{pre}}$ 
    be the presheaf
     $ S \mapsto \on{L} \Omega_{B(S)/A(S)}^{+, \bullet}  $
    on light profintie sets with values in $\on{cMod}_{\mathbb{D}_+^{\vee}}  (\on{DAlg} \on{Gr}D(\mathbb{Z}))$
    as in remark \ref{735082730854}. Then $\on{L} \Omega_{B/A}^{+, \bullet} = \alpha_{\sharp} \on{L} \Omega_{B/A}^{+, \bullet, \on{pre}} $, where 
    \[
    \alpha_{\sharp}: 
    \on{Fun}(\on{Prof}^{\on{light}, \on{op}},  \on{cMod}_{\mathbb{D}_+^{\vee}}  (\on{DAlg} \on{Gr}D(\mathbb{Z}))) 
    \longrightarrow 
    \on{Cond} \on{cMod}_{\mathbb{D}_+^{\vee}}  (\on{DAlg} \on{Gr}D(\mathbb{Z}))
    \]
is the hypersheafification from the presheaf of light profinite sets with values in $\on{cMod}_{\mathbb{D}_+^{\vee}}  (\on{DAlg} \on{Gr}D(\mathbb{Z}))$
to hypersheaves, which is the left adjoint to the embedding 
\[
\on{Fun}(\on{Prof}^{\on{light}, \on{op}},  \on{cMod}_{\mathbb{D}_+^{\vee}}  (\on{DAlg} \on{Gr}D(\mathbb{Z}))) 
    \longleftarrow 
    \on{Cond} \on{cMod}_{\mathbb{D}_+^{\vee}}  (\on{DAlg} \on{Gr}D(\mathbb{Z}))
: \iota
\]
 The $A$-algebra structure on 
    $\on{L} \Omega_{B/A}^{+, \bullet, \on{pre}}$ gives rise to 
    \[
    \on{L} \Omega_{B/A}^{+, \bullet, \on{pre}} \in \on{Fun}(\on{Prof}^{\on{light}, \on{op}},  \on{cMod}_{\mathbb{D}_+^{\vee}}  (\on{DAlg} \on{Gr}D(\mathbb{Z})))_{ \iota \on{ins}^0 A /}
    \]
    where we use (\ref{80234803}), and gives $\on{ins}^0 A$ the trivial $\mathbb{D}_+^{\vee}$-comodule structure. We claim that 
    \[
\on{Map}_{\on{Fun}(\on{Prof}^{\on{light}, \on{op}},  \on{cMod}_{\mathbb{D}_+^{\vee}}  (\on{DAlg} \on{Gr}D(\mathbb{Z})))_{\iota\on{ins}^0 A/} }(\on{L} \Omega_{B/A}^{+, \bullet, \on{pre}}, \iota X^{\bullet} )
   \cong 
   \on{Map}_{\on{Fun}(\on{Prof}^{\on{light}, \on{op}},   \on{DAlg} D(\mathbb{Z}))_{\iota A/} } (\iota B, \on{ev}^0 \iota X^{\bullet}),
    \]
  from which  we can finish the proof by observing
    \begin{align*}
       \on{Map}_{\on{cMod}_{\mathbb{D}_+^{\vee}}  (\on{DAlg}\on{Cond} \on{Gr}D(\mathbb{Z})_{ \on{ins}^0 A/})} (\on{L} \Omega_{B/A}^{+, \bullet}, X^{\bullet} )
       &
    \cong 
    \on{Map}_{\on{Fun}(\on{Prof}^{\on{light}, \on{op}},  \on{cMod}_{\mathbb{D}_+^{\vee}}  (\on{DAlg} \on{Gr}D(\mathbb{Z})))_{\iota\on{ins}^0 A/} }(\on{L} \Omega_{B/A}^{+, \bullet, \on{pre}}, \iota X^{\bullet} )
    \\
    & 
    \cong 
    \on{Map}_{\on{Fun}(\on{Prof}^{\on{light}, \on{op}},   \on{DAlg} D(\mathbb{Z}))_{\iota A/} } (\iota B, \on{ev}^0 \iota X^{\bullet})
    \\
    & 
    \cong 
     \on{Map}_{\on{Fun}(\on{Prof}^{\on{light}, \on{op}},   \on{DAlg} D(\mathbb{Z}))_{\iota A/} } ( \iota B, \iota \on{ev}^0  X^{\bullet})
     \\
     & 
     \cong 
     \on{Map}_{\on{Cond}\on{DAlg} D(\mathbb{Z})_{ A/}  }
     (B, \on{ev}^0 X^{\bullet} )
    \end{align*}
    where we use $\iota$ for the inclusion of sheaves into presheaves with values in both $\on{DAlg} D(\mathbb{Z})$ and 
    $\on{cMod}_{\mathbb{D}_+^{\vee}}  (\on{DAlg} \on{Gr}D(\mathbb{Z}))) $.

To prove the claim, we first observe that if $A= \mathbb{Z}$, the initial object, then both sides are just mapping spaces of functor categories. Since functor categories preserves adjunction on the target, the claim in this case holds by the adjunction between $\on{DAlg} D(\mathbb{Z})$ and  $\on{cMod}_{\mathbb{D}_+^{\vee}}  (\on{DAlg} \on{Gr}D(\mathbb{Z}))) $, which is proved in \cite[proposition 5.3.2.]{Raksit2026}. 

Next, we observe that 
\begin{equation} \label{0935foe}
\on{L} \Omega_{B/A}^{+, \bullet, \on{pre}} 
\cong 
\on{L} \Omega_{B/\mathbb{Z}}^{+, \bullet, \on{pre}} \coprod_{\on{L} \Omega_{A/\mathbb{Z}}^{+, \bullet, \on{pre}}} \iota \on{ins}^0 A
\end{equation}
where $\on{ins}^0 A$ is equipped with the trivial $\mathbb{D}_{+}^{\vee}$-comodules structure, and the map 
$\on{L} \Omega_{A/\mathbb{Z}}^{+, \bullet, \on{pre}} \longrightarrow \iota \on{ins}^0 A$
is by adjunction the identity map $A = \on{ev}^0 \on{ins}^0 A $. On each section, it is the standard change of base formula of derived de Rham complex, which can be deduced directly from the universal property of \cite[definition 5.3.3.]{Raksit2026}, and it immediately passes to presheaves. 

The claim follows from  
\begin{align*}
&
    \on{Map}_{\on{Fun}(\on{Prof}^{\on{light}, \on{op}},  \on{cMod}_{\mathbb{D}_+^{\vee}}  (\on{DAlg} \on{Gr}D(\mathbb{Z})))_{\iota\on{ins}^0 A/} }(\on{L} \Omega_{B/A}^{+, \bullet, \on{pre}}, \iota X^{\bullet} ) 
   \\
   \cong
   &
   \on{Map}_{\on{Fun}(\on{Prof}^{\on{light}, \on{op}},  \on{cMod}_{\mathbb{D}_+^{\vee}}  (\on{DAlg} \on{Gr}D(\mathbb{Z})))_{\on{L} \Omega_{A/\mathbb{Z}}^{+, \bullet, \on{pre}} /} }(\on{L} \Omega_{B/\mathbb{Z
   }}^{+, \bullet, \on{pre}}, \iota X^{\bullet} ) 
  \\
 \cong & 
 \on{Map}_{\on{Fun}(\on{Prof}^{\on{light}, \on{op}},  \on{DAlg} D(\mathbb{Z}))_{\iota A /}}(\iota B, \on{ev}^0 \iota X^{\bullet} )
\end{align*}
where in the first equivalence we use (\ref{0935foe}), where the $\on{L} \Omega_{A/\mathbb{Z}}^{+, \bullet, \on{pre}} $-algebra structure on $\iota X^{\bullet} $ is the composition of its $A$-algebra structure with the map $\on{L} \Omega_{A/\mathbb{Z}}^{+, \bullet, \on{pre}} \longrightarrow A$. In the second equivalence, we use the known adjunction with absolute base $\mathbb{Z}$, and the formula $\on{Map_{C_{a/}}(X,Y)} \cong 
\on{Map}_{C}(X,Y) \times_{ \on{Map}_C(a, Y) } * $
for  mapping space of slice categories, where the left structure map is composition with structure map from $a$ to $X$, and the right map is the point corresponds to the structure map from $a$ to $Y$. 
\end{proof}

We can now define the Hodge completed derived de Rham complex for maps of analytic animated rings. 

\begin{definition} \label{zasdlfowei}
    Let $\mathcal{A} \longrightarrow \mathcal{B}$ be a map of analytic animated rings, define
    \[
    \widehat{\mathbb{DR}}_{\mathcal{B}/\mathcal{A}}^{\on{un}} = ( L\Omega_{\mathcal{B}^{\triangleright}/\mathcal{A}^{\triangleright}}^{+,\bullet}, D(L\Omega_{\mathcal{B}/\mathcal{A}}^{+,\bullet})) \in \on{GrAnDAG}^{\on{un}}
    \]
where the uncompleted graded analytic derived ring structure $D(L\Omega_{\mathcal{B}^{\triangleright}/\mathcal{A}^{\triangleright}}^{+,\bullet} )$
comes from the analytic ring structure $D(\mathcal{B})$ on $\on{ev}^0 L\Omega_{\mathcal{B}^{\triangleright}/\mathcal{A}^{\triangleright}}^{+,\bullet} = \mathcal{B}^{\triangleright}$ 
using proposition \ref{wury9387r9}. 

    Let 
    \[
    \widehat{\mathbb{DR}}_{\mathcal{B}/\mathcal{A}} := (L\Omega_{\mathcal{B}/\mathcal{A}}^{+,\bullet}, D(L\Omega_{\mathcal{B}/\mathcal{A}}^{+,\bullet}))  \in \on{GrAnDAG}
    \]
    be the completion of $ \widehat{\mathbb{DR}}_{\mathcal{B}/\mathcal{A}}^{\on{un}}$ 
    as defined in proposition \ref{wyrw3874}.
\end{definition}

From the construction of proposition \ref{wyrw3874}, we have 
 \begin{equation}  \label{zvdsueowr03}
    L\Omega_{\mathcal{B}/\mathcal{A}}^{+,\bullet} = L_{L\Omega_{\mathcal{B}^{\triangleright}/\mathcal{A}^{\triangleright}}^{+,\bullet}} (L\Omega_{\mathcal{B}^{\triangleright}/\mathcal{A}^{\triangleright}}^{+,\bullet} )  
    \end{equation}
where 
\[
L_{L\Omega_{\mathcal{B}^{\triangleright}/\mathcal{A}^{\triangleright}}^{+,\bullet}} :  \on{Mod}_{L\Omega_{\mathcal{B}^{\triangleright}/\mathcal{A}^{\triangleright}}^{+,\bullet}} \on{Cond} \on{Gr} D(\mathbb{Z})
\longrightarrow 
D(  L\Omega_{\mathcal{B}/\mathcal{A}}^{+,\bullet})
\]
is the left adjoint of the fully faithful embedding of $D(  L\Omega_{\mathcal{B}/\mathcal{A}}^{+,\bullet})$ 
into 
$\on{Mod}_{L\Omega_{\mathcal{B}^{\triangleright}/\mathcal{A}^{\triangleright}}^{+,\bullet}} \on{Cond} \on{Gr} D(\mathbb{Z})$, and 
\[
D(  L\Omega_{\mathcal{B}/\mathcal{A}}^{+,\bullet}) = \on{Mod}_{L\Omega_{\mathcal{B}/\mathcal{A}}^{+,\bullet}} \on{Gr} D(\mathcal{B}). 
\]
Moreover, using the last claim of proposition \ref{wury9387r9}, we have that 
\[
\on{ev}^i   L\Omega_{\mathcal{B}/\mathcal{A}}^{+,\bullet} \cong L_{\mathcal{B}} \circ  \on{ev}^i L\Omega_{\mathcal{B}^{\triangleright}/\mathcal{A}^{\triangleright}}^{+,\bullet} 
\cong 
L_B \on{Sym}^i_{\mathcal{B}^{\triangleright}} (\Sigma \mathbb{L}_{\mathcal{B}^{\triangleright}/\mathcal{A}^{\triangleright}} [1]) 
\cong 
\on{Sym}^i_{\mathcal{B}} (\Sigma \mathbb{L}_{\mathcal{B}/\mathcal{A}} [1]). 
\]

\begin{lemma}
Let $\mathcal{A} \longrightarrow \mathcal{B}$ be a map of analytic animated rings, then the comodule structure of $L\Omega_{\mathcal{B}^{\triangleright}/\mathcal{A}^{\triangleright}}^{+,\bullet} \in \on{cMod}_{\mathbb{D}_+^{\vee}}(\on{DAlg}\on{Cond} \on{Gr}D(\mathbb{Z})_{\on{ins}^0 A^{\vee}/})$ 
gives rise to
    \[
    \widehat{\mathbb{DR}}_{\mathcal{B}/\mathcal{A}}^{\on{un}}  \in \on{cMod}_{\underline{\mathbb{D}_+^{\vee}}} (\on{GrAnDAG}_{\on{ins}^0 \mathcal{A}/}^{\on{un}})
    \]
    where $\underline{\mathbb{D}_+^{\vee}}$ 
    is defined in example \ref{xfower870}, and $\on{ins}^0 \mathcal{A}$ is equipped with the graded uncompleted analytic ring structure from the analytic structure $D(\mathcal{A})$ on $\on{ev}^0 \on{ins}^0 \mathcal{A} = \mathcal{A}$ using proposition \ref{wury9387r9}.
\end{lemma}

\begin{proof}
We need to check all the coherent comodule structure maps of 
$ L\Omega_{\mathcal{B}^{\triangleright}/\mathcal{A}^{\triangleright}}^{+,\bullet} $
over $\underline{\mathbb{D}_+^{\vee}}$
lies in $ \on{FilAnDAG}^{\on{un}}$, but mapping spaces of $\on{GrAnDAG}^{\on{un}}$ by its very definition
are connected components of corresponding mapping spaces in $\on{DAlg}\on{Cond}\on{Gr}_{}D_{} (\mathbb{Z})$, in other words, it is enough to prove the claim by passing to $\pi_0$ of the categories. In this case, comodule structures are specified by the single comodule structure map 
    \begin{equation} \label{giyjo}
     L\Omega_{\mathcal{B}^{\triangleright}/\mathcal{A}^{\triangleright}}^{+,\bullet}
     \longrightarrow
\underline{\mathbb{D}_+^{\vee}}
\bigsqcup_{}  L\Omega_{\mathcal{B}^{\triangleright}/\mathcal{A}^{\triangleright}}^{+,\bullet}, 
    \end{equation}
    and we need to prove it
    is actually a map in 
$\on{FilAnDAG}^{\on{un}}$, in other words, we need to prove that the forgetful functor maps the complete modules on the right hand side to complete modules of the left hand side.

We observe that the analytic ring structure on the right hand side is also 
given by proposition \ref{wury9387r9}. Since $
\underline{\mathbb{D}_{+}^{\vee}}  $
is given the trivial analytic structure, by the description of colimits in $\on{GrAnDAG}^{\on{un}}$ as in lemma \ref{weuryeiwyrue}, we see that 
\[
D(
\underline{\mathbb{D}_{+}^{\vee}} \bigsqcup  L\Omega_{\mathcal{B}^{\triangleright}/\mathcal{A}^{\triangleright}}^{+,\bullet}) = \on{Mod}_{\mathbb{D}_{+}^{\vee}} D(L\Omega_{\mathcal{B}^{\triangleright}/\mathcal{A}^{\triangleright}}^{+,\bullet}) 
\]
where we view 
$D(L\Omega_{\mathcal{B}^{\triangleright}/\mathcal{A}^{\triangleright}}^{+,\bullet} )$
as tensored over $\on{Gr} \on{Cond}D(\mathbb{Z})$
through the obvious symmetric monoidal functor. We have that
\begin{align*}
\on{Mod}_{\mathbb{D}_{+}^{\vee}} D(L\Omega_{\mathcal{B}^{\triangleright}/\mathcal{A}^{\triangleright}}^{+,\bullet}) 
&
=
     \on{Mod}_{\mathbb{D}_{+}^{\vee}} \on{Mod}_{L\Omega_{\mathcal{B}^{\triangleright}/\mathcal{A}^{\triangleright}}^{+,\bullet}}  \on{Gr} D(\mathcal{B})
    \\
    &    
=
\on{Mod}_{\underline{\mathbb{D}_{+}^{\vee}} \otimes  L\Omega_{\mathcal{B}^{\triangleright}/\mathcal{A}^{\triangleright}}^{+,\bullet}} \on{Gr} D(\mathcal{B})
\\
&
=
\on{Mod}_{\underline{\mathbb{D}_{+}^{\vee}} \otimes L\Omega_{\mathcal{B}^{\triangleright}/\mathcal{A}^{\triangleright}}^{+,\bullet}}
\on{Mod}_{\on{ins}^0\on{ev}^0 (\underline{\mathbb{D}_{+}^{\vee}} \otimes  L\Omega_{\mathcal{B}^{\triangleright}/\mathcal{A}^{\triangleright}}^{+,\bullet}) } \on{Gr}  D(\mathcal{B})
\\
&
= \on{Mod}_{\underline{\mathbb{D}_{+}^{\vee}} \otimes L\Omega_{\mathcal{B}^{\triangleright}/\mathcal{A}^{\triangleright}}^{+,\bullet}}
\on{Gr}
\on{Mod}_{\on{ev}^0 (\underline{\mathbb{D}_{+}^{\vee}} \otimes  L\Omega_{\mathcal{B}^{\triangleright}/\mathcal{A}^{\triangleright}}^{+,\bullet}) }  D(\mathcal{B})
\end{align*}
where the first equality is proposition \ref{wury9387r9}. 
It means that the analytic structure on
$
\underline{\mathbb{D}_{+}^{\vee}} \bigsqcup  L\Omega_{\mathcal{B}^{\triangleright}/\mathcal{A}^{\triangleright}}^{+,\bullet}$
is induced by the analytic structure on 
$\on{ev}^0 (\underline{\mathbb{D}_{+}^{\vee}} \otimes  L\Omega_{\mathcal{B}^{\triangleright}/\mathcal{A}^{\triangleright}}^{+,\bullet}) $
through proposition \ref{wury9387r9}, where we give the analytic structure on 
$\on{ev}^0 (\underline{\mathbb{D}_{+}^{\vee}} \otimes  L\Omega_{\mathcal{B}^{\triangleright}/\mathcal{A}^{\triangleright}}^{+,\bullet})$
through the canonical map 
\[
\mathcal{B}^{\triangleright} = \on{ev}^0  L\Omega_{\mathcal{B}^{\triangleright}/\mathcal{A}^{\triangleright}}^{+,\bullet}
\longrightarrow
\on{ev}^0 (\underline{\mathbb{D}_{+}^{\vee}} \otimes  L\Omega_{\mathcal{B}^{\triangleright}/\mathcal{A}^{\triangleright}}^{+,\bullet}), 
\]
which is the identity map on the second factor. 

By lemma \ref{wityuoiewuroe}, it is enough to prove that the comodule structure map induces on $\on{ev}^0$ a map of analytic animated rings. Using $ \mathbb{D}_+^{\vee} = \mathbb{Z} \oplus  \Sigma^{-1}\mathbb{Z}[-1]$, we compute that 
\[
\on{ev}^0 (\underline{\mathbb{D}_{+}^{\vee}} \otimes  L\Omega_{\mathcal{B}^{\triangleright}/\mathcal{A}^{\triangleright}}^{+,\bullet})
= \mathcal{B}^{\triangleright} \oplus \mathbb{L}_{\mathcal{B}^{\triangleright}/\mathcal{A}^{\triangleright}}
\]
with the square zero algebra structure on $\mathbb{L}_{\mathcal{B}^{\triangleright}/\mathcal{A}^{\triangleright}}$, 
and the comodule structure map (\ref{giyjo}) on $\on{ev}^0$ is given by the canonical map 
\[
(\on{id},d): \mathcal{B}^{\triangleright} \longrightarrow \mathcal{B}^{\triangleright} \oplus \mathbb{L}_{\mathcal{B}^{\triangleright}/\mathcal{A}^{\triangleright}}
\]
where $d$ is universal map classifying universal square zero extensions, 
this is lemma \ref{4702374903}. The analytic structure is nilpotence invariant, so depends only on the quotient $\mathcal{B}^{\triangleright}$ of the square zero extension 
$\mathcal{B}^{\triangleright} \oplus \mathbb{L}_{\mathcal{B}^{\triangleright}/\mathcal{A}^{\triangleright}}$, and the composition of 
$(\on{id},d)$ with it is identity, so it induces a map of analytic animated rings, proving that the comodule structure map is a map in $\on{GrAnDAG}^{\on{un}}$. 
\end{proof}

\begin{corollary}
Let $\mathcal{A} \longrightarrow \mathcal{B}$ be a map of analytic animated rings,    then 
  \[
    \widehat{\mathbb{DR}}_{\mathcal{B}/\mathcal{A}} = (L\Omega_{\mathcal{B}/\mathcal{A}}^{+,\bullet}, D(L\Omega_{\mathcal{B}/\mathcal{A}}^{+,\bullet}))   \in \on{cMod}_{\underline{\mathbb{D}_+^{\vee}}} (\on{GrAnDAG}_{\on{ins}^0 \mathcal{A}/}), 
    \]
    and in particular
    \[
    L\Omega_{\mathcal{B}/\mathcal{A}}^{+,\bullet}\in \on{cMod}_{\mathbb{D}_+^{\vee}} \on{Gr} D(\mathcal{A})
    \]
\end{corollary}

\begin{proof}
    The is simply applying the left adjoint 
    $\on{GrAnDAG} 
    \longrightarrow \on{GrAnDAG}^{\on{un}}$
    of proposition \ref{wyrw3874} to 
    $ \widehat{\mathbb{DR}}_{\mathcal{B}/\mathcal{A}}^{\on{un}}$. It preserves comodule structures since they are given cocartesian symmetric monoidal structures. 
\end{proof}

We have defined $L\Omega_{\mathcal{B}/\mathcal{A}}^{+,\bullet}$
as a graded module with extra structures, but the usual derived de Rham complex is a filtered complex with Hodge filtration.  The underlying graded module of $L\Omega_{\mathcal{B}/\mathcal{A}}^{+,\bullet}$
is actually the (shifted) graded pieces of 
the sought-after derived de Rham complex with Hodge filtration. We can recover the original filtered complex from the graded pieces if we remember the $\mathbb{D}_-$-modules structure as in (\ref{nxcsfd}) up to completing the filtration, and this is why we have only the Hodge completed derived de Rham complex here. We need to first relate the $\mathbb{D}_+^{\vee}$-comodule structure to $\mathbb{D}_-$-module structure, and this is given by the equivalence
\[
\Sigma^{-2*}: \on{cMod}_{\mathbb{D}_+^{\vee}} \on{Gr} D(\mathcal{A}) \cong 
\on{cMod}_{\mathbb{D}_-^{\vee}} \on{Gr} D(\mathcal{A})
\]
as in \cite[remark 5.1.12.]{Raksit2026}, which sends a graded object $X^{\bullet}$ to $\Sigma^{-2\bullet} X^{\bullet}$. Composing with the equivalence
\[
\on{cMod}_{\mathbb{D}_-^{\vee}} \on{Gr} D(\mathcal{A}) 
\cong 
\on{Mod}_{\mathbb{D}_-} \on{Gr} D(\mathcal{A}) 
\cong 
\on{Fil}^{\wedge} D(\mathcal{A})
\]
in \cite[theorem 3.2.14.]{Raksit2026}, we have the equivalence 
 \begin{equation} \label{zcbksueife}
    |-|^{\geq *} : \on{cMod}_{\mathbb{D}_+^{\vee}}  \on{Gr} D(\mathcal{A})
    \cong 
    \on{Fil}^{\wedge} D(\mathcal{A})
    \end{equation}
    which is  defined in \cite[definition 5.2.4.]{Raksit2026}. We now define what we call the Hodge completed derived de Rham complex for analytic animated rings. We setup the theory using  graded objects with $\mathbb{D}_+^{\vee}$-comodule structure as in \cite{Raksit2026} since this is more easily comparable with the filtered Hochschild homology to be considered in the next section. 

    \begin{definition} \label{sfdsfjekl}
Let $\mathcal{A} \longrightarrow \mathcal{B}$ be a map of analytic animated rings, we define the Hodge completed derived de Rham complex of $\mathcal{B}$ over $\mathcal{A}$    to be 
\[
\widehat{\on{dR}}_{\mathcal{B}/\mathcal{A}}^{\geq \bullet} := |L\Omega_{\mathcal{B}/\mathcal{A}}^{+,\bullet} |^{\geq *} \in  \on{Fil}^{\wedge} D(\mathcal{A}). 
\]
    \end{definition}

From the definition and lemma \ref{4702374903}, we have 
\begin{equation}  \label{zdfhweuoofd}
\on{gr}^{\bullet}\widehat{\on{dR}}_{\mathcal{B}/\mathcal{A}}^{\geq \bullet} 
\cong 
\on{Und}(
L\Omega_{\mathcal{B}/\mathcal{A}}^{+,\bullet}) [-2*]
\cong
\Sigma^{-2*}
\on{Sym}^{\bullet}_{\mathcal{B} } (\Sigma \mathbb{L}_{\mathcal{B}/\mathcal{A}}[1])  
\end{equation}
where $\on{Und}$ denotes the underlying graded algebra. 

\begin{remark} \label{2307krieo}
    We note that $\widehat{\on{dR}}_{\mathcal{B}/\mathcal{A}} $ 
    is also a commutative algebra object in the filtered category since the equivalence $|-|^{\geq *}$ is symmetric monoidal. We can actually lift the graded derived algebra structure on $L\Omega_{\mathcal{B}/\mathcal{A}}^{+,\bullet}$ 
    to a derived algebra structure on $\widehat{\on{dR}}_{\mathcal{B}/\mathcal{A}}$, but the corresponding monad on filtered objects is not the usual one, but a suitably shifted version of the usual filtered derived algebra monad. 
\end{remark}

\begin{lemma} \label{2830948302}
Let $\mathcal{A} \longrightarrow \mathcal{B} \longrightarrow \mathcal{R}$ be  maps of analytic animated rings, there exists a canonical equivalence
\[
\widehat{\on{dR}}_{\mathcal{R}/\mathcal{A}}^{\geq \bullet} 
\coprod_{\widehat{\on{dR}}_{\mathcal{B}/\mathcal{A}}^{\geq \bullet} }
B^{\triangleright}
\cong
\widehat{\on{dR}}_{\mathcal{R}/\mathcal{B}}^{\geq \bullet}  
\]    
as elements of $\on{CAlg}  \on{Fil}^{\wedge} D(\mathcal{A})$, where we give $\mathcal{B}^{\triangleright}$ the trivial filtration, i.e. the only non-zero graded pieces is at degree 0, which has value $\mathcal{B}^{\triangleright}$. 
\end{lemma}

\begin{proof}
We only need to prove that 
    \[
    L\Omega_{\mathcal{R}/\mathcal{A}}^{+,\bullet} 
\coprod_{L\Omega_{\mathcal{B}/\mathcal{A}}^{+,\bullet} }
\mathcal{B}^{\triangleright}
\cong
L\Omega_{\mathcal{R}/\mathcal{B}}^{+,\bullet} 
    \]
    as elements in $\on{cMod}_{\mathbb{D}_+^{\vee}} \on{CAlg} \on{Gr}  D(\mathcal{A})$, where $\mathcal{B}^{\triangleright}$ 
    is concentrated at degree 0 with trivial $\mathbb{D}_{+}^{\vee}$-comodule structure.  We have the equivalence
     \[
    L\Omega_{\mathcal{R}^{\triangleright}/\mathcal{A}^{\triangleright}}^{+,\bullet} 
\coprod_{L\Omega_{\mathcal{B}^{\triangleright}/\mathcal{A}^{\triangleright}}^{+,\bullet} }
\mathcal{B}^{\triangleright}
\cong
L\Omega_{\mathcal{R}^{\triangleright}/\mathcal{B}^{\triangleright}}^{+,\bullet}, 
    \]
    which is the usual base change formula for Hodge completed derived de Rham complex of algebras and sheafifying. 
    
    We observe that this is also a pushout diagram in  $\on{GrAnDAG}^{\on{un}}$, where we give $\mathcal{B}^{\triangleright}$ the graded analytic structure induced from the analytic structure of $\mathcal{B}$ using proposition \ref{wury9387r9}, $L\Omega_{\mathcal{R}^{\triangleright}/\mathcal{A}^{\triangleright}}^{+,\bullet} $
    is given
    the analytic structure $\widehat{\mathbb{DR}}_{\mathcal{R}/\mathcal{A}}^{\on{un}}$ in definition \ref{zasdlfowei}, and similarly for 
    $L\Omega_{\mathcal{R}^{\triangleright}/\mathcal{B}^{\triangleright}}^{+,\bullet}$ and 
    $L\Omega_{\mathcal{B}^{\triangleright}/\mathcal{A}^{\triangleright}}^{+,\bullet} $. 
    This is clear from the description of colimits in $\on{GrAnDAG}^{\on{un}}$ in lemma \ref{weuryeiwyrue}, from which we see that a module is complete with respect to the analytic structure of the left hand side if and only if it lies in both $\on{Gr} D(\mathcal{R}) $ and $\on{Gr} D(\mathcal{B})$ under the relevant forgetful maps. This is equivalent to requiring it lies in $\on{Gr} D(\mathcal{R})$ alone since the $\mathcal{B}^{\triangleright}$-module comes from the $\mathcal{R}^{\triangleright}$-module structure under the map of analytic rings $\mathcal{B} \rightarrow \mathcal{R}$, which is exactly the analtic structure on the right hand side. Now applying the completion functor from uncompleted graded analytic rings to completed graded analytic rings in proposition \ref{wyrw3874}, and taking the underlying condensed graded rings gives the desired formula. It is also a pushout in the $\mathbb{D}^{\vee}_+$-comodule category since the   natural map from the left to right is clearly compatible with  $\mathbb{D}^{\vee}_+$-comodule  structure, and it is then enough to check it is an equivalence after forgetting the comodule structure. 
\end{proof}

\section{Filtered Hochschild Homology of analytic animated rings}

Let $A \in \on{DAlg}\on{Cond}\on{Fil}D_{} (\mathbb{Z})$, 
then we map 
$\underline{\mathbb{T}_{\on{Fil}}^{\vee}}$ 
to an element 
\[
A[\underline{\mathbb{T}_{\on{Fil}}^{\vee}}] := A \otimes_{\underline{\mathbb{Z}}} \underline{\mathbb{T}_{\on{Fil}}^{\vee}}
\in
\on{cCAlg}(
\on{DAlg}\on{Cond}\on{Fil}D_{} (\mathbb{Z})_{A/})
\]
through the left adjoint 
\[
- \amalg A:
\on{DAlg}\on{Cond}\on{Fil}D_{} (\mathbb{Z})_{}
\longrightarrow
\on{DAlg}\on{Cond}\on{Fil}D_{} (\mathbb{Z})_{A/}
\]
to the forgetful functor, which we see that it is also a coalgebra object since the left adjoint is symmetric monoidal when we give both sides the cocartesian product.
 By \cite[construction 4.2.32]{Raksit2026}, the monad $U \circ \on{Sym}^{\bullet}_A$ 
 on 
 $\on{Mod}_{A} \on{Cond}\on{Fil}D_{} (\mathbb{Z})$
 as defined in lemma \ref{mxcnaaadfs}
 induces a monad $T_A$ on $\on{cMod}_{A[\underline{\mathbb{T}_{\on{Fil}}^{\vee}}]} \on{Mod}_{A} \on{Cond}\on{Fil}D_{} (\mathbb{Z})$
 such that
 \begin{equation} \label{lseoir}
 \on{LMod}_{T_A} (\on{cMod}_{A[\underline{\mathbb{T}_{\on{Fil}}^{\vee}}]} \on{Mod}_{A} \on{Cond}\on{Fil}D_{} (\mathbb{Z})) \cong \on{cMod}_{A[\underline{\mathbb{T}_{\on{Fil}}^{\vee}}]} 
( \on{DAlg}\on{Cond}\on{Fil}D_{} (\mathbb{Z})_{A/})
 \end{equation}
 we have the commutative diagram
 \begin{equation} \label{zddgsged}
 \begin{tikzcd}
 \on{cMod}_{A[\underline{\mathbb{T}_{\on{Fil}}^{\vee}}]} \on{Mod}_{A} \on{Cond}\on{Fil}D_{} (\mathbb{Z})
 \arrow[r,"T_A"]
 \arrow[d]
 &
 \on{cMod}_{A[\underline{\mathbb{T}_{\on{Fil}}^{\vee}}]} \on{Mod}_{A} \on{Cond}\on{Fil}D_{} (\mathbb{Z})
 \arrow[d]
\\
\on{Mod}_{A} \on{Cond}\on{Fil}D_{} (\mathbb{Z})
\arrow[r,"U \circ \on{Sym}^{\bullet}_A"]
&
\on{Mod}_{A} \on{Cond}\on{Fil}D_{} (\mathbb{Z})
 \end{tikzcd}
 \end{equation}
where the vertical arrows are forgetful functors. Indeed, the monadic adjunction
\[
\on{Sym}_A^{\bullet}:
\on{Mod}_{A} \on{Cond}\on{Fil}D_{} (\mathbb{Z}) 
\leftrightarrows
\on{DAlg}\on{Cond}\on{Fil}D_{} (\mathbb{Z})_{A/}: U
\]
in lemma \ref{mxcnaaadfs}
has symmetric monoidal right adjoint,  where we give $\on{DAlg}\on{Cond}\on{Fil}D_{} (\mathbb{Z})_{A/}$ the cocartesian product. 
This holds since by lemma \ref{areyurwo}
the coproducts in $\on{DAlg}\on{Cond}\on{Fil}D_{} (\mathbb{Z})_{A/}$ 
is given by tensor product over $A$ which is the same as the tensor product on the left hand side. Then the left adjoint $\on{Sym}_A^{\bullet}$ 
is oplax symmetric monoidal, so indcues a monadic adjunction on the comodule categories 
\[
\on{Sym}_A^{\bullet}:
\on{cMod}_{A[\underline{\mathbb{T}_{\on{Fil}}^{\vee}}]}
(\on{Mod}_{A} \on{Cond}\on{Fil}D_{} (\mathbb{Z}) )
\leftrightarrows
\on{cMod}_{A[\underline{\mathbb{T}_{\on{Fil}}^{\vee}}]}
(\on{DAlg}\on{Cond}\on{Fil}D_{} (\mathbb{Z})_{A/}): U. 
\]

Let $A \in \on{DAlg}\on{Cond}D_{} (\mathbb{Z})$, i.e. $A$ is not a filtered condensed derived ring, but only a condensed drived ring, then 
we have a commutative diagram 
\begin{equation} \label{q3473649374}
\begin{tikzcd} 
\on{Mod}_{A} \on{Cond}D_{} (\mathbb{Z}) 
& 
\on{DAlg}\on{Cond}D_{} (\mathbb{Z})_{A/}
\arrow[l,"U"]
\\
\on{Mod}_{\on{ins}^0A} \on{Cond}\on{Fil}D_{} (\mathbb{Z}) 
\arrow[u,"\on{ev}^0"]
    & 
\on{DAlg}\on{Cond}\on{Fil}D_{} (\mathbb{Z})_{\on{ins}^0A/}
\arrow[u,"\on{ev}^0"]
\arrow[l,"U"]
    \\
\on{cMod}_{\on{ins}^0A[\underline{\mathbb{T}_{\on{Fil}}^{\vee}}]}
(\on{Mod}_{\on{ins}^0A} \on{Cond}\on{Fil}D_{} (\mathbb{Z}) )
\arrow[u,"\beta'"]
    & \on{cMod}_{\on{ins}^0A[\underline{\mathbb{T}_{\on{Fil}}^{\vee}}]}
(\on{DAlg}\on{Cond}\on{Fil}D_{} (\mathbb{Z})_{\on{ins}^0A/})
\arrow[u,"\beta"]
\arrow[l,"U"]
\end{tikzcd}
\end{equation}
where 
\[
\on{ins}^0: \on{Cond}D_{} (\mathbb{Z})
\longrightarrow
\on{Cond}\on{Fil}D_{} (\mathbb{Z})
\]
is the left adjoint of the evaluation 
$\on{ev}^0$, which induces
\[
\on{ins}^0: \on{DAlg}\on{Cond}D_{} (\mathbb{Z})
\longrightarrow
\on{DAlg}\on{Cond}\on{Fil}D_{} (\mathbb{Z})
\]
by \cite[remark 4.3.6.]{Raksit2026} and passing to condensed objects. The first vertical arrows are induced by the unit map
$A \longrightarrow \on{ev}^0 \on{ins}^0 A$, 
and the lower vertical arrows are  forgetting the comodule structure. The commutativity of the diagram is clear by what we have already seen. We claim that it is also a commutative diagram of right adjoint, this is clear except for the two lower vertical arrows, since the forgetful functor from a comodule category to the underlying category is in general only a left adjoint. We argue instead using that $\mathbb{T}_{\on{Fil}}^{\vee}$ 
is dualizable in $\on{Fil}D(\mathbb{Z})$, so 
$\on{ins}^0A[\underline{\mathbb{T}_{\on{Fil}}^{\vee}}]$ 
is dualizable in 
$\on{Mod}_{\on{ins}^0A} \on{Cond}\on{Fil}D_{} (\mathbb{Z})$
since the corresponding functors are symmetric monoidal. We know from \cite[corollary 2.3.3.]{Raksit2026}
that 
\begin{equation} \label{sfddfe}
\on{cMod}_{\on{ins}^0A[\underline{\mathbb{T}_{\on{Fil}}^{\vee}}]}
(\on{Mod}_{\on{ins}^0A} \on{Cond}\on{Fil}D_{} (\mathbb{Z}) ) 
\cong 
\on{Mod}_{\on{ins}^0A[\underline{\mathbb{T}_{\on{Fil}}}]}
(\on{Mod}_{\on{ins}^0A} \on{Cond}\on{Fil}D_{} (\mathbb{Z}) )
\end{equation}
where $\mathbb{T}_{\on{Fil}}:= \tau_{\geq *} \mathbb{Z}[\mathbb{S}^1]$ is the dual of 
$\mathbb{T}_{\on{Fil}}^{\vee}$, which is also a cocommutative commutative bialgebra in $\on{Fil}D(\mathbb{Z})$. Under this equivalence, the forgetful functor $\beta'$ of comodules corresponds to forgetful functor of modules, and is clearly a right adjoint. This also implies the veritcal arrow
 $\beta$ is a right adjoint. Indeed, it is enough to prove that $\beta$ preserves limits, but then  $U \circ \beta$ preserves limits since the commutative diagram tells us that this is $\beta' \circ U$ which we have already seen is a right adjoint. Now the claim follows from the conservativity of $U$ since the horizontal arrows are monadic right adjoints.

Taking the left adjoint the commutative diagram, we have the commutative diagram 
\begin{equation} \label{zcnbdbfd}
\begin{tikzcd} 
\on{Mod}_{A} \on{Cond}D_{} (\mathbb{Z}) 
\arrow[r,"\on{Sym}_A^{\bullet}"]
\arrow[d,"\on{ins}^0"]
& 
\on{DAlg}\on{Cond}D_{} (\mathbb{Z})_{A/}
\arrow[d,"\on{ins}^0"]
\arrow[dd,bend left=80,"\on{HH}_{\on{Fil}}(-/A)"]
\\
\on{Mod}_{\on{ins}^0A} \on{Cond}\on{Fil}D_{} (\mathbb{Z}) 
\arrow[r,"\on{Sym}_{\on{ins}^0A}^{\bullet}"]
\arrow[d," -\otimes_{ \on{ins}^0A} \on{ins}^0A\lbrack \underline{\mathbb{T}_{\on{Fil}}}\rbrack"]
    & 
\on{DAlg}\on{Cond}\on{Fil}D_{} (\mathbb{Z})_{\on{ins}^0A/}
\arrow[d,"\alpha"]
    \\
\on{cMod}_{\on{ins}^0A[\underline{\mathbb{T}_{\on{Fil}}^{\vee}}]}
(\on{Mod}_{\on{ins}^0A} \on{Cond}\on{Fil}D_{} (\mathbb{Z}) )
\arrow[r,"\on{Sym}_{\on{ins}^0A}^{\bullet}"]
    & \on{cMod}_{\on{ins}^0A[\underline{\mathbb{T}_{\on{Fil}}^{\vee}}]}
(\on{DAlg}\on{Cond}\on{Fil}D_{} (\mathbb{Z})_{\on{ins}^0A/})
\end{tikzcd}
\end{equation}
where  we use the equivalence (\ref{sfddfe}) in the lower left vertical arrow. 

Similar to \cite[definition 6.2.5.]{Raksit2026}, we make the following definition of filtered Hochschild homology. 

\begin{definition}
Let
$A \in \on{DAlg}\on{Cond}\on{Fil}D_{} (\mathbb{Z})$, 
we define the filtered Hochschild homology over $A$ to be the functor 
\[
\on{HH}_{\on{Fil}}(-/A): 
\on{DAlg}\on{Cond}D_{} (\mathbb{Z})_{A/}
\overset{\alpha \circ \on{ins}^0}{\longrightarrow} 
\on{cMod}_{\on{ins}^0A[\underline{\mathbb{T}_{\on{Fil}}^{\vee}}]}
(\on{DAlg}\on{Cond}\on{Fil}D_{} (\mathbb{Z})_{\on{ins}^0A/}).
\]
\end{definition}

\begin{remark}
    Compared to \cite[definition 6.2.5.]{Raksit2026}, we do not restrict to non-negatively filtered derived rings, this does not affect the result since we will see  below that it is automatic non-negatively filtered. We can also see that directly from the definition as follows. We   observe the fully faithful embedding
    \[
    \on{ins}^{\geq 0}:\on{Cond} \on{Fil}^{\geq 0} D(\mathbb{Z}) \longrightarrow \on{Cond}\on{Fil} D(\mathbb{Z}) 
    \]
    induces 
    \[
    \on{ins}^{\geq 0}: 
    \on{DAlg}\on{Cond} \on{Fil}^{\geq 0} D(\mathbb{Z}) 
    \longrightarrow 
    \on{DAlg}\on{Cond} \on{Fil}^{} D(\mathbb{Z})
    \]
    by \cite[remark 4.3.6.]{Raksit2026} and passing to condensed objects, 
    which is the left adjoint to the $\on{ev}^{\geq 0}$, the restriction to non-negative filtration.
    It induces 
    \[
    \on{ins}^{\geq 0}: 
    \on{cMod}_{\on{ins}^0A[\underline{\mathbb{T}_{\on{Fil}}^{\vee}}]}
(\on{DAlg}\on{Cond}\on{Fil}^{\geq 0}D_{} (\mathbb{Z})_{\on{ins}^0A/})
\longrightarrow
\on{cMod}_{\on{ins}^0A[\underline{\mathbb{T}_{\on{Fil}}^{\vee}}]}
(\on{DAlg}\on{Cond}\on{Fil}^{}D_{} (\mathbb{Z})_{\on{ins}^0A/}).
    \]
where we 
 follow \cite[notation 6.2.3.]{Raksit2026} to define the left hand side by the pullback diagram
    \[
    \begin{tikzcd}
    \on{cMod}_{\on{ins}^0A[\underline{\mathbb{T}_{\on{Fil}}^{\vee}}]}
(\on{DAlg}\on{Cond}\on{Fil}^{\geq 0}D_{} (\mathbb{Z})_{\on{ins}^0A/})
\arrow[r,hook,"\on{ins}^{\geq 0}"]
\arrow[d,"\beta^{\geq 0}"]
\arrow[dr, phantom, "\ulcorner", very near start]
&
        \on{cMod}_{\on{ins}^0A[\underline{\mathbb{T}_{\on{Fil}}^{\vee}}]}
(\on{DAlg}\on{Cond}\on{Fil}D_{} (\mathbb{Z})_{\on{ins}^0A/})
\arrow[d,"\beta"]
        \\
  \on{DAlg}\on{Cond}\on{Fil}^{\geq 0}D_{} (\mathbb{Z})_{\on{ins}^0A/}
  \arrow[r,hook, "\on{ins}^{\geq 0}"]
  &
  \on{DAlg}\on{Cond}\on{Fil}D_{} (\mathbb{Z})_{\on{ins}^0A/}
    \end{tikzcd}
    \]
    with the right vertical arrow forgetful functor. Note that we cannot directly define it as the comodule category in $ \on{DAlg}\on{Cond}\on{Fil}^{\geq 0}D_{} (\mathbb{Z})_{\on{ins}^0A/}$
    since $\underline{\mathbb{T}_{\on{Fil}}^{\vee}}$
    is negatively filtered, so not in this category. 

    We claim that it also induces a commutative diagram
    \[
    \begin{tikzcd}
    \on{cMod}_{\on{ins}^0A[\underline{\mathbb{T}_{\on{Fil}}^{\vee}}]}
(\on{DAlg}\on{Cond}\on{Fil}^{\geq 0}D_{} (\mathbb{Z})_{\on{ins}^0A/})
\arrow[d,"\beta^{\geq 0}"]
&
        \on{cMod}_{\on{ins}^0A[\underline{\mathbb{T}_{\on{Fil}}^{\vee}}]}
(\on{DAlg}\on{Cond}\on{Fil}D_{} (\mathbb{Z})_{\on{ins}^0A/})
\arrow[d,"\beta"]
\arrow[l,"\on{ev}^{\geq 0}"]
        \\
  \on{DAlg}\on{Cond}\on{Fil}^{\geq 0}D_{} (\mathbb{Z})_{\on{ins}^0A/}
  &
  \on{DAlg}\on{Cond}\on{Fil}D_{} (\mathbb{Z})_{\on{ins}^0A/}
  \arrow[l,"\on{ev}^{\geq 0}"]
    \end{tikzcd}
    \]
    this is because (\ref{lseoir}) and (\ref{zddgsged}) and the observation that the monad $U \circ \on{Sym}_A^{\bullet}$ preserves non-negatively filtered objects when $A$ is so, so  
    $\on{cMod}_{\on{ins}^0A[\underline{\mathbb{T}_{\on{Fil}}^{\vee}}]}
(\on{DAlg}\on{Cond}\on{Fil}^{\geq 0}D_{} (\mathbb{Z})_{\on{ins}^0A/})$
is the module category of the restriction of the monad $T_{\on{ins}^0A}$ on $\on{cMod}_{\on{ins}^0A[\underline{\mathbb{T}_{\on{Fil}}^{\vee}}]}\on{Mod}_{\on{ins^0A}}\on{Cond}\on{Fil}^{\geq 0}D_{} (\mathbb{Z})$, and the commutative diagram follows. 

Now we have the commutative diagram
    \[
    \begin{tikzcd}
    \on{cMod}_{\on{ins}^0A[\underline{\mathbb{T}_{\on{Fil}}^{\vee}}]}
(\on{DAlg}\on{Cond}\on{Fil}^{\geq 0}D_{} (\mathbb{Z})_{\on{ins}^0A/})
\arrow[d,"\beta^{\geq 0}"]
&
        \on{cMod}_{\on{ins}^0A[\underline{\mathbb{T}_{\on{Fil}}^{\vee}}]}
(\on{DAlg}\on{Cond}\on{Fil}D_{} (\mathbb{Z})_{\on{ins}^0A/})
\arrow[d,"\beta"]
\arrow[l,"\on{ev}^{\geq 0}"]
        \\
  \on{DAlg}\on{Cond}\on{Fil}^{\geq 0}D_{} (\mathbb{Z})_{\on{ins}^0A/}
  \arrow[d,"\on{ev}^0"]
  &
  \on{DAlg}\on{Cond}\on{Fil}D_{} (\mathbb{Z})_{\on{ins}^0A/}
  \arrow[l,"\on{ev}^{\geq 0}"]
  \arrow[d,"\on{ev}^0"]
  \\
\on{DAlg}\on{Cond}D_{} (\mathbb{Z})_{A/}
\arrow[r,equal] 
&  \on{DAlg}\on{Cond}D_{} (\mathbb{Z})_{A/}
    \end{tikzcd}
    \] 
    taking left adjoint provides us with the commutative diagram
       \[
    \begin{tikzcd}
    \on{cMod}_{\on{ins}^0A[\underline{\mathbb{T}_{\on{Fil}}^{\vee}}]}
(\on{DAlg}\on{Cond}\on{Fil}^{\geq 0}D_{} (\mathbb{Z})_{\on{ins}^0A/})
\arrow[r,"\on{ins}^{\geq 0}"]
&
        \on{cMod}_{\on{ins}^0A[\underline{\mathbb{T}_{\on{Fil}}^{\vee}}]}
(\on{DAlg}\on{Cond}\on{Fil}D_{} (\mathbb{Z})_{\on{ins}^0A/})
        \\
  \on{DAlg}\on{Cond}\on{Fil}^{\geq 0}D_{} (\mathbb{Z})_{\on{ins}^0A/}
  \arrow[u,"\alpha^{\geq 0}"]
  \arrow[r,"\on{ins}^{\geq 0}"]
  &
  \on{DAlg}\on{Cond}\on{Fil}D_{} (\mathbb{Z})_{\on{ins}^0A/}
  \arrow[u,"\alpha"]
  \\
\on{DAlg}\on{Cond}D_{} (\mathbb{Z})_{A/}
\arrow[r,equal] 
\arrow[u,"\on{ins}^0"]
\arrow[uu,bend left=80,"\on{HH}^{\geq 0}_{\on{Fil}}(-/A)"]
&  \on{DAlg}\on{Cond}D_{} (\mathbb{Z})_{A/}
\arrow[u,"\on{ins}^0"]
\arrow[uu,bend right=80,"\on{HH}_{\on{Fil}}(-/A)"']
    \end{tikzcd}
    \] 
   where we denote $\on{HH}^{\geq 0}_{\on{Fil}}(-/A)$ to be the composite $\alpha^{\geq 0} \circ \on{ins}^0$, which is exactly how the filtered Hochschild homology in \cite[definition 6.2.5.]{Raksit2026} is defined.  The commutative diagram tells us that 
   \[
\on{HH}_{\on{Fil}}(-/A)) = \on{ins}^{\geq 0} \circ \on{HH}_{\on{Fil}}^{\geq 0}(-/A),
\]
which is the desired comparison. 
\end{remark}

We check it has the expected properties. 

\begin{lemma}
    The functor 
    $\on{HH}_{\on{Fil}}(-/A)$ commutes with colimits. 
\end{lemma}

\begin{proof}
    It is defined by composite of left adjoint functors. 
\end{proof}

\begin{lemma} \label{70237450834}
For any $B \in  \on{DAlg}\on{Cond}D_{} (\mathbb{Z})_{A/}$,
    the filtration on $\on{HH}_{\on{Fil}}(B/A)$ is (non-equivariantly) $B$-linear. More precisely, there is a natural map 
    \[
    \on{ins}^0 B \longrightarrow \beta(\on{HH}_{\on{Fil}}(B/A))
    \]
    in 
    $\on{DAlg}\on{Cond}\on{Fil}D_{} (\mathbb{Z})_{\on{ins}^0A/}$,
    where $\beta$ is forgetting the comodule structure functor as in (\ref{q3473649374}). 
\end{lemma}

\begin{proof}
    This is simply the unit of the adjunction $\alpha$ and $\beta$, i.e. 
    $\on{ins}^0 \longrightarrow \beta \circ \alpha \circ \on{ins}^0 = \beta \circ \on{HH}_{\on{Fil}}(-/A) $. 
\end{proof}

We can now prove the Hochschild-Kostant-Rosenberg for condensed animated rings, without the analytic structure. 

\begin{proposition}[Hochschild-Kostant-Rosenberg for condensed animated rings] \label{3074375372}
    Let $A \in \on{DAlg} \on{Cond} D(\mathbb{Z})$, and we assume that $A$ is connective, in other words, $A$ is a condensed animated ring, then for any condensed animated ring $B$ over $A$, we have 
    \[
    \on{gr}^{\bullet}\on{HH}_{\on{Fil}}(B/A) 
    \cong
    \on{L} \Omega_{B/A}^{+, \bullet} \in 
   \on{cMod}_{\underline{\mathbb{D}_{+}^{\vee}}}
(\on{DAlg}\on{Cond}\on{Gr}D_{} (\mathbb{Z})_{\on{ins}^0A/}). 
    \]

    Moreover, there exists an equivalence
    \[
    \on{colim}(\on{HH}_{\on{Fil}}(B/A) ) \cong B \otimes_{B\otimes_AB}B \in (\on{DAlg} \on{Cond} D(\mathbb{Z})_{A/})^{B\mathbb{S}^1}
    \]
    the usual Hochschild homology, where 
    \[
    \on{colim}: \on{Cond}\on{Fil} D(\mathbb{Z}) \longrightarrow 
    \on{Cond}D(\mathbb{Z}) 
    \]
    is the colimit of the filtration, and it induces 
    \[
    \on{colim}: \on{DAlg} \on{Cond}\on{Fil} D(\mathbb{Z})
    \longrightarrow 
    \on{DAlg} \on{Cond}D(\mathbb{Z})
    \]
    by \cite[remark 4.3.6.]{Raksit2026}
    and passing to condensed objects, which is colimit preserving by \cite[remark 4.2.25.]{Raksit2026}, so it induces
    \[
    \on{colim}: 
    \on{cMod}_{\on{ins}^0A[\underline{\mathbb{T}_{\on{Fil}}^{\vee}}]}
(\on{DAlg}\on{Cond}\on{Fil}D_{} (\mathbb{Z})_{\on{ins}^0A/}) 
\longrightarrow 
\on{cMod}_{A^{\mathbb{S}^1}}
(\on{DAlg}\on{Cond}D_{} (\mathbb{Z})_{A/})
    \]
    where $A^{\mathbb{S}^1}:= \underline{\mathbb{Z}^{\mathbb{S}^1}} \bigsqcup A \in \on{cCAlg}(\on{DAlg}\on{Cond}D_{} (\mathbb{Z})_{A/})$. 
\end{proposition}

\begin{proof}
This is exactly the same as classical computation of Hochschild homology. 
    We know that 
    $\on{Mod}_{A} \on{Cond}D_{} (\mathbb{Z})$
    is generated by $A[S]$ for $S$ light profinite sets, then condensed animated rings over $A$ is generated by $\on{Sym}^{\bullet}_A(A[S])$, and by the commutative diagram (\ref{zcnbdbfd}), we have
    \[
     \on{HH}_{\on{Fil}}(\on{Sym}^{\bullet}_A(A[S])/A) 
     \cong 
     \on{Sym}^{\bullet}_{\on{ins}^0A}(\on{ins}^0A[S]
    \otimes_{ \on{ins}^0A} \on{ins}^0A\lbrack \underline{\mathbb{T}_{\on{Fil}}}\rbrack ).
    \]
    We can compute the graded pieces
    \begin{align*}
   \on{gr}^{\bullet} \on{HH}_{\on{Fil}}(\on{Sym}^{\bullet}_A(A[S])/A) 
    & 
    \cong 
   \on{gr}^{\bullet} \on{Sym}^{\bullet}_{\on{ins}^0A}(\on{ins}^0A[S]
    \otimes_{ \on{ins}^0A} \on{ins}^0A\lbrack \underline{\mathbb{T}_{\on{Fil}}}\rbrack )
    \\
    &
    \cong 
    \on{Sym}^{\bullet}_{\on{ins}^0A}(\on{ins}^0A[S] \oplus \Sigma \on{ins}^0A[S][1])
    \\
    & 
    \cong  \on{Sym}^{\bullet}_{\on{ins}^0A}(\on{ins}^0A[S] ) \otimes_{\on{ins}^0A}
    \on{Sym}^{\bullet}_{\on{ins}^0A}(\Sigma \on{ins}^0A[S][1])
    \\
    & 
    \cong 
    \on{Sym}^{\bullet}_{\on{ins}^0\on{Sym}^{\bullet}_{A}(A[S] )}(\Sigma \on{ins}^0A[S] \otimes_A \on{Sym}^{\bullet}_{A}(A[S] ) [1]
    )
    \\
    & 
    \cong \on{Sym}^{\bullet}_{\on{ins}^0\on{Sym}^{\bullet}_{A}(A[S] )}(\Sigma \mathbb{L}_{\on{Sym}^{\bullet}_{A}(A[S] )/A}[1])
     \end{align*}
where we use that $\on{gr}^{\bullet} \mathbb{T}_{\on{Fil}} = \mathbb{Z} \oplus \Sigma\mathbb{Z}[1]$, in the second equivalences, and we have used the same notation of $\on{ins}^0$ denoting both the functor from condensed modules to graded and filtered condensed modules, and it should be clear that from the second equivalence, we are working in the graded setting. We have also used the standard base change formula for cotangent complex
\[
 \mathbb{L}_{\on{Sym}_{A}(A[S])/A} = A[S] \otimes_{A} \on{Sym}_{A}(A[S]). 
 \]
 Note that this is an equivalence of $\on{gr}^{\bullet} \mathbb{T}_{\on{Fil}}=\mathbb{D}_+=\mathbb{Z}[\epsilon]/\epsilon^2$-modules, where $\epsilon$ acts by the universal differential $d : \on{Sym}_{A}(A[S] )\longrightarrow \mathbb{L}_{\on{Sym}_{A}(A[S])/A}$.

 Using proposition \ref{423087403284}, we see that the identity map on $\on{Sym}_{A}(A[S])$ 
 induces a map 
 \[
 \on{L} \Omega_{B/A}^{+, \bullet} 
 \longrightarrow 
 \on{gr}^{\bullet} \on{HH}_{\on{Fil}}(\on{Sym}^{\bullet}_A(A[S])/A)
 \]
in
$   \on{cMod}_{\underline{\mathbb{D}_{+}^{\vee}}}
(\on{DAlg}\on{Cond}\on{Gr}D_{} (\mathbb{Z})_{\on{ins}^0A/}). 
$
It follows from by lemma \ref{4702374903} that both sides have the same underlying graded algebra 
$\on{Sym}^{\bullet}_{\on{ins}^0\on{Sym}^{\bullet}_{A}(A[S] )}(\Sigma \mathbb{L}_{\on{Sym}^{\bullet}_{A}(A[S] )/A}[1])$, and by construction we know it is the identity map at degree 0. Since they are symmetric algebra on degree 1 elements, to prove it is an isomorphism it is enough to prove that it is an isomorphism in degree 1. Truncating both sides in degree 0 and 1, we are left with a commutative diagram
\[
\begin{tikzcd}
 \mathbb{L}_{\on{Sym}_{A}(A[S])/A}
 \arrow[r]
 & 
 \mathbb{L}_{\on{Sym}_{A}(A[S])/A} 
    \\
  \on{Sym}_{A}(A[S]) 
  \arrow[u,"d"]
  \arrow[r,equal]
  & 
  \on{Sym}_{A}(A[S])
  \arrow[u,"d"]
\end{tikzcd}
\]
where we use lemma \ref{4702374903} for the description of $\mathbb{D}_+$-module structure of $\on{L} \Omega_{B/A}^{+, \bullet} $ on degree 0 and 1.
The universal property of cotangent complex forces the horizontal map to be an isomorphism. 

 Moreover, we can compute the underlying object 
 \begin{align*}
     \on{colim}      \on{HH}_{\on{Fil}}(\on{Sym}^{\bullet}_A(A[S])/A) 
    &
    \cong 
     \on{colim}
     \on{Sym}^{\bullet}_{\on{ins}^0A}(\on{ins}^0A[S]
    \otimes_{ \on{ins}^0A} \on{ins}^0A\lbrack \underline{\mathbb{T}_{\on{Fil}}}\rbrack )
    \\
    & 
    \cong 
     \on{Sym}^{\bullet}_{A}(A[S]
    \otimes_{ A} A\lbrack \underline{\mathbb{S}^1}\rbrack )
    \\
    & 
    \cong 
     \on{Sym}^{\bullet}_{A}(A[S] \oplus \Sigma A[S])
     \\
     & 
     \cong 
      \on{Sym}^{\bullet}_{A}(A[S] \bigsqcup_{ A[S] \oplus  A[S]} A[S])
      \\
      & 
      \cong 
       \on{Sym}^{\bullet}_{A}(A[S]) \otimes_{ \on{Sym}^{\bullet}_{A}(A[S]) \otimes_A  \on{Sym}^{\bullet}_{A}(A[S]) } 
        \on{Sym}^{\bullet}_{A}(A[S]) 
        \\
        &
        \cong 
        \on{HH}( \on{Sym}^{\bullet}_{A}(A[S]) /A)
 \end{align*}
which we see is actually $\mathbb{S}^1$-equivariant, namely, the $\mathbb{Z}[\mathbb{S}^1] = \mathbb{Z}[\epsilon]/\epsilon^2$-module structure is specified by $\epsilon$ sending $A[S]$ to $\Sigma A[S]$ in the third equivalence. Note that we use the canonical equivalence
\[
(\on{DAlg} \on{Cond} D(\mathbb{Z})_{A/})^{B\mathbb{S}^1} 
\cong 
\on{cMod}_{A^{\mathbb{S}^1}}
(\on{DAlg}\on{Cond}D_{} (\mathbb{Z})_{A/})
\]
see \cite[remark 6.1.3.]{Raksit2026}, where the proof goes through verbally here. 

Passing to colimits of $\on{Sym}^{\bullet}_{A}(A[S])$
we obtain the results for general condensed animated ring $B$. 
\end{proof}

\begin{remark}
    We can also define the filtered Hochschild homology for condensed rings directly by sheafifying the sectionwise filtered Hochschild homology as defined in \cite[definition 6.2.5.]{Raksit2026}, and deduce the universal property by sheafifying the universal property established in $loc.cit.$. This is the approach we take in section \ref{3958203834}, and we choose to phrase slightly differently here. 
\end{remark}

\begin{lemma}
 Let $A \in \on{DAlg} \on{Cond} D(\mathbb{Z})$, and we assume that $A$ is connective, in other words, $A$ is a condensed animated ring, then for any condensed animated ring $B$ over $A$, then the filtration on $ \on{HH}_{\on{Fil}}(B/A)$ is complete.     
\end{lemma}

\begin{proof}
This is exactly the same as the proof of \cite[proposition 6.2.9.]{Raksit2026}. We observe that since both $B$ and $A$ are connective, so does $\mathbb{L}_{B/A}$, and the graded pieces $\on{gr}^i\on{HH}_{\on{Fil}}(B/A)= \Sigma^i \wedge^i \mathbb{L}_{B/A}$ 
is $i$-connective. It implies that the $i$-th filtration of $\on{HH}_{\on{Fil}}(B/A)$
is $i$-connective as well since this is true for $B= \on{Sym}^{\bullet}_{A}(A[S]) $
because the filtration is split as we see in 
the proof of proposition \ref{3074375372}, and then taking colimits of $\on{Sym}^{\bullet}_{A}(A[S])$ we have that in general. Now the $t$-structure on (light) condensed $D(\mathbb{Z})$ is left complete by \cite[theorem 4.11]{2025arXiv250618174C}, so the limit of the filtration of $\on{HH}_{\on{Fil}}(B/A)$ is $0$, proving that it is complete. 
\end{proof}

We now proceed to define the HKR filtered Hochschild homology for  analytic ring. To fix notation,  Let us first define the obvious notion of Hochschild homology for analytic animated rings.  

\begin{definition}
     Let $\mathcal{A} \rightarrow \mathcal{B}$ be a map of analytic animated rings, we define 
     \[
     \mathbb{HH}(\mathcal{B}/\mathcal{A}) := \mathcal{B} \otimes_{\mathcal{B}\otimes_{\mathcal{A}} \mathcal{B}} \mathcal{B} \in (\on{AnRing}_{\mathcal{A}/})^{\on{B} \mathbb{S}^1}
     \]
   where the tensor product is the colimit taken in the category of analytic animated rings.   
   We define 
   \[
    \on{HH}(\mathcal{B}/\mathcal{A}) :=  \mathbb{HH}(\mathcal{B}/\mathcal{A})^{\triangleright} \in \on{CAlg}D(\mathcal{A})^{B\mathbb{S}^1}
   \]
   as the underlying condensed ring in, which also lies in $D(\mathcal{B})$.
\end{definition}

\begin{remark}
    It is not hard to see from the definition of colimits in analtic rings that 
    \[
    \on{HH}(\mathcal{B}/\mathcal{A})
    \cong 
    \on{HH}(\mathcal{B}_{\mathcal{A}/}/\mathcal{A}) \otimes_{\mathcal{B}_{\mathcal{A}/}} \mathcal{B}
    \]
    where $\mathcal{B}_{\mathcal{A}/}$ 
    is the induced analytic structure from $\mathcal{A}$ on $\mathcal{B}^{\triangleright}$. 
\end{remark}

Let us now define the uncompleted HKR filtered Hochschild homology for  analytic ring.

\begin{definition}
    Let $\mathcal{A} \rightarrow \mathcal{B}$ be a map of analytic animated rings,  we define 
    \[
    \mathbb{HH}^{\on{un}}_{\on{fil}}(\mathcal{B}/\mathcal{A})=(\mathbb{HH}^{\on{un}}_{\on{fil}}(\mathcal{B}/\mathcal{A})^{\triangleright}, D(\mathbb{HH}^{\on{un}}_{\on{fil}}(\mathcal{B}/\mathcal{A}))) \in  \on{FilAnDAG}^{\on{un}}
    \]
    with the underlying  filtered condensed derived ring
    \[
    \mathbb{HH}^{\on{un}}_{\on{fil}}(\mathcal{B}/\mathcal{A})^{\triangleright}:=
    \on{HH}_{\on{Fil}}(\mathcal{B}^{\triangleright}/\mathcal{A}^{\triangleright})  \in  \on{DAlg}\on{Cond}\on{Fil}_{}D_{} (\mathbb{Z})_{\on{ins}^0 \mathcal{A}^{\triangleright}/}
    \]
    and the analytic structure
    $D(\mathbb{HH}^{\on{un}}_{\on{fil}}(\mathcal{B}/\mathcal{A}))$
    coming from 
    $\on{gr}^0 \on{HH}_{\on{Fil}}(\mathcal{B}^{\triangleright}/\mathcal{A}^{\triangleright})  \cong \mathcal{B}^{\triangleright} $
    and the analytic structure on $\mathcal{B}$
    using proposition \ref{xnvmnvd}. 
\end{definition}

We first observe that the analytic structure on $\mathbb{HH}^{\on{un}}_{\on{fil}}(\mathcal{B}/\mathcal{A})$ 
is induced from $\on{ins}^0 \mathcal{B} $. 

\begin{lemma} \label{669780789798}
     Let $\mathcal{A} \rightarrow \mathcal{B}$ be a map of analytic animated rings, then the analytic structure on $\mathbb{HH}^{\on{un}}_{\on{fil}}(\mathcal{B}/\mathcal{A})$
   is induced from (in the sense of corollary \ref{74329874389}) the analytic structure on 
   $\on{ins}^0 \mathcal{B}^{\triangleright}$
   through the map 
   \[
   \on{ins}^0 \mathcal{B}^{\triangleright} 
   \longrightarrow 
    \on{HH}_{\on{Fil}}(\mathcal{B}^{\triangleright}/\mathcal{A}^{\triangleright})
   \]
  as in lemma \ref{70237450834}, where the analytic structure on $\on{ins}^0 \mathcal{B}^{\triangleright}$ 
  comes from the analytic structure on $\mathcal{B}$ and  $\on{gr}^0\on{ins}^0 \mathcal{B}^{\triangleright} \cong \mathcal{B}^{\triangleright}$
 using proposition \ref{xnvmnvd}.
\end{lemma}

\begin{proof}
    By proposition \ref{xnvmnvd}, we have 
    \[
    D(\on{ins}^0 \mathcal{B}^{\triangleright} ) = \on{Mod}_{\mathbb{D}_-} \on{Gr}D(\mathcal{B})
    \]
    so 
    \[
    \on{Mod}_{C\tau} 
    \on{Mod}_{ \on{HH}_{\on{Fil}}(\mathcal{B}^{\triangleright}/\mathcal{A}^{\triangleright})}
    D(\on{ins}^0 \mathcal{B}^{\triangleright} )
    = 
    \on{Mod}_{\on{gr}^{\bullet}  \on{HH}_{\on{Fil}}(\mathcal{B}^{\triangleright}/\mathcal{A}^{\triangleright}) }\on{Gr}D(\mathcal{B})
    \]
    which is exactly 
    \[
    \on{Mod}_{C\tau} 
    D( \mathbb{HH}^{\on{un}}_{\on{fil}}(\mathcal{B}/\mathcal{A}))
    \]
    again by proposition \ref{xnvmnvd}, and the claim follows from proposition \ref{2378937583}. 
\end{proof}

\begin{proposition} \label{7234703570329}
    The $\underline{\mathbb{T}_{\on{Fil}}^{\vee}}$-comodule structure on $\on{HH}_{\on{Fil}}(\mathcal{B}^{\triangleright}/\mathcal{A}^{\triangleright}) $
    induces 
    \[
    \mathbb{HH}^{\on{un}}_{\on{fil}}(\mathcal{B}/\mathcal{A}) \in  \on{cMod}_{\underline{\mathbb{T}_{\on{Fil}}^{\vee}} } \on{FilAnDAG}^{\on{un}}
    \]
    where we give $\underline{\mathbb{T}_{\on{Fil}}^{\vee}} $
    the trivial analytic structure as in example \ref{xfower870}, which also provides the coalgebra structure of $\underline{\mathbb{T}_{\on{Fil}}^{\vee}} $ in $\on{FilAnDAG}^{\on{un}}$ with cocartesian symmetric monoidal structure. 
\end{proposition}

\begin{proof}
We need to check all the coherent comodule structure maps of 
$ \on{HH}_{\on{Fil}}(\mathcal{B}^{\triangleright}/\mathcal{A}^{\triangleright}) $
over $\underline{\mathbb{T}_{\on{Fil}}^{\vee}} $
lies in $ \on{FilAnDAG}^{\on{un}}$, but mapping spaces of $\on{FilAnDAG}^{\on{un}}$ by its very definition
are connected components of corresponding mapping spaces in $\on{DAlg}\on{Cond}\on{Fil}_{}D_{} (\mathbb{Z})$, in other words, it is enough to prove the claim by passing to $\pi_0$ of the categories. In this case, comodule structures are specified by the single comodule structure map 
    \begin{equation} \label{nvcjgefiue}
     \on{HH}_{\on{Fil}}(\mathcal{B}^{\triangleright}/\mathcal{A}^{\triangleright}) 
     \longrightarrow
\underline{\mathbb{T}_{\on{Fil}}^{\vee}}  \bigsqcup_{}  \on{HH}_{\on{Fil}}(\mathcal{B}^{\triangleright}/\mathcal{A}^{\triangleright}), 
    \end{equation}
    and we need to prove it
    is actually a map in 
$\on{FilAnDAG}^{\on{un}}$, in other words, we need to prove that the forgetful functor maps the complete modules on the right hand side to complete modules of the left hand side. 

We observe that the analytic ring structure on the right hand side is also 
given by proposition \ref{xnvmnvd}. Since $
\underline{\mathbb{T}_{\on{Fil}}^{\vee}}  $
is given the trivial analytic structure, by the description of colimits in $\on{FilAnDAG}^{\on{un}}$ as in lemma \ref{weuryeiwyrue}, we see that 
\[
D(
\underline{\mathbb{T}_{\on{Fil}}^{\vee}} \bigsqcup  \on{HH}_{\on{Fil}}(\mathcal{B}^{\triangleright}/\mathcal{A}^{\triangleright})) = \on{Mod}_{\mathbb{T}_{\on{Fil}}^{\vee}} D(\on{HH}_{\on{Fil}}(\mathcal{B}^{\triangleright}/\mathcal{A}^{\triangleright}))
\]
where we view 
$D(\on{HH}_{\on{Fil}}(\mathcal{B}^{\triangleright}/\mathcal{A}^{\triangleright}))$
as tensored over $\on{Fil} \on{Cond}D(\mathbb{Z})$
through the obvious symmetric monoidal functor. From the definition of $D(\on{HH}_{\on{Fil}}(\mathcal{B}^{\triangleright}/\mathcal{A}^{\triangleright}))$, we have that 
\begin{align*}
\on{Mod}_{C\tau}
\on{Mod}_{\mathbb{T}_{\on{Fil}}^{\vee}} D(\on{HH}_{\on{Fil}}(\mathcal{B}^{\triangleright}/\mathcal{A}^{\triangleright}))
& = \on{Mod}_{\on{gr}^{\bullet}\mathbb{T}_{\on{Fil}}^{\vee}} 
\on{Mod}_{\on{gr}^{\bullet}\on{HH}_{\on{Fil}}(\mathcal{B}^{\triangleright}/\mathcal{A}^{\triangleright})}
\on{Gr} D(\mathcal{B})
\\
& = \on{Mod}_{\on{gr}^{\bullet}(\mathbb{T}_{\on{Fil}}^{\vee} \otimes \on{HH}_{\on{Fil}}(\mathcal{B}^{\triangleright}/\mathcal{A}^{\triangleright}))}
\on{Gr} D(\mathcal{B})
\\
& = \on{Mod}_{\on{gr}^{\bullet}(\mathbb{T}_{\on{Fil}}^{\vee} \otimes \on{HH}_{\on{Fil}}(\mathcal{B}^{\triangleright}/\mathcal{A}^{\triangleright}))}
\on{Mod}_{\on{ins}^0\on{gr}^{0}(\mathbb{T}_{\on{Fil}}^{\vee} \otimes \on{HH}_{\on{Fil}}(\mathcal{B}^{\triangleright}/\mathcal{A}^{\triangleright}))}
\on{Gr} D(\mathcal{B})
\\
& = \on{Mod}_{\on{gr}^{\bullet}(\mathbb{T}_{\on{Fil}}^{\vee} \otimes \on{HH}_{\on{Fil}}(\mathcal{B}^{\triangleright}/\mathcal{A}^{\triangleright}))}
\on{Gr}
\on{Mod}_{\on{gr}^{0}(\mathbb{T}_{\on{Fil}}^{\vee} \otimes \on{HH}_{\on{Fil}}(\mathcal{B}^{\triangleright}/\mathcal{A}^{\triangleright}))}D(\mathcal{B})
\end{align*}
which means that the analytic structure on
$\underline{\mathbb{T}_{\on{Fil}}^{\vee}} \bigsqcup  \on{HH}_{\on{Fil}}(\mathcal{B}^{\triangleright}/\mathcal{A}^{\triangleright})$
is given by the analytic structure on 
$\on{gr}^{0}(\mathbb{T}_{\on{Fil}}^{\vee} \otimes \on{HH}_{\on{Fil}}(\mathcal{B}^{\triangleright}/\mathcal{A}^{\triangleright}))$
through proposition \ref{xnvmnvd}, where we give the analytic structure on 
$\on{gr}^{0}(\mathbb{T}_{\on{Fil}}^{\vee} \otimes \on{HH}_{\on{Fil}}(\mathcal{B}^{\triangleright}/\mathcal{A}^{\triangleright}))$
induced from $\mathcal{B}$
through the canonical map 
\[
\mathcal{B}^{\triangleright} = \on{gr}^{0}\on{HH}_{\on{Fil}}(\mathcal{B}^{\triangleright}/\mathcal{A}^{\triangleright})
\longrightarrow
\on{gr}^{0}(\mathbb{T}_{\on{Fil}}^{\vee} \otimes \on{HH}_{\on{Fil}}(\mathcal{B}^{\triangleright}/\mathcal{A}^{\triangleright})), 
\]
which is the identity map on the second factor. 

By lemma \ref{79832403983}, it is enough to prove that the comodule structure map induces on $\on{gr}^0$ a map of analytic animated rings. Using $\on{gr}^{\bullet} \mathbb{T}_{\on{Fil}}^{\vee} = \mathbb{Z} \oplus  \Sigma^{-1}\mathbb{Z}[-1]$, we compute that 
\[
\on{gr}^{0}(\mathbb{T}_{\on{Fil}}^{\vee} \otimes \on{HH}_{\on{Fil}}(\mathcal{B}^{\triangleright}/\mathcal{A}^{\triangleright}))
= \mathcal{B}^{\triangleright} \oplus \mathbb{L}_{\mathcal{B}^{\triangleright}/\mathcal{A}^{\triangleright}}
\]
with the square zero algebra structure on $\mathbb{L}_{\mathcal{B}^{\triangleright}/\mathcal{A}^{\triangleright}}$, 
and the comodule structure map (\ref{nvcjgefiue}) on $\on{gr}^0$ is given by the canonical map 
\[
(\on{id},d): \mathcal{B}^{\triangleright} \longrightarrow \mathcal{B}^{\triangleright} \oplus \mathbb{L}_{\mathcal{B}^{\triangleright}/\mathcal{A}^{\triangleright}}
\]
where $d$ is universal map classifying universal square zero extensions, 
this is easy to see for $\mathcal{B}^{\triangleright}= \on{Sym}^{\bullet}(\mathcal{A}^{\triangleright}[S])$ 
by going through the proof of proposition \ref{3074375372}, and the general case follows by taking colimits of it. The analytic structure is nilpotence invariant, so depends only on the quotient $\mathcal{B}^{\triangleright}$ of the square zero extension 
$\mathcal{B}^{\triangleright} \oplus \mathbb{L}_{\mathcal{B}^{\triangleright}/\mathcal{A}^{\triangleright}}$, and the compose of 
$(\on{id},d)$ with it is identity, so it induces a map of analytic animated rings, proving that the comodule structure map is a map in $\on{FilAnDAG}^{\on{un}}$. 
\end{proof}

We can now define the filtered Hochschild homology of analytic rings. 

\begin{definition}
    Let $\mathcal{A} \rightarrow \mathcal{B}$ be a map of analytic animated rings,  we define 
    \[
     \mathbb{HH}_{\on{Fil}}(\mathcal{B}/\mathcal{A})=( \mathbb{HH}_{\on{Fil}}(\mathcal{B}/\mathcal{A})^{\triangleright}, D( \mathbb{HH}_{\on{Fil}}(\mathcal{B}/\mathcal{A}))) \in  \on{FilAnDAG}
    \]
    to be the image of $\mathbb{HH}^{\on{un}}_{\on{Fil}}(\mathcal{B}/\mathcal{A})$
    under the left adjoint functor 
    $\on{FilAnDAG}^{\on{un}} \longrightarrow \on{FilAnDAG}$
    in
    proposition \ref{wyrw3874}. 
\end{definition}

\begin{theorem}[Analytic Hochschild-Kostant-Rosenberg] \label{7987034jlkje}
     Let $\mathcal{A} \rightarrow \mathcal{B}$ be a map of analytic animated rings, then 
     \[
    \mathbb{HH}_{\on{Fil}}(\mathcal{B}/\mathcal{A}) \in  \on{cMod}_{\underline{\mathbb{T}_{\on{Fil}}^{\vee}} } \on{FilAnDAG}_{\on{ins}^0 \mathcal{A}/}
    \]
    where we give the cocartesian symmetric monoidal structure on $\on{FilAnDAG}$. Moreover, we have
    \[
    \on{gr}^{\bullet}
     \mathbb{HH}_{\on{Fil}}(\mathcal{B}/\mathcal{A})
     \cong 
     \widehat{\mathbb{DR}}_{\mathcal{B}/\mathcal{A}}
     \in \on{cMod}_{\underline{\mathbb{D}_{+}^{\vee} }}  \on{GrAnDAG}_{\on{ins}^0 \mathcal{A} /} 
    \]
    and
    \[
    \on{colim}  \mathbb{HH}_{\on{Fil}}(\mathcal{B}/\mathcal{A})^{\triangleright} \cong 
    \mathbb{HH}(\mathcal{B}/\mathcal{A})^{\triangleright} = \on{HH}(\mathcal{B}/\mathcal{A}). 
    \]
    In particular, we have
    \[
     \on{gr}^{\bullet}
     \mathbb{HH}_{\on{Fil}}(\mathcal{B}/\mathcal{A})^{\triangleright}
     \cong 
     \widehat{\on{dR}}_{\mathcal{B}/\mathcal{A}}
     \in \on{cMod}_{\underline{\mathbb{D}_{+}^{\vee} }} \on{Gr} D(\mathcal{A})
    \]
\end{theorem}

\begin{proof}
    The left adjoint $\on{FilAnDAG}^{\on{un}} \longrightarrow \on{FilAnDAG}$  
 in proposition \ref{wyrw3874}
    preserves comodules since it is symmetric monoidal as both sides are given the cocartesian symmetric monoidal structures and left adjoints preserves colimits, so proposition \ref{7234703570329} gives rise to $ \mathbb{HH}_{\on{Fil}}(\mathcal{B}/\mathcal{A}) \in  \on{cMod}_{\underline{\mathbb{T}_{\on{Fil}}^{\vee}} } \on{FilAnDAG}_{\on{ins}^0 \mathcal{A}/}$. 
    
    We know from proposition \ref{3074375372} that 
    \[
    \on{gr}^{\bullet}
    \mathbb{HH}_{\on{Fil}}(\mathcal{B}/\mathcal{A})^{\on{un}} \cong 
    \widehat{\mathbb{DR}}_{\mathcal{B}/\mathcal{A}}^{\on{un}}
    \in  \on{cMod}_{\underline{\mathbb{D}_{+}^{\vee}} } \on{GrAnDAG}^{\on{un}}_{\on{ins}^0 \mathcal{A}/}.
    \]
Indeed, both sides have the same condensed graded ring, and the analytic structure are both induced from the analytic structure on the degree 0 part, which is $\mathcal{B}$. Taking completions we have 
 \[
    \on{gr}^{\bullet}
     \mathbb{HH}_{\on{Fil}}(\mathcal{B}/\mathcal{A})
     \cong 
     \widehat{\mathbb{DR}}_{\mathcal{B}/\mathcal{A}}
     \in \on{cMod}_{\underline{\mathbb{D}_{+}^{\vee}} }  \on{GrAnDAG}_{\on{ins}^0 \mathcal{A}/}. 
    \]

   It follows that the filtration on 
   $\mathbb{HH}_{\on{Fil}}(\mathcal{B}/\mathcal{A})^{\triangleright}$ 
   is non-negatively filtered, so $\on{colim} \mathbb{HH}_{\on{Fil}}(\mathcal{B}/\mathcal{A})^{\triangleright}$
   is the $0$-th filtration on 
   $\mathbb{HH}_{\on{Fil}}(\mathcal{B}/\mathcal{A})^{\triangleright}$. It follows from lemma \ref{669780789798} that 
   \[
     \mathbb{HH}_{\on{Fil}}(\mathcal{B}/\mathcal{A})^{\triangleright} = 
     L_{ \mathbb{HH}^{\on{un}}_{\on{Fil}}(\mathcal{B}/\mathcal{A})} ( \on{HH}_{\on{Fil}}(\mathcal{B}^{\triangleright}/\mathcal{A}^{\triangleright})   )
     = L_{\on{ins}^0 \mathcal{B}} ( \on{HH}_{\on{Fil}}(\mathcal{B}^{\triangleright}/\mathcal{A}^{\triangleright})   )
     \]
     since $\on{ins}^0 \mathcal{B}^{\triangleright}$
     is concentrated at degree $0$,
     \[
     D(\on{ins}^0 \mathcal{B}^{\triangleright}) = \on{Fil} D(\mathcal{B})
     \subset
      \on{Fil} \on{Mod}_{\mathcal{B}^{\triangleright}} \on{Cond} D(\mathbb{Z}) =
      \on{Mod}_{\on{ins}^0 \mathcal{B}} \on{Cond} \on{Fil} D(\mathbb{Z}),
     \]
     and the localization $L_{\on{ins}^0 \mathcal{B}}$ is pointwise the localization $L_{\mathcal{B}}$, so the 0-th filtration of $ \mathbb{HH}_{\on{Fil}}(\mathcal{B}/\mathcal{A})^{\triangleright}$ 
     is 
     \[
     L_{\on{ins}^0 \mathcal{B}} ( \on{HH}_{\on{Fil}}(\mathcal{B}^{\triangleright}/\mathcal{A}^{\triangleright}))^0  
     = 
     L_{\mathcal{B}} ( \on{HH}(\mathcal{B}^{\triangleright}/\mathcal{A}^{\triangleright}) )
     = \on{HH}(\mathcal{B}/\mathcal{A})
     \]
     where we use proposition \ref{3074375372}
     that 
     $\on{colim}\on{HH}_{\on{Fil}}(\mathcal{B}^{\triangleright}/\mathcal{A}^{\triangleright}) = \on{HH}(\mathcal{B}^{\triangleright}/\mathcal{A}^{\triangleright})$
     is the 0-th filtration of $\on{HH}_{\on{Fil}}(\mathcal{B}^{\triangleright}/\mathcal{A}^{\triangleright})$, and the analytic structure is induced from $\mathbb{HH}(\mathcal{B}/\mathcal{A})$ is induced from $\mathcal{B}$. 
\end{proof}

\begin{theorem} \label{weirdzdfd}
    Let $\mathcal{A} \rightarrow \mathcal{B}$ be a map of analytic animated rings, then there exists a filtration on the $\mathbb{S}^1$-fixed point $ \on{HH}(\mathcal{B}/\mathcal{A})^{h\mathbb{S}^1} $
    with graded pieces 
    \[
    \on{gr}^{i}
    (\on{HH}(\mathcal{B}/\mathcal{A})^{h\mathbb{S}^1} )
    = \Sigma^{2i}\widehat{\on{dR}}_{\mathcal{B}/\mathcal{A}}^{\geq i}, 
    \]
    where 
    $
    \widehat{\on{dR}}_{\mathcal{B}/\mathcal{A}}^{\geq i}
    $
    is the i-th filtration of the filtered  $\widehat{\on{dR}}_{\mathcal{B}/\mathcal{A}}^{\geq \bullet} \in \on{Fil}^{\wedge} D(\mathcal{A})$
  defined   in 
 definition \ref{sfdsfjekl}. 
    
    Similarly, there exists a filtration on the Tate construction 
    $ \on{HH}(\mathcal{B}/\mathcal{A})^{t\mathbb{S}^1} $
    with graded pieces
    \[
     \on{gr}^{i}
    (\on{HH}(\mathcal{B}/\mathcal{A})^{t\mathbb{S}^1} )
    = \Sigma^{2i}\widehat{\on{dR}}_{\mathcal{B}/\mathcal{A}}, 
    \]
    where 
    \[
    \widehat{\on{dR}}_{\mathcal{B}/\mathcal{A}} := \on{colim} \widehat{\on{dR}}_{\mathcal{B}/\mathcal{A}}^{\geq \bullet} \in D(\mathcal{A}). 
    \]
 Moreover,   the canonical map 
    \[
    \on{HH}(\mathcal{B}/\mathcal{A})^{h\mathbb{S}^1} 
    \longrightarrow 
    \on{HH}(\mathcal{B}/\mathcal{A})^{t\mathbb{S}^1} 
    \]
    induces the filtration map (shifted by $\Sigma^{2i}$)
    \[
    \widehat{\on{dR}}_{\mathcal{B}/\mathcal{A}}^{\geq i}
    \longrightarrow
\widehat{\on{dR}}_{\mathcal{B}/\mathcal{A}}
    \]
    on graded pieces, and the universal map of fixed points
    \[
    \on{HH}(\mathcal{B}/\mathcal{A})^{h\mathbb{S}^1} 
    \longrightarrow 
    \on{HH}(\mathcal{B}/\mathcal{A})
    \]
    induces the quotient map
    \[
     \Sigma^{2i} \widehat{\on{dR}}_{\mathcal{B}/\mathcal{A}}^{\geq i}
    \longrightarrow
     \Sigma^{2i}
\on{gr}^i\widehat{\on{dR}}_{\mathcal{B}/\mathcal{A}}^{\geq *} =\on{Sym}^i( \Sigma \mathbb{L}_{\mathcal{B}/\mathcal{A}})
    \]
    on graded pieces.
\end{theorem}

\begin{proof}
   We apply \cite[proposition 6.3.3.]{Raksit2026}  to 
    $ \mathbb{HH}_{\on{Fil}}(\mathcal{B}/\mathcal{A})^{\triangleright} \in \on{cMod}_{\mathbb{T}^{\vee}_{\on{Fil}}}  \on{Fil} D(\mathcal{A}) $, where we forget the derived algebra structure, and we get filtered objects $ \mathbb{HH}_{\on{Fil}}(\mathcal{B}/\mathcal{A})^{\triangleright, h \mathbb{T}_{\on{Fil}}} $ and
    $\mathbb{HH}_{\on{Fil}}(\mathcal{B}/\mathcal{A})^{\triangleright, t \mathbb{T}_{\on{Fil}}} $
    with graded pieces as displayed on the right hand side by theorem \ref{7987034jlkje}.  Then \cite[proposition 6.3.5.]{Raksit2026} tells us that the filtered object is a filtration on $ (\on{colim}\mathbb{HH}_{\on{Fil}}(\mathcal{B}/\mathcal{A})^{\triangleright})^{ h \mathbb{S}^1} = \on{HH}(\mathcal{B}/\mathcal{A})^{h\mathbb{S}^1}$
    and 
    $(\on{colim}\mathbb{HH}_{\on{Fil}}(\mathcal{B}/\mathcal{A})^{\triangleright})^{ t \mathbb{S}^1} = \on{HH}(\mathcal{B}/\mathcal{A})^{t\mathbb{S}^1} $ respectively, since we have seen  
    $\mathbb{HH}_{\on{Fil}}(\mathcal{B}/\mathcal{A})^{\triangleright})$
    is non-negatively filtered. 
\end{proof}

 \begin{remark} \label{xueoiwuroew}
     It is not hard to see from proof that the equivalence $\on{gr}^{\bullet}
     \mathbb{HH}_{\on{Fil}}(\mathcal{B}/\mathcal{A})^{\triangleright}
     \cong 
     L\Omega_{\mathcal{B}/\mathcal{A}}^{+,\bullet}$
     is an equivalence of commutative algebras in $\on{Gr} D(\mathcal{B})$.   Moreover, the identifications in theorem \ref{weirdzdfd} is also multiplicative by \cite[remark 6.3.4.]{Raksit2026}.  More precisely, we have an equivalence of graded commutative algebra
     \[
     \on{gr}^{\bullet}
    (\on{HH}(\mathcal{B}/\mathcal{A})^{h\mathbb{S}^1} )
    = \Sigma^{2*} \on{Und}(\widehat{\on{dR}}_{\mathcal{B}/\mathcal{A}}^{\geq \bullet})
     \]
     where $\on{Und}: \on{Fil} D(\mathcal{A}) \longrightarrow \on{Gr} D(\mathcal{A})$ is the underlying graded objects of filtered category, i.e. we forget the transition maps of the filtration and retains only the indexed filtration, and $\Sigma^{2*}: \on{Gr} D(\mathcal{A}) \longrightarrow \on{Gr} D(\mathcal{A})$
     sends a graded object $X^{\bullet}$ to 
     $\Sigma^{2\bullet}X^{\bullet}$. 
     They are lax symmetric monoidal functors, so the commutative algebra structure on the filtered object $\widehat{\on{dR}}_{\mathcal{B}/\mathcal{A}}^{\geq \bullet}$ induces a commutative algebra structure on 
     $\Sigma^{2*} \on{Und}(\widehat{\on{dR}}_{\mathcal{B}/\mathcal{A}}^{\geq \bullet})$. 

     Similarly, we have an equivalence of graded commutative algebra objects
      \[
     \on{gr}^{\bullet}
    (\on{HH}(\mathcal{B}/\mathcal{A})^{t\mathbb{S}^1} )
    = \Sigma^{2*} \delta( \widehat{\on{dR}}_{\mathcal{B}/\mathcal{A}}) = \widehat{\on{dR}}_{\mathcal{B}/\mathcal{A}} [\sigma, \sigma^{-1}], 
    \]
    where $\delta: D(\mathcal{A}) \longrightarrow \on{Gr} D(\mathcal{A})$
    is the diagonal functor, which sends an object to the constant graded object. It is lax symmetric monoidal, so the commutative algebra structure on $ \widehat{\on{dR}}_{\mathcal{B}/\mathcal{A}}$ 
    induces one on 
    $\delta( \widehat{\on{dR}}_{\mathcal{B}/\mathcal{A}})$, and so does $\Sigma^{2*} \delta( \widehat{\on{dR}}_{\mathcal{B}/\mathcal{A}})$. The graded algebra structure can be more explicitly expressed as $\widehat{\on{dR}}_{\mathcal{B}/\mathcal{A}} [\sigma, \sigma^{-1}]$,  where $\sigma$ is a free polynomial variable of homological degree 2 and graded degree 1. 
 \end{remark}

\section{Bhatt-Morrow-Scholze \rom{2} nonsense}

Before constructing the concrete cohomology theories, we first give an axiomatic treatment of cohomology theory out of derived de Rham complex and Hochschild homology. This is formally very similar to \cite{BMS2}. 

Through this section, we fix a $\mathcal{A} \longrightarrow  \mathcal{B}$ be a map of static analytic animated rings over $\mathbb{Q}$ such that $\mathbb{L}_{\mathcal{B}/\mathcal{A}} \cong \Sigma \mathcal{B}^{\triangleright}$ as an element of $D(\mathcal{B})$. 

\begin{lemma} \label{2374032493}
Let $\mathcal{A} \longrightarrow  \mathcal{B}$ be a map of static analytic animated rings over $\mathbb{Q}$ such that $\mathbb{L}_{\mathcal{B}/\mathcal{A}} \cong \Sigma \mathcal{B}^{\triangleright}$ as an element of $D(\mathcal{B})$,
then 
$\widehat{\on{dR}}_{\mathcal{B}/\mathcal{A}} \in \on{CAlg} D(\mathcal{B})^{\heartsuit}$ is static and so does the filtration 
$\widehat{\on{dR}}^{\geq *}_{\mathcal{B}/\mathcal{A}}$. Moreover, $\widehat{\on{dR}}^{\geq 1}_{\mathcal{B}/\mathcal{A}}$  defines a principal ideal $I$ of $\widehat{\on{dR}}_{\mathcal{B}/\mathcal{A}} $
with quotient $\mathcal{B}^{\triangleright}$, and $\widehat{\on{dR}}^{\geq n}_{\mathcal{B}/\mathcal{A}} = I^n$.

For any analytic animated ring $\mathcal{R}$ over $\mathcal{B}$, 
$
\widehat{\on{dR}}_{\mathcal{R}/\mathcal{A}}
$
is a $\widehat{\on{dR}}_{\mathcal{B}/\mathcal{A}}$-algebra such that 
\[
\widehat{\on{dR}}_{\mathcal{R}/\mathcal{A}} \otimes_{\widehat{\on{dR}}_{\mathcal{B}/\mathcal{A}}} \widehat{\on{dR}}_{\mathcal{B}/\mathcal{A}}/I 
\cong 
\widehat{\on{dR}}_{\mathcal{R}/\mathcal{A}} \otimes_{\widehat{\on{dR}}_{\mathcal{B}/\mathcal{A}}} \mathcal{B}^{\triangleright} \cong \widehat{\on{dR}}_{\mathcal{R}/\mathcal{B}}.
\]
Moreover, the Katz-Oda filtration 
\[
F_{\on{KO}}^i := \widehat{\on{dR}}_{\mathcal{R}/\mathcal{A}} \otimes_{\widehat{\on{dR}}_{\mathcal{B}/\mathcal{A}}} 
\widehat{\on{dR}}_{\mathcal{B}/\mathcal{A}}^{\geq i}
\]
on $\widehat{\on{dR}}_{\mathcal{R}/\mathcal{A}}$ has graded pieces
\[
\on{gr}_{\on{KO}}^i \widehat{\on{dR}}_{\mathcal{R}/\mathcal{A}}
\cong 
\widehat{\on{dR}}_{\mathcal{R}/\mathcal{B}} \otimes_{\mathcal{B}} 
\Sigma^{-2i}
 \on{Sym}^i_{\mathcal{B}}( \Sigma
\mathbb{L}_{\mathcal{B}/\mathcal{A}}) 
\cong 
\widehat{\on{dR}}_{\mathcal{R}/\mathcal{B}}.  
\] 
\end{lemma}

\begin{proof}
  It follows from   (\ref{zdfhweuoofd}) that the i-th graded pieces of $\widehat{\on{dR}}^{\geq *}_{\mathcal{B}/\mathcal{A}}$ 
  is 
  \[
  \Sigma^{-2i}
  \on{Sym}^i_{\mathcal{B}}( \Sigma
\mathbb{L}_{\mathcal{B}/\mathcal{A}})
=
\Sigma^{-2i}
  \on{Sym}^i_{\mathcal{B}}( \Sigma^2
\mathcal{B}^{\triangleright})
=
\on{Sym}^i_{\mathcal{B}}(
\mathcal{B}^{\triangleright})
\]
where we use the standard Illusie's formula $\on{Sym}^i (\Sigma (-)) = \Sigma^i \on{Sym}^i(-)$ over characteristic 0, namely the free divided power algebra coincides with symmetric power algebra in characteristic 0, which holds in the condensed setting by sheafifying the classical algebraic formula and passing to completion. Now our assumption that $\mathcal{B}$ is static tells us that this graded piece is static, which proves $\widehat{\on{dR}}^{\geq *}_{\mathcal{B}/\mathcal{A}}$ and 
$\widehat{\on{dR}}_{\mathcal{B}/\mathcal{A}}$ is static
 since the filtration is  complete. 
 
 It follows that the filtration  $\widehat{\on{dR}}^{\geq 1}_{\mathcal{B}/\mathcal{A}}$ embeds into $\widehat{\on{dR}}_{\mathcal{B}/\mathcal{A}}$ since its cofiber is static and we have a cofiber sequence of static objecs. It  defines an ideal $I$ since the $\widehat{\on{dR}}^{\geq *}_{\mathcal{B}/\mathcal{A}}$ is a commutative algebra in the filtered category as observed in remark \ref{2307krieo}. It follows from the computation of graded pieces that any lifting of a generator of $
 \on{gr}^1  \widehat{\on{dR}}^{\geq *}_{\mathcal{B}/\mathcal{A}}  \cong \mathcal{B}^{\triangleright}$
 to $\widehat{\on{dR}}^{\geq 1}_{\mathcal{B}/\mathcal{A}}$ is a generator of the ideal $I$, and $I^n = \widehat{\on{dR}}^{\geq n}_{\mathcal{B}/\mathcal{A}}$.

 The claims on $\mathcal{R}$ are immediate from the change of bases formula of Hodge completed derived Rham complexes in lemma \ref{2830948302}.  
\end{proof}

\begin{lemma} \label{qrehjglsg}
Let $\mathcal{A} \longrightarrow  \mathcal{B}$ be a map of static analytic animated rings over $\mathbb{Q}$ such that $\mathbb{L}_{\mathcal{B}/\mathcal{A}} \cong \Sigma \mathcal{B}^{\triangleright}$ as an element of $D(\mathcal{B})$, then there exists a filtration on $\on{HH}(\mathcal{B}/\mathcal{A})^{t\mathbb{S}^1}$ with graded pieces 
\[
\on{gr}^{\bullet} \on{HH}(\mathcal{B}/\mathcal{A})^{t\mathbb{S}^1} \cong \pi_{2\bullet} (\on{HH}(\mathcal{B}/\mathcal{A})^{t\mathbb{S}^1})
\cong
\widehat{\on{dR}}_{\mathcal{B}/\mathcal{A}}[\sigma^{\pm 1}] \]
where $\sigma$ is a generator of $\Sigma\mathbb{L}_{\mathcal{B}/\mathcal{A}}$ 
with graded degree 1 and homological degree 2. Moreover, there exists a filtration on $\on{HH}(\mathcal{B}/\mathcal{A})^{h\mathbb{S}^1}$ with graded pieces
\[
\on{gr}^{\bullet}\on{HH}(\mathcal{B}/\mathcal{A})^{h\mathbb{S}^1} = \pi_{2\bullet} (\on{HH}(\mathcal{B}/\mathcal{A})^{h\mathbb{S}^1})
= \widehat{\on{dR}}_{\mathcal{B}/\mathcal{A}}[\Tilde{u}, v]/(\Tilde{u}v- t)
\]
where $\Tilde{u}$ has graded degree 1 and homological degree 2,  $v$ has graded degree $-1$ and homological degree $-2$, and $t$ is a generator of the ideal $I$ of $\widehat{\on{dR}}_{\mathcal{B}/\mathcal{A}}$.

Let $\mathcal{R}$ be an analytic animated ring over $\mathcal{B}$, then there exists filtration on $\on{HH}(\mathcal{R}/\mathcal{A})^{t\mathbb{S}^1}$
and $\on{HH}(\mathcal{R}/\mathcal{A})^{h\mathbb{S}^1}$
with graded pieces
\[
\on{gr}^{i}
\on{HH}(\mathcal{R}/\mathcal{A})^{t\mathbb{S}^1}
= 
\Sigma^{2i} \widehat{\on{dR}}_{\mathcal{R}/\mathcal{A}}
\]
\[
\on{gr}^{i}
\on{HH}(\mathcal{R}/\mathcal{A})^{h\mathbb{S}^1}
= \Sigma^{2i}\widehat{\on{dR}}_{\mathcal{R}/\mathcal{A}}^{\geq i}
\]
which are algebras over $\on{gr}^{\bullet}\on{HH}(\mathcal{B}/\mathcal{A})^{t\mathbb{S}^1} $
and 
$\on{gr}^{\bullet}\on{HH}(\mathcal{B}/\mathcal{A})^{h\mathbb{S}^1} $, and
the multiplication by $v$ on $\on{gr}^{\bullet}
\on{HH}(\mathcal{R}/\mathcal{A})^{h\mathbb{S}^1}$
induces the filtration map $\widehat{\on{dR}}_{\mathcal{R}/\mathcal{A}}^{\geq i} \longrightarrow \widehat{\on{dR}}_{\mathcal{R}/\mathcal{A}}^{\geq i-1}$, and the multiplication map by $\Tilde{u}$ induces the map $\widehat{\on{dR}}_{\mathcal{R}/\mathcal{A}}^{\geq i} \longrightarrow \widehat{\on{dR}}_{\mathcal{R}/\mathcal{A}}^{\geq i+1} $ 
which is multiplication by $t \in \widehat{\on{dR}}_{\mathcal{B}/\mathcal{A}}^{\geq 1}$. Moreover, we have canonical equivalence
\[
\on{HH}(\mathcal{R}/\mathcal{A})^{t\mathbb{S}^1} 
\otimes_{\on{HH}(\mathcal{B}/\mathcal{A})^{t\mathbb{S}^1}} 
\mathcal{B}^{\triangleright, t\mathbb{S}^1}
\cong 
\on{HH}(\mathcal{R}/\mathcal{B})^{t\mathbb{S}^1} 
\]
\[
\on{HH}(\mathcal{R}/\mathcal{A})^{h\mathbb{S}^1} 
\otimes_{\on{HH}(\mathcal{B}/\mathcal{A})^{h\mathbb{S}^1}} 
\mathcal{B}^{\triangleright, h\mathbb{S}^1}
\cong 
\on{HH}(\mathcal{R}/\mathcal{B})^{h\mathbb{S}^1} 
\]
of filtered objects. 

\end{lemma}

\begin{proof}
    This follows directly from theorem \ref{weirdzdfd}. The computation of $\on{gr}^{\bullet} \on{HH}(\mathcal{B}/\mathcal{A})^{t\mathbb{S}^1} $
    is theorem \ref{weirdzdfd} and the commutative algebra structure is as in remark \ref{xueoiwuroew}. The $\sigma$ is chosen as in remark \ref{xueoiwuroew}. It is also the homotopy groups since we know from lemma \ref{2374032493}
    that $\widehat{\on{dR}}_{\mathcal{B}/\mathcal{A}}$
    is static. 
    
    We take $\Tilde{u}$ to be a generator $t$ of $I$ shifted by homological degree 2, and $v$ to be the image of $\sigma$ under the canonical map from homotopy fixed points to Tate construction, then $\Tilde{u}v=t$ by construction since multiplication by $v$ induces the filtration map $\widehat{\on{dR}}_{\mathcal{B}/\mathcal{A}}^{\geq 1} \longrightarrow \widehat{\on{dR}}_{\mathcal{B}/\mathcal{A}}$ by construction of $v$ and the description of the natural map from homotopy fixed points to Tate construction as specified in theorem \ref{weirdzdfd}. The claim on homotopy fixed points is then also clear. The claims on $\mathcal{R}$ is straightforward. 

    For the last claim, we note that there is a filtered $\mathbb{T}_{\on{Fil}}$-equivariant  isomorphism 
    \[
\on{HH}(\mathcal{R}/\mathcal{A}) 
\otimes_{\on{HH}(\mathcal{B}/\mathcal{A})} 
\mathcal{B}^{\triangleright}
\cong 
\on{HH}(\mathcal{R}/\mathcal{B})
\]
where we give $\mathcal{B}^{\triangleright}$ the trivial filtration, namely its only non-zero graded piece is at degree 0. There is a natural map from the left to right, and we check it is an equivalence by looking at the graded pieces, which is 
\[
\on{Sym}^{\bullet}_{\mathcal{R}}
(\Sigma \mathbb{L}_{\mathcal{R}/\mathcal{A}} [1]) \otimes_{\on{Sym}^{\bullet}_{\mathcal{B}}
(\Sigma \mathbb{L}_{\mathcal{B}/\mathcal{A}} [1])} 
\mathcal{B}^{\triangleright} 
\cong 
\on{Sym}^{\bullet}_{\mathcal{R}}
(\Sigma \mathbb{L}_{\mathcal{R}/\mathcal{B}} [1]).
\]
It holds by applying $\on{Sym}^{\bullet}_{\mathcal{R}}$ to the cofiber sequence 
\[
\mathbb{L}_{\mathcal{B}/\mathcal{A}} \otimes_{\mathcal{B}} \mathcal{R}
\longrightarrow \mathbb{L}_{\mathcal{R}/\mathcal{A}}
\longrightarrow 
\mathbb{L}_{\mathcal{R}/\mathcal{B}}
\]
where we use that $\on{Sym}^{\bullet}_{\mathcal{R}}$ is left adjoint, so maps colimits to colimits in algebras. 
We then observe that $\mathcal{B}^{\triangleright}$ 
is perfect as $\on{HH}(\mathcal{B}/\mathcal{A})$-modules, which follows from the computation 
$\on{gr}^{\bullet}_{\on{HKR}} \on{HH}(\mathcal{B}/\mathcal{A}) \cong \pi_{2\bullet } \on{HH}(\mathcal{B}/\mathcal{A}) 
\cong \mathcal{B}^{\triangleright} [u]$
for a generator of $\Sigma\mathbb{L}_{\mathcal{B}/\mathcal{A}}$. Since both the Tate construction and homotopy fixed points are exact lax symmetric monoidal functors, the last claim follows. 
\end{proof}

\begin{lemma} \label{242039htwo}
Let $\mathcal{A} \longrightarrow  \mathcal{B}$ be a map of static analytic animated rings over $\mathbb{Q}$ such that $\mathbb{L}_{\mathcal{B}/\mathcal{A}} \cong \Sigma \mathcal{B}^{\triangleright}$ as an element of $D(\mathcal{B})$, and $\Tilde{u}, v, \sigma$ and $t$ as in the previous lemma.  Let $u$ be the image of 
$\Tilde{u} $ under the natural map 
$\on{HH}(\mathcal{B}/\mathcal{A})^{h\mathbb{S}^1}
\longrightarrow 
\on{HH}(\mathcal{B}/\mathcal{A})$, which specifies a trivialization $ \mathbb{L}_{\mathcal{B}/\mathcal{A}} \cong \Sigma \mathcal{B}^{\triangleright}$, and 
\[
\on{gr}^{\bullet}_{\on{HKR}} \on{HH}(\mathcal{B}/\mathcal{A}) \cong \pi_{2\bullet } \on{HH}(\mathcal{B}/\mathcal{A}) 
\cong \mathcal{B}^{\triangleright} [u]. 
\]
Moreover, the limit map 
\[
\widehat{\on{dR}}_{\mathcal{B}/\mathcal{A}}[\Tilde{u}, v]/(\Tilde{u}v- t)
\cong 
\on{gr}^{\bullet}
\on{HH}(\mathcal{B}/\mathcal{A})^{h\mathbb{S}^1}
\longrightarrow 
\on{gr}^{\bullet}\on{HH}(\mathcal{B}/\mathcal{A})
\cong \mathcal{B}^{\triangleright} [u]
\]
is simply quotient out by $v$. 
\end{lemma}

\begin{proof}
    This is immediate from theorem \ref{weirdzdfd}. Indeed, we know 
  $\on{gr}^{\bullet}
\on{HH}(\mathcal{B}/\mathcal{A})^{h\mathbb{S}^1}
\longrightarrow 
\on{gr}^{\bullet}\on{HH}(\mathcal{B}/\mathcal{A})$
is the graded quotient map
\[
\Sigma^{2i}\widehat{\on{dR}}_{\mathcal{B}/\mathcal{A}}^{\geq i} \longrightarrow \Sigma^{2i}\on{gr}^{i}\widehat{\on{dR}}_{\mathcal{B}/\mathcal{A}}^{\geq \bullet}
= \on{Sym}^i_{\mathcal{B}} \Sigma \mathbb{L}_{\mathcal{B}/\mathcal{A}} 
= \on{Sym}^i_{\mathcal{B}} \Sigma^2 \mathcal{B}^{\triangleright},
\]
but this is the quotient by $v$ map since $v$ corresponds to the inclusion of the filtration as lemma \ref{qrehjglsg} shows.
\end{proof}

One should compare the next two lemmas with the construction of Frobenius descent of prismatic cohomology in \cite[proposition 11.15]{BMS2}.

\begin{lemma} \label{2f937r03fs}
Let $\mathcal{R}$ be an analytic animated ring over $\mathcal{B}$, then 
$ \on{gr}^{\bullet}_{\on{HKR}}\on{HH}(\mathcal{R}/\mathcal{A}) = \on{Sym}^{\bullet}_{\mathcal{R}}
\Sigma \mathbb{L}_{\mathcal{R}/\mathcal{A}}  $
is an algebra over $\on{gr}^{\bullet}_{\on{HKR}}\on{HH}(\mathcal{B}/\mathcal{A}) =  \mathcal{B}^{\triangleright} [u]$, which defines a filtered object  
\[ 
\mathcal{R}^{\triangleright} 
\overset{u}{\longrightarrow} 
\Sigma^{-2}\on{Sym}^{1}_{\mathcal{R}}
\Sigma \mathbb{L}_{\mathcal{R}/\mathcal{A}} \overset{u}{\longrightarrow}   
\Sigma^{-4}\on{Sym}^{2}_{\mathcal{R}}
\Sigma \mathbb{L}_{\mathcal{R}/\mathcal{A}}
\overset{u}{\longrightarrow} \cdots 
\]
where the transition map is multiplication by the generator of $ \on{Sym}^{1}_{\mathcal{B}}
\Sigma \mathbb{L}_{\mathcal{B}/\mathcal{A}} \cong \Sigma \mathbb{L}_{\mathcal{B}/\mathcal{A}}  $
specified by $u$. 
Let 
\[
\overline{\prism}_{\mathcal{R}/\widehat{\on{dR}}_{\mathcal{B}/\mathcal{A}}}
:=
\on{colim} ( 
\mathcal{R}^{\triangleright} 
\overset{u}{\longrightarrow} 
\Sigma^{-2}
\on{Sym}^{1}_{\mathcal{R}}
\Sigma \mathbb{L}_{\mathcal{R}/\mathcal{A}} \overset{u}{\longrightarrow}   
\Sigma^{-4}
\on{Sym}^{2}_{\mathcal{R}}
\Sigma \mathbb{L}_{\mathcal{R}/\mathcal{A}}
\overset{u}{\longrightarrow} \cdots),
\]
the filtration it inherits is called the Hodge-Tate filtration, i.e.
\[
F^{i}_{\on{HT}}\overline{\prism}_{\mathcal{R}/\widehat{\on{dR}}_{\mathcal{B}/\mathcal{A}}}
:= 
\Sigma^{-2i}
\on{Sym}^{i}_{\mathcal{R}}
\Sigma \mathbb{L}_{\mathcal{R}/\mathcal{A}} 
\]
and we have 
\[
\on{gr}^{\bullet}_{\on{HT}} 
\overline{\prism}_{\mathcal{R}/\widehat{\on{dR}}_{\mathcal{B}/\mathcal{A}}} 
\cong
\Sigma^{-2*}\on{Sym}_{\mathcal{R}}^{\bullet} (\Sigma \mathbb{L}_{\mathcal{R}/\mathcal{B}}[1])
\]
Moreover, we have
\[
\on{gr}_{\on{HKR}}^{\bullet}\on{HH}(\mathcal{R}/\mathcal{A})[u^{-1}]
\cong
\overline{\prism}_{\mathcal{R}/\widehat{\on{dR}}_{\mathcal{B}/\mathcal{A}}}[u, u^{-1}],
\]
and 
\[
\on{gr}_{\on{HKR}}^{i}\on{HH}(\mathcal{R}/\mathcal{A})
\cong
\Sigma^{2i}
F^{i}_{\on{HT}}\overline{\prism}_{\mathcal{R}/\widehat{\on{dR}}_{\mathcal{B}/\mathcal{A}}}.
\]
\end{lemma}

\begin{proof}
This is a standard computation. 
To compute the graded pieces of Hodge-Tate filtration, first note that the 0th graded piece is clear. For degree one, this is simply the cofiber sequence
\[
\mathbb{L}_{\mathcal{B}/\mathcal{A}} \otimes_{\mathcal{B}} \mathcal{R}
\longrightarrow \mathbb{L}_{\mathcal{R}/\mathcal{A}}
\longrightarrow 
\mathbb{L}_{\mathcal{R}/\mathcal{B}}
\]
shifted by $\Sigma$
together with the trivialization 
$\Sigma\mathbb{L}_{\mathcal{B}/\mathcal{A}} \cong \Sigma^2 \mathcal{B}^{\triangleright} $
specified by $u$. The higher degree terms are simply the formula 
\[
\on{Sym}^{\bullet}_{\mathcal{R}}
(\Sigma \mathbb{L}_{\mathcal{R}/\mathcal{A}} [1]) \otimes_{\on{Sym}^{\bullet}_{\mathcal{B}}
(\Sigma \mathbb{L}_{\mathcal{B}/\mathcal{A}} [1])} 
\mathcal{B}^{\triangleright} 
\cong 
\on{Sym}^{\bullet}_{\mathcal{R}}
(\Sigma \mathbb{L}_{\mathcal{R}/\mathcal{B}} [1]),
\]
as  in the proof of lemma \ref{qrehjglsg}, and the observation that the left hand side is the quotient out by $u$. 
The last claim follows directly from $ \on{gr}_{\on{HKR}}^{\bullet}\on{HH}(\mathcal{R}/\mathcal{A}) =
    \on{Sym}^{\bullet}_{\mathcal{R}}
(\Sigma \mathbb{L}_{\mathcal{R}/\mathcal{A}} [1])$. 
\end{proof}

\begin{lemma} \label{92385093283fsd}
The limit map $\on{HH}(\mathcal{B}/\mathcal{A})^{h\mathbb{S}^1}
\longrightarrow 
\on{HH}(\mathcal{B}/\mathcal{A})$
gives rise to
\[
\on{gr}^{\bullet}\on{HH}(\mathcal{R}/\mathcal{A})^{h\mathbb{S}^1}[\Tilde{u}^{-1}]
\longrightarrow 
\on{gr}^{\bullet}\on{HH}(\mathcal{R}/\mathcal{A})[u^{-1}]
\]
which induces a pushout diagram
\begin{equation} \label{cxweriw0cds}
\begin{tikzcd}
  \on{gr}^{\bullet}\on{HH}(\mathcal{R}/\mathcal{A})^{h\mathbb{S}^1}[\Tilde{u}^{-1}]
  \arrow[r]
  &
  \on{gr}^{\bullet}\on{HH}(\mathcal{R}/\mathcal{A})[u^{-1}]
  \arrow[dl, phantom, "\llcorner", very near start]
    \\
   \on{gr}^{\bullet}\on{HH}(\mathcal{B}/\mathcal{A})^{h\mathbb{S}^1}[\Tilde{u}^{-1}]
   \arrow[r]
   \arrow[u]
   &
    \on{gr}^{\bullet}\on{HH}(\mathcal{B}/\mathcal{A})[u^{-1}]
    \arrow[u]
\end{tikzcd}
\end{equation}
where the vertical maps are the functoriality of the map from $\mathcal{B}$ to $\mathcal{R}$, and the horizontal maps are limit maps of fixed points. 

Let 
\[
\widehat{\prism}_{\mathcal{R}/\widehat{\on{dR}}_{\mathcal{B}/\mathcal{A}}} := 
\on{colim}( \widehat{\on{dR}}_{\mathcal{R}/\mathcal{A}}^{\geq 0} \overset{t}{\longrightarrow} \widehat{\on{dR}}_{\mathcal{R}/\mathcal{A}}^{\geq 1} 
    \overset{t}{\longrightarrow} \widehat{\on{dR}}_{\mathcal{R}/\mathcal{A}}^{\geq 2}
    \overset{t}{\longrightarrow}
    \cdots ),
\] 
which is naturally an algebra over 
\[
\on{colim}( \widehat{\on{dR}}_{\mathcal{B}/\mathcal{A}}^{\geq 0} \overset{t}{\cong} \widehat{\on{dR}}_{\mathcal{B}/\mathcal{A}}^{\geq 1} 
    \overset{t}{\cong} \widehat{\on{dR}}_{\mathcal{B}/\mathcal{A}}^{\geq 2}
    \overset{t}{\cong}
    \cdots )
    \cong 
   \widehat{ \on{dR}}_{\mathcal{B}/\mathcal{A}}.
\]

 (1) We have a natural equivalence of $ \widehat{ \on{dR}}_{\mathcal{B}/\mathcal{A}}$-algebras
\[
\on{gr}^{\bullet}\on{HH}(\mathcal{R}/\mathcal{A})^{h\mathbb{S}^1}[\Tilde{u}^{-1}] 
\cong
\widehat{\prism}_{\mathcal{R}/\widehat{\on{dR}}_{\mathcal{B}/\mathcal{A}}} [\Tilde{u},\Tilde{u}^{-1}],  
\]
where the $\widehat{ \on{dR}}_{\mathcal{B}/\mathcal{A}}$-algebra structure on the left hand side is through its algebra structure over
\[
\on{gr}^{\bullet}\on{HH}(\mathcal{B}/\mathcal{A})^{h\mathbb{S}^1}[\Tilde{u}^{-1}]
\cong 
\widehat{ \on{dR}}_{\mathcal{B}/\mathcal{A}}[\Tilde{u}, v, \Tilde{u}^{-1}]/(\Tilde{u}v-t )
\cong 
\widehat{ \on{dR}}_{\mathcal{B}/\mathcal{A}}[\Tilde{u}, \Tilde{u}^{-1}],
\]
and
$\widehat{\prism}_{\mathcal{R}/\widehat{\on{dR}}_{\mathcal{B}/\mathcal{A}}} $
sits in degree 0. 

(2) The pushout diagram (\ref{cxweriw0cds}) gives rise to an equivalence
\[ \widehat{\prism}_{\mathcal{R}/\widehat{\on{dR}}_{\mathcal{B}/\mathcal{A}}} \otimes_{\widehat{\on{dR}}_{\mathcal{B}/\mathcal{A}}} \widehat{\on{dR}}_{\mathcal{B}/\mathcal{A}}/I
\cong 
\overline{\prism}_{\mathcal{R}/\widehat{\on{dR}}_{\mathcal{B}/\mathcal{A}}},
\]
with 
$\overline{\prism}_{\mathcal{R}/\widehat{\on{dR}}_{\mathcal{B}/\mathcal{A}}}$
defined as in lemma \ref{2f937r03fs}.

(3) The localization map 
\[
\on{HH}(\mathcal{R}/\mathcal{A})^{h\mathbb{S}^1}
\longrightarrow
\on{HH}(\mathcal{R}/\mathcal{A})^{h\mathbb{S}^1}[\Tilde{u}^{-1}]
\]
gives rise to a map of $\widehat{\on{dR}}_{\mathcal{B}/\mathcal{A}} $-algebras 
\[
\widehat{\on{dR}}_{\mathcal{R}/\mathcal{A}}
\longrightarrow 
\widehat{\prism}_{\mathcal{R}/\widehat{\on{dR}}_{\mathcal{B}/\mathcal{A}}}
\]
on the 0-th graded pieces, which is the natural map from the first element $\widehat{\on{dR}}_{\mathcal{R}/\mathcal{A}}^{\geq 0}$ 
of the colimit diagram defining $\widehat{\prism}_{\mathcal{R}/\widehat{\on{dR}}_{\mathcal{B}/\mathcal{A}}} $
into the colimit. 
It induces an equivalence after inverting $t$, i.e. 
\[
\widehat{\on{dR}}_{\mathcal{R}/\mathcal{A}}[t^{-1}]
\cong 
\widehat{\prism}_{\mathcal{R}/\widehat{\on{dR}}_{\mathcal{B}/\mathcal{A}}}[t^{-1}].
\]
\end{lemma}

\begin{proof}
We prove the diagram is a pushout diagram before inverting $\Tilde{u}$. We note first that the lower horizontal arrow
$\on{gr}^{\bullet}\on{HH}(\mathcal{B}/\mathcal{A})^{h\mathbb{S}^1}
   \longrightarrow
    \on{gr}^{\bullet}\on{HH}(\mathcal{B}/\mathcal{A})$
is the quotient by $v$ map
\[
\widehat{\on{dR}}_{\mathcal{B}/\mathcal{A}}[\Tilde{u}, v]/(\Tilde{u}v-  t)
\longrightarrow 
\mathcal{B}^{\triangleright} [u]
\]
by lemma \ref{242039htwo}. The upper horizontal map 
$ \on{gr}^{\bullet}\on{HH}(\mathcal{R}/\mathcal{A})^{h\mathbb{S}^1}
\longrightarrow
  \on{gr}^{\bullet}\on{HH}(\mathcal{R}/\mathcal{A})$
is the graded quotient map $\Sigma^{2i}\widehat{\on{dR}}_{\mathcal{R}/\mathcal{A}}^{\geq i} \longrightarrow \Sigma^{2i} \on{gr}^{i}\widehat{\on{dR}}_{\mathcal{R}/\mathcal{A}}^{\geq \bullet} $ 
by theorem \ref{weirdzdfd}. We know from lemma \ref{qrehjglsg} that $v$ corresponds to the (shifted) canonical filtration map $\widehat{\on{dR}}_{\mathcal{R}/\mathcal{A}}^{\geq i+1} \longrightarrow \widehat{\on{dR}}_{\mathcal{R}/\mathcal{A}}^{\geq i}$, it follows that quotient out by $v$ is exactly the graded quotient map, proving the claim.

    By definition
    \begin{align*}
    (\on{gr}^{\bullet}\on{HH}(\mathcal{R}/\mathcal{A})^{h\mathbb{S}^1}[\Tilde{u}^{-1}])^0
   & = 
    \on{colim}( \on{gr}^{0}\on{HH}(\mathcal{R}/\mathcal{A})^{h\mathbb{S}^1} \overset{\Tilde{u}}{\longrightarrow} 
    \Sigma^{-2} \on{gr}^{1}\on{HH}(\mathcal{R}/\mathcal{A})^{h\mathbb{S}^1} 
    \overset{\Tilde{u}}{\longrightarrow}
  \Sigma^{-4}\on{gr}^2  \cdots )
    \\
    &
    =\on{colim}( \widehat{\on{dR}}_{\mathcal{R}/\mathcal{A}}^{\geq 0} \overset{\Tilde{u}}{\longrightarrow}  \widehat{\on{dR}}_{\mathcal{R}/\mathcal{A}}^{\geq 1} 
    \overset{\Tilde{u}}{\longrightarrow}
    \widehat{\on{dR}}_{\mathcal{R}/\mathcal{A}}^{\geq 2}
    \overset{\Tilde{u}}{\longrightarrow}
    \cdots )
    \end{align*}
    using 
    $\on{gr}^{i}\on{HH}(\mathcal{R}/\mathcal{A})^{h\mathbb{S}^1} = 
    \Sigma^{2i}\widehat{\on{dR}}_{\mathcal{R}/\mathcal{A}}^{\geq i}$
    proved in theorem \ref{weirdzdfd},  but we see in lemma \ref{qrehjglsg} that   $\Tilde{u}$ corresponds to multiplication by $t$ under this identification, which proves (1). (2) follows  by looking at the degree 0 terms of the pushout diagram. (3) is also immediate from the definition. 
\end{proof}

\begin{remark}
As the notation suggests,
the algebra $\widehat{\on{dR}}_{\mathcal{B}/\mathcal{A}}$ 
plays the role of a prism and 
$\widehat{\prism}_{\mathcal{R}/\widehat{\on{dR}}_{\mathcal{B}/\mathcal{A}}}  $ 
is the analogue of prismatic cohomology of $\mathcal{R}$ with respect to the prism $\widehat{\on{dR}}_{\mathcal{B}/\mathcal{A}}$. Moreover, 
$\widehat{\on{dR}}_{\mathcal{R}/\mathcal{A}}$
is the analogue of "Frobenius" twisted prismatic cohomology of $\mathcal{R}$.  There is no Frobenius in this setting, and  $\widehat{\on{dR}}_{\mathcal{R}/\mathcal{A}}$
is not determined by $ \widehat{\prism}_{\mathcal{R}/\widehat{\on{dR}}_{\mathcal{B}/\mathcal{A}}}$. 

The HKR filtration used here is the analogue of even filtration, and it is proved in \cite{Evenfil} that under the above condition of cotangent complexes and if the cotangent complex of $\mathcal{R}$ over $\mathcal{B}$ is static, then the even filtration is exactly the HKR filtration. The statements here matches perfectly with \cite{BMS2}. 
\end{remark}

\begin{remark}
    The choice of $t$ corresponds to fixing a trivialization of $\mathbb{L}_{\mathcal{B}/\mathcal{A}}$, which is the analogue of Breuil-Kisin twists in prismatic cohomology. We can make the formulas above more canonical by keeping track of the  Breuil-Kisin twists.  
\end{remark}

\section{$\mathbb{B}_{\on{dR}}^+$-cohomology of solid Huber rings}

We first show that the Hodge completed derived de Rham complex of characteristic 0 perfectoid ring over solid $\mathbb{Z}$ is the usual de Rham period ring. This is a classical result observed by  Beilinson, Bhatt, Fargues, and Guo-Li (\cite{GuoLi2021}). We present a direct proof in the condensed setting making use of Raksit's universal property of Hodge completed derived de Rham complex  established in previous sections.

\begin{definition}
 Let $(K, K^+)$ be a perfectoid ring over $(\mathbb{Q}_p, \mathbb{Z}_p)$,
 we define the de Rham period ring of $(K, K^+)$ as 
    \[
    \mathbb{B}_{\on{dR}}^+(K) :=\underline{ W(K^{+, \flat})} [p^{-1}]^{\wedge}_{\xi} \in \on{CAlg} D((\mathbb{Q}_p, \mathbb{Z}_p)_{\blacksquare})^{\heartsuit},
    \]
    where we give $W(K^{+, \flat})$ the $(p, [\omega])$-adic topology with $\omega$ a pseudo-uniformizer of $K^{+, \flat}$, and $\underline{ W(K^{+, \flat})}$ is the condensed ring with respect to the topology. $\xi$ is a generator for the canonical map $ W(K^{+, \flat}) \longrightarrow K^+$, which is viewed as an element of the underlying set $\underline{ W(K^{+, \flat})}(*)$ of the condensed ring $\underline{ W(K^{+, \flat})}$,  and the localization and completion are taken in the category of condensed rings. 

    When $K$ is obvious from the context, we will write
    \[
    \on{B}_{\on{dR}}^+ :=  \mathbb{B}_{\on{dR}}^+(K) (*)
    \]
    for the underlying ring of the condensed ring $\mathbb{B}_{\on{dR}}^+(K) $, with the usual topology. 
\end{definition}

Note that $ \mathbb{B}_{\on{dR}}^+(K)$ depends only on $K$ and not on $K^+$. It is solid since adically complete rings are solid, and completion and localization preserves solid modules. 

\begin{proposition} \label{afhsjfheosdf}
    Let $(K, K^+)$ be a perfectoid ring over $(\mathbb{Q}_p, \mathbb{Z}_p)$, and $(K,K^+)_{\blacksquare}$ the corresponding analytic animated rings as defined in \cite[theorem 3.2.8.]{Andreychev2021}, then 
    \[
    \widehat{\on{dR}}_{(K,K^+)_{\blacksquare}/\mathbb{Z}_{\blacksquare}} 
    \cong
    \widehat{\on{dR}}_{(K,K^+)_{\blacksquare}/(\mathbb{Q}_p,\mathbb{Z}_p)_{\blacksquare}}
    \cong 
    \widehat{\on{dR}}_{(K,\mathbb{Z}_p)_{\blacksquare}/(\mathbb{Q}_p,\mathbb{Z}_p)_{\blacksquare}}
    \cong 
    \mathbb{B}_{\on{dR}}^+(K) 
    \in  \on{CAlg} D((\mathbb{Q}_p, \mathbb{Z}_p)_{\blacksquare})
    \]
    and the filtration 
    $\widehat{\on{dR}}_{(K,K^+)_{\blacksquare}/\mathbb{Z}_{\blacksquare}}^{\geq *}$ 
    corresponds to the $\xi$-adic filtration on $ \mathbb{B}_{\on{dR}}^+(K)$. 
\end{proposition}

\begin{proof}
    The first equivalence is clear since $(\mathbb{Q}_p, \mathbb{Z}_p)_{\blacksquare}$
    is idempotent over $\mathbb{Z}_{\blacksquare}$. It is well-known that $\mathbb{L}_{(K,K^+)_{\blacksquare}/(\mathbb{Q}_p,\mathbb{Z}_p)_{\blacksquare}} \cong \Sigma K $
    after trivializing the Breuil-Kisin twist,
    see \cite{Scholze2021SolidBdR}. Indeed, we can compute that 
    \begin{equation} \label{vc4r83902fhd}
    \mathbb{L}_{(K,\mathbb{Z}_p)_{\blacksquare}/(\mathbb{Q}_p,\mathbb{Z}_p)_{\blacksquare}} \cong 
     \mathbb{L}_{(K^+,\mathbb{Z}_p)_{\blacksquare}/\mathbb{Z}_{p, \blacksquare}}[p^{-1}] 
     \cong 
    \mathbb{L}_{(K^+,\mathbb{Z}_p)_{\blacksquare}/(W(K^{+, \flat}), \mathbb{Z}_p)_{\blacksquare}}[p^{-1}]
    \cong \Sigma \underline{(\xi)/(\xi^2)}
    \end{equation}
    where the second isomorphism we use that $\mathbb{L}_{(W(K^{+, \flat}),\mathbb{Z}_p)_{\blacksquare}/\mathbb{Z}_{p, \blacksquare}}= 0 $, and this is because it is already $p$-complete, which is a miracle for solid category,  so modulo $p$ we are reduced to that the cotangent complex of a perfect ring is 0. This also tells us that $\mathbb{L}_{(K,K^+)_{\blacksquare}/(\mathbb{Q}_p,\mathbb{Z}_p)_{\blacksquare}} \cong \mathbb{L}_{(K,\mathbb{Z}_p)_{\blacksquare}/(\mathbb{Q}_p,\mathbb{Z}_p)_{\blacksquare}}$
    because the right hand side is already complete with respect to $(K,K^+)_{\blacksquare}$ as the computation shows it is equivalent to $\Sigma K$.

    The computation tells us also that the universal differential $d: K \longrightarrow \mathbb{L}_{(K,K^+)_{\blacksquare}/(\mathbb{Q}_p,\mathbb{Z}_p)_{\blacksquare}}\cong \Sigma \underline{(\xi)/(\xi^2)}$
    corresponds to the extension 
    \begin{equation} \label{afdweore}
    0 \longrightarrow 
(\xi)/(\xi^2)
    \longrightarrow
    W(K^{+, \flat})[p^{-1}]/\xi^2 \longrightarrow K \longrightarrow 0. 
     \end{equation}

    Let us view $\mathbb{B}_{\on{dR}}^+(K)$ with $\xi$-adic filtration as an element of $\on{Fil}^{\wedge}D((\mathbb{Q}_p,\mathbb{Z}_p)_{\blacksquare}) $, which defines an element of $ \on{cMod}_{\mathbb{D}_+^{\vee}} \on{Gr} D((\mathbb{Q}_p,\mathbb{Z}_p)_{\blacksquare}) $
    under the equivalence (\ref{zcbksueife}). It is also a commutative algebra object since the equivalence (\ref{zcbksueife}) is symmetric monoidal, so it defines an element of $\on{cMod}_{\mathbb{D}_+^{\vee}} (\on{DAlg} \on{Cond} \on{Gr} D(\mathbb{Z})_{\mathbb{Q}_p/})  $
    since derived algebra objects coincides with commutative algebra objects in characteristic 0. The universal property in proposition \ref{423087403284} induces a map 
    \[
    L\Omega^{+,\bullet}_{\underline{K}/\underline{ \mathbb{Q}_p}} \longrightarrow 
    \mathbb{B}_{\on{dR}}^+(K)
    \]
    of $\on{cMod}_{\mathbb{D}_+^{\vee}} (\on{DAlg} \on{Cond} \on{Gr} D(\mathbb{Z})_{\mathbb{Q}_p/})$
    which on degree 0 is the identity map of $K$. Give $ \mathbb{B}_{\on{dR}}^+(K)$ the analytic graded structure induced from its degree 0 part $K$ and its analytic structure $(K, K^+)_{\blacksquare}$ using proposition \ref{wury9387r9}, the functoriality as in  
   lemma \ref{wityuoiewuroe} gives rise to a map 
     \[
    L\Omega^{+,\bullet}_{(K,K^+)_{\blacksquare}/(\mathbb{Q}_p, \mathbb{Z}_p)_{\blacksquare}} \longrightarrow 
    \mathbb{B}_{\on{dR}}^+(K).
    \]
Indeed,  lemma \ref{wityuoiewuroe} induces a map $\widehat{\mathbb{DR}}^{\on{un}}_{(K,K^+)_{\blacksquare}/(\mathbb{Q}_p, \mathbb{Z}_p)_{\blacksquare}} \longrightarrow \mathbb{B}_{\on{dR}}^+(K)$
in $\on{GrAnDAG}^{\on{un}}$, passing to completion gives rise to 
$\widehat{\mathbb{DR}}_{(K,K^+)_{\blacksquare}/(\mathbb{Q}_p, \mathbb{Z}_p)_{\blacksquare}} \longrightarrow \mathbb{B}_{\on{dR}}^+(K)$, where we use that 
$\mathbb{B}_{\on{dR}}^+(K)$ is already complete as analytic graded derived ring since its underlying graded ring is $\on{Sym}^{\bullet}_{K} (\Sigma^2 K[1])$ that lies in $\on{Gr} D((K,K^+)_{\blacksquare})$. Taking the underlying condensed ring gives the desired map. 
    
    Since the underlying graded ring of both sides coincides with $\on{Sym}^{\bullet}_{K} (\Sigma^2 K[1])$ by lemma \ref{4702374903} and the definition of the equivalence (\ref{zcbksueife}), to prove it is an isomorphism it is enough to prove that it induces an isomorphism of degree 1, this is the same as in the proof proposition \ref{3074375372} where we observe that the $\mathbb{D}_+$-module structures on 1-truncated parts of both sides give rise to a commutative diagram 
    \[
    \begin{tikzcd}
     \mathbb{L}_{(K,K^+)_{\blacksquare}/(\mathbb{Q}_p, \mathbb{Z}_p)_{\blacksquare}} 
     \arrow[r]
     &
     \Sigma (\xi)/(\xi^2)
        \\
     K \arrow[u,"d"] \arrow[r,equal]  & K \arrow[u]
    \end{tikzcd}
    \]
    we have observed that the left hand side corresponds to the extension (\ref{afdweore}) and it follows directly from the definition that the right hand side also corresponds to the same extension, so the horizontal arrow is an isomorphism by the universal property of cotangent complex. 

    Lastly, translating everything from comodules over $\mathbb{D}_{+}^{\vee}$ to completely filtered objects as in definition \ref{sfdsfjekl} gives the equivalence we want to prove. 
\end{proof}

We can now define two $\mathbb{B}_{\on{dR}}^+(K)$-cohomology theories of any analytic animated rings over perfectoid solid ring $(K,K^+)_{\blacksquare}$
with expected properties. 

\begin{definition}
    Let $(K, K^+)$ be a perfectoid ring over $(\mathbb{Q}_p, \mathbb{Z}_p)$, and $(K,K^+)_{\blacksquare} \longrightarrow \mathcal{R}$
    a map of analytic animated rings, then we define
    \[
    R\Gamma_{\on{inf}}(\mathcal{R}/\mathbb{B}_{\on{dR}}^+(K) ):= \widehat{\on{dR}}_{\mathcal{R}/ (\mathbb{Q}_p, \mathbb{Z}_p)_{\blacksquare}}
    \]
    with the filtration
    \[
    R\Gamma_{\on{inf}}(\mathcal{R}/\mathbb{B}_{\on{dR}}^+(K) )^{\geq \bullet}:= \widehat{\on{dR}}_{\mathcal{R}/ (\mathbb{Q}_p, \mathbb{Z}_p)_{\blacksquare}}^{\geq \bullet}
    \]
    as a $\mathbb{B}_{\on{dR}}^+(K) $-algebra in solid rings together with a filtration that is compatible with the $\xi$-adic filtration on $\mathbb{B}_{\on{dR}}^+(K) $. The  $\mathbb{B}_{\on{dR}}^+(K) $-algebra structure is through the natural $\widehat{\on{dR}}_{(K,K^+)_{\blacksquare}/ (\mathbb{Q}_p, \mathbb{Z}_p)_{\blacksquare}}$-algebra structure, which is equivalent to $\mathbb{B}_{\on{dR}}^+(K)$ by  proposition \ref{afhsjfheosdf}. 
\end{definition}

The notation is meant as a generalization of the infinitesmial $\on{B}_{\on{dR}}^+$-cohomology of \cite{Guo2025Crystalline}, and we will make it precise soon.  

\begin{proposition}
   Let $(K, K^+)$ be a perfectoid ring over $(\mathbb{Q}_p, \mathbb{Z}_p)$, and $(K,K^+)_{\blacksquare} \longrightarrow \mathcal{R}$
    a map of analytic animated rings,   there exists filtration on $\on{HH}((K, K^+)_{\blacksquare}/(\mathbb{Q}_p, \mathbb{Z}_p)_{\blacksquare})^{t\mathbb{S}^1}$
and $\on{HH}((K, K^+)_{\blacksquare}/(\mathbb{Q}_p, \mathbb{Z}_p)_{\blacksquare})^{h\mathbb{S}^1}$
with graded pieces
    \[
\on{gr}^{\bullet}
\on{HH}((K, K^+)_{\blacksquare}/(\mathbb{Q}_p, \mathbb{Z}_p)_{\blacksquare})^{t\mathbb{S}^1}
\cong
\mathbb{B}_{\on{dR}}^+(K) [\sigma,\sigma^{-1}]
\]
with $\sigma$ of graded degree 1 and homological degree 2, and 
\[
\on{gr}^{\bullet}
\on{HH}((K, K^+)_{\blacksquare}/(\mathbb{Q}_p, \mathbb{Z}_p)_{\blacksquare})^{h\mathbb{S}^1}
\cong \mathbb{B}_{\on{dR}}^+(K) [\Tilde{u},v]/(\Tilde{u}v-\xi)
\]
with $\Tilde{u}$ of homological degree 2 and graded degree 1. 
Moreover, the limit map 
\[
\mathbb{B}_{\on{dR}}^+(K)[\Tilde{u}, v]/(\Tilde{u}v- \xi)
\cong 
\on{gr}^{\bullet}
\on{HH}((K, K^+)_{\blacksquare}/(\mathbb{Q}_p, \mathbb{Z}_p)_{\blacksquare})^{h\mathbb{S}^1}
\longrightarrow 
\on{gr}^{\bullet}\on{HH}((K, K^+)_{\blacksquare}/(\mathbb{Q}_p, \mathbb{Z}_p)_{\blacksquare})
\cong K [u]
\]
with $u$ the image of $\Tilde{u}$,
is simply quotient out by $v$. 
\end{proposition}

\begin{proof}
    Lemma \ref{242039htwo} and (\ref{vc4r83902fhd}). 
\end{proof}

\begin{theorem}
    Let $(K, K^+)$ be a perfectoid ring over $(\mathbb{Q}_p, \mathbb{Z}_p)$, and $(K,K^+)_{\blacksquare} \longrightarrow \mathcal{R}$
    a map of analytic animated rings, then we have the de Rham comparison
    \[
     R\Gamma_{\on{inf}}(\mathcal{R}/\mathbb{B}_{\on{dR}}^+(K) ) \otimes_{(\mathbb{B}_{\on{dR}}^+(K),\mathbb{Z}_p)_{\blacksquare}  } K 
     \cong 
     \widehat{\on{dR}}_{\mathcal{R}/ (K, K^+)_{\blacksquare}}
    \]

    Moreover,  there exists filtration on $\on{HH}(\mathcal{R}/(\mathbb{Q}_p, \mathbb{Z}_p)_{\blacksquare})^{t\mathbb{S}^1}$
and $\on{HH}(\mathcal{R}/(\mathbb{Q}_p, \mathbb{Z}_p)_{\blacksquare})^{h\mathbb{S}^1}$
with graded pieces
\[
\on{gr}^{i}
\on{HH}(\mathcal{R}/(\mathbb{Q}_p, \mathbb{Z}_p)_{\blacksquare})^{t\mathbb{S}^1}
= 
\Sigma^{2i}  R\Gamma_{\on{inf}}(\mathcal{R}/\mathbb{B}_{\on{dR}}^+(K) ) 
\]
\[
\on{gr}^{i}
\on{HH}(\mathcal{R}/(\mathbb{Q}_p, \mathbb{Z}_p)_{\blacksquare})^{h\mathbb{S}^1}
= \Sigma^{2i} R\Gamma_{\on{inf}}(\mathcal{R}/\mathbb{B}_{\on{dR}}^+(K) )^{\geq i}
\]
such that the canonical map from homotopy fixed points to Tate construction induces the filtration maps on graded pieces. 
\end{theorem}

\begin{proof}
The de Rham comparison is lemma \ref{2830948302}, and the rest is lemma \ref{qrehjglsg}.
\end{proof}

\begin{definition}
 Let $(K, K^+)$ be a perfectoid ring over $(\mathbb{Q}_p, \mathbb{Z}_p)$, and $(K,K^+)_{\blacksquare} \longrightarrow \mathcal{R}$
    a map of analytic animated rings, then we define 
    \[
    \widehat{\prism}_{\mathcal{R}/\mathbb{B}_{\on{dR}}^+(K)   }: = \on{colim}( \widehat{\on{dR}}_{\mathcal{R}/(\mathbb{Q}_p, \mathbb{Z}_p)_{\blacksquare}}^{\geq 0} \overset{\xi}{\longrightarrow} \widehat{\on{dR}}_{\mathcal{R}/(\mathbb{Q}_p, \mathbb{Z}_p)_{\blacksquare}}^{\geq 1} 
    \overset{\xi}{\longrightarrow} \widehat{\on{dR}}_{\mathcal{R}/(\mathbb{Q}_p, \mathbb{Z}_p)_{\blacksquare}}^{\geq 2}
    \overset{\xi}{\longrightarrow}
    \cdots ),
    \]
    as in lemma \ref{92385093283fsd}, which is a $\mathbb{B}_{\on{dR}}^+(K) $-algebra in solid modules,  and 
    \[
    \overline{\prism}_{\mathcal{R}/\mathbb{B}_{\on{dR}}^+(K)   } 
    := \widehat{\prism}_{\mathcal{R}/\mathbb{B}_{\on{dR}}^+(K)   } \otimes_{ (\mathbb{B}_{\on{dR}}^+(K), \mathbb{Z}_p )_{\blacksquare}} K. 
    \]
\end{definition}

\begin{theorem}[Hodge-Tate filtration]
     Let $(K, K^+)$ be a perfectoid ring over $(\mathbb{Q}_p, \mathbb{Z}_p)$, and $(K,K^+)_{\blacksquare} \longrightarrow \mathcal{R}$
    a map of analytic animated rings, then there exists a natural Hodge-Tate filtration
    \[
F_{\on{HT}}^{i}\overline{\prism}_{\mathcal{R}/\mathbb{B}_{\on{dR}}^+(K)   }
    \cong 
     \Sigma^{-2i}
    \on{Sym}^i_{\mathcal{R}} \Sigma \mathbb{L}_{\mathcal{R}/(\mathbb{Q}_p, \mathbb{Z}_p)_{\blacksquare} }
    \]
    on $\overline{\prism}_{\mathcal{R}/\mathbb{B}_{\on{dR}}^+(K)   } $
    with graded pieces
    \[
    \on{gr}_{\on{HT}}^{\bullet} \overline{\prism}_{\mathcal{R}/\mathbb{B}_{\on{dR}}^+(K)   }  
    \cong 
    \Sigma^{-2*}\on{Sym}_{\mathcal{R}}^{\bullet} (\Sigma \mathbb{L}_{\mathcal{R}/(K,K^+)_{\blacksquare}}[1])
    \]

  Moreover, we have
\[
\on{gr}_{\on{HKR}}^{\bullet}\on{HH}(\mathcal{R}/(\mathbb{Q}_p, \mathbb{Z}_p)_{\blacksquare})[u^{-1}]
\cong
\overline{\prism}_{\mathcal{R}/\mathbb{B}_{\on{dR}}^+(K)  } [u, u^{-1}]
\]
with $u$ in homological degree 2 and graded degree 1,
and 
\begin{equation} \label{30r538fsd}
\on{gr}_{\on{HKR}}^{i}\on{HH}(\mathcal{R}/(\mathbb{Q}_p, \mathbb{Z}_p)_{\blacksquare})
\cong
\Sigma^{2i}
F^{i}_{\on{HT}}\overline{\prism}_{\mathcal{R}/\mathbb{B}_{\on{dR}}^+(K)  },
\end{equation}
and $u$ acts by the filtration map of Hodge-Tate filtration. Lastly, the limit map 
\[
\on{HH}(\mathcal{R}/(\mathbb{Q}_p, \mathbb{Z}_p)_{\blacksquare})^{h\mathbb{S}^1}
\longrightarrow
\on{HH}(\mathcal{R}/(\mathbb{Q}_p, \mathbb{Z}_p)_{\blacksquare})
\]
induces on graded objects
\begin{equation}  \label{cnwerwidjfi}
\on{gr}^i R\Gamma_{\on{inf}}(\mathcal{R}/\mathbb{B}_{\on{dR}}^+(K) )
\cong 
F^{i}_{\on{HT}}\overline{\prism}_{\mathcal{R}/\mathbb{B}_{\on{dR}}^+(K)  }
\end{equation}
\end{theorem}

\begin{proof}
    This is lemma \ref{2f937r03fs} and lemma \ref{92385093283fsd}. 
\end{proof}

\begin{remark}
    One should interpret the infinitesmial cohomology $R\Gamma_{\on{inf}}(\mathcal{R}/\mathbb{B}_{\on{dR}}^+(K) )$
    as the analogue of Frobenius twisted
    prismatic cohomology, and the infinitesmial filtration as the analogue of the Nygaard filtration. The equivalences (\ref{30r538fsd}) and (\ref{cnwerwidjfi}) are then the analogue of Segal conjecture in this case. 
\end{remark}

\begin{theorem}
Let $(K, K^+)$ be a perfectoid ring over $(\mathbb{Q}_p, \mathbb{Z}_p)$, and $(K,K^+)_{\blacksquare} \longrightarrow \mathcal{R}$
    a map of analytic animated rings, then there exists a canonical equivalence
    \[
    \on{gr}^{\bullet}
\on{HH}(\mathcal{R}/(\mathbb{Q}_p, \mathbb{Z}_p)_{\blacksquare})^{h\mathbb{S}^1}[\Tilde{u}^{-1}]
\cong 
\widehat{\prism}_{\mathcal{R}/\mathbb{B}_{\on{dR}}^+(K)   } [ \Tilde{u}, \Tilde{u}^{-1}]
    \]
with $\Tilde{u}$ of homological degree 2 and graded degree 1, and the natural limit map 
\[
\on{HH}
(\mathcal{R}/(\mathbb{Q}_p, \mathbb{Z}_p)_{\blacksquare})^{h\mathbb{S}^1} 
\longrightarrow 
\on{HH}(\mathcal{R}/(\mathbb{Q}_p, \mathbb{Z}_p)_{\blacksquare})
\]
induces after inverting $\Tilde{u}$
the quotient map 
$\widehat{\prism}_{\mathcal{R}/\mathbb{B}_{\on{dR}}^+(K)   } 
\longrightarrow \overline{\prism}_{\mathcal{R}/\mathbb{B}_{\on{dR}}^+(K)   } $
on  0-th graded. 
  
   Moreover, the localization map 
   \[
   \on{HH}(\mathcal{R}/(\mathbb{Q}_p, \mathbb{Z}_p)_{\blacksquare})^{h\mathbb{S}^1}
   \longrightarrow
   \on{HH}(\mathcal{R}/(\mathbb{Q}_p, \mathbb{Z}_p)_{\blacksquare})^{h\mathbb{S}^1}[\Tilde{u}^{-1}]
   \]
   induces a comparison map 
   \[
     R\Gamma_{\on{inf}}(\mathcal{R}/\mathbb{B}_{\on{dR}}^+(K) ) 
     \longrightarrow 
     \widehat{\prism}_{\mathcal{R}/\mathbb{B}_{\on{dR}}^+(K)   }
   \]
   and it becomes an equivalence after inverting $\xi$, i.e. 
   \[
     R\Gamma_{\on{inf}}(\mathcal{R}/\mathbb{B}_{\on{dR}}^+(K) ) [\xi^{-1}]
     \cong 
     \widehat{\prism}_{\mathcal{R}/\mathbb{B}_{\on{dR}}^+(K)   }[\xi^{-1}].
   \]
\end{theorem}

\begin{proof}
  This is   lemma \ref{92385093283fsd}. 
\end{proof}

\begin{remark}
Throughout the discussion above,
     we have implicitly fixed  generator $\xi$ and $\Tilde{u} $
    to trivialize the Breuil-Kisin twists.
\end{remark}

\section{Comparison with Guo's Crystalline $\on{B}_{\on{dR}}^+$-cohomology}

In this section, we compare the cohomology theory with crystalline $\on{B}_{\on{dR}}^+$-cohomology constructed in \cite{Guo2025Crystalline} in smooth cases. This cohomology theory is first construct by hand in \cite[section 13]{Bhatt2018}, and we will follow mostly the treatment \cite{Guo2025Crystalline}.  

\begin{lemma} \label{weruoincoe35}
Let $(K, K^+)$ be a perfectoid ring over $(\mathbb{Q}_p, \mathbb{Z}_p)$, then 
\[
\mathbb{L}_{(\mathbb{B}_{\on{dR}}^+(K),  \mathbb{Z}_p)_{\blacksquare}/(\mathbb{Q}_p, \mathbb{Z}_p)_{\blacksquare}} 
\cong 0. 
\]
\end{lemma}

\begin{proof}
    We have already seen in the proof of proposition \ref{afhsjfheosdf} that 
    \[
    \mathbb{L}_{(W(K^{+, \flat})[p^{-1}], \mathbb{Z}_p)_{\blacksquare}/(\mathbb{Q}_p,\mathbb{Z}_{p})_{ \blacksquare}}= \mathbb{L}_{(W(K^{+, \flat}), \mathbb{Z}_p)_{\blacksquare}/\mathbb{Z}_{p, \blacksquare}}[p^{-1}]=
    0.
    \]
    We know from the lemma below that 
    $
   W(K^{+, \flat})[p^{-1}] \longrightarrow \mathbb{B}_{\on{dR}}^+(K)
    $
    is an idempotent algebra in solid rings, which forces
    \[
    \mathbb{L}_{(\mathbb{B}_{\on{dR}}^+(K) , \mathbb{Z}_p)_{\blacksquare}/   (W(K^{+, \flat})[p^{-1}],  \mathbb{Z}_p)_{\blacksquare}
  } = 0,
  \]
  and the lemma follows from these two observations. 
\end{proof}

\begin{lemma} \label{fhowufew0u9f0ew}
Let $(K, K^+)$ be a perfectoid ring over $(\mathbb{Q}_p, \mathbb{Z}_p)$, then there exists a canonical isomorophism
\[
 W(K^{+, \flat})[p^{-1}][[T]]
 \otimes_{  (W(K^{+, \flat})[p^{-1}], \mathbb{Z}_p)_{\blacksquare}}
 W(K^{+, \flat})[p^{-1}][[U]]
  \cong 
   W(K^{+, \flat})[p^{-1}][[T, U]],
\]
and the natural map 
\[
  W(K^{+, \flat})[p^{-1}] \longrightarrow \mathbb{B}_{\on{dR}}^+(K)
\]
is an ideompotent algebra in $D(\mathbb{Z}_{p, \blacksquare})$. 
\end{lemma}

\begin{proof}
    The second claim follows from the first claim. Indeed, we know the $\xi$-adic completion can be expressed as 
    \[
    \mathbb{B}_{\on{dR}}^+(K) =   W(K^{+, \flat})[p^{-1}]^{\wedge}_{\xi} =   W(K^{+, \flat})[p^{-1}][[T]]/(T-\xi)
    \]
    see the proof of \cite[Lemma 2.12.9]{Mann} for example. It follows that 
    \begin{align*}
    \mathbb{B}_{\on{dR}}^+(K)
     \otimes_{  (W(K^{+, \flat})[p^{-1}], \mathbb{Z}_p)_{\blacksquare}}
     \mathbb{B}_{\on{dR}}^+(K)
     &
     \cong 
     W(K^{+, \flat})[p^{-1}][[T]]
     \otimes_{  (W(K^{+, \flat})[p^{-1}], \mathbb{Z}_p)_{\blacksquare}}
     W(K^{+, \flat})[p^{-1}][[U]]/(T-\xi, U-\xi)
     \\
     & 
     \cong   W(K^{+, \flat})[p^{-1}][[T, U]] /(T-\xi, U-\xi)
     \\
     & \cong  \mathbb{B}_{\on{dR}}^+(K)
     \end{align*}
     where we use the first claim on the second isomorphism. 

     We now prove the first claim. Let $ \mathfrak{C} $ be the  set of closed subrings of $W(K^{+, \flat})$ such that its image in the quotient map $ W(K^{+, \flat}) \longrightarrow W(K^{+, \flat})/(p, [\omega])^n $
     is finite for all $n$.  Let  $\mathfrak{D}$ be the set of subrings of $W(K^{+, \flat})[p^{-1}]$ such that $ R \in \mathfrak{D}$ if and only if $R= C[p^{-1}]$ for some $C \in \mathfrak{C}$, and  $\mathfrak{D}$ is ordered by inclusion. 
     
     We claim that $\mathfrak{D}$ is a $\omega_1$-filtered set. Indeed, given countably many $R_i=C_i[p^{-1}] \in \mathfrak{D}$ with $C_i \in \mathfrak{C}$ and $i \geq 1$, let $A = \cup p^i C_i$ be the subset of $W(K^{+, \flat}) $ defined as the union of the subsets $p^i C_i$, and let 
     $C = \overline{\mathbb{Z}[A]}$ be the closed subring of $W(K^{+, \flat})  $ generated by $A$. More precisely, $C$ is the closure of the sub-$\mathbb{Z}$-algebra of $W(K^{+, \flat}) $ generated by the set $A$. By definition, we have $p^i C_i \subset C$, so $C_i[p^{-1}] \subset C[p^{-1}]$
     for every $i$. Moreover, we claim that $C \in \mathfrak{C}$. For every $n$, the image of $A$ in the quotient $W(K^{+, \flat})/(p, [\omega])^n$ is the union of the image of $ p^i C_i$ for $i \leq n$, so is finite since $C_i \in \mathfrak{C}$. The algebra it generates in $W(K^{+, \flat})/(p, [\omega])^n$ is also finite because the image of $ p^i C_i$ in $W(K^{+, \flat})/(p, [\omega])^n$ 
     is nilpotent  and $p^n$-torsion. It follows that the image of $\mathbb{Z}[A]$ in  $W(K^{+, \flat})/(p, [\omega])^n$ 
     is finite, which coincides with the image of $C$  by definition of closure. We have now seen that $C[p^{-1}] \in \mathfrak{D}$ and contains all $C_i[p^{-1}]$.

     We observe that 
     \begin{equation} \label{2374035325fg}
      W(K^{+, \flat})[p^{-1}] = \on{colim}_{R \in \mathfrak{D}} R
     \end{equation}
     as condensed rings. This is because every compact subset $ V $ of $ W(K^{+, \flat})[p^{-1}]$ is contained in $ p^{-i}   W(K^{+, \flat}) $ 
     for some $n$, and after multiplication by $p^n $ we can assume that it is contained in $  W(K^{+, \flat})$. Let $ C = \overline{\mathbb{Z}[p^{i+1} V]}$, the closed subalgebra of $ W(K^{+, \flat})$ generated by $p^{n+1} V$. Since the image of any compact subset of $  W(K^{+, \flat})$ in $W(K^{+, \flat})/(p, [\omega])^n$ must be finite, we see that the image of $ p^{i+1} V$ in 
     $W(K^{+, \flat})/(p, [\omega])^n$
 is finite, nilpotent, and $p^n$-torsion, so the image of $C$ is also finite, i.e. $C \in \mathfrak{C}$. Now $C[p^{-1}] \in \mathfrak{D}$ and $ V \subset  C[p^{-1}] $, proving the claim. 

 Moreover, we have that 
 \begin{equation} \label{zcveeuoro9}
 \prod_{\mathbb{N}}
  W(K^{+, \flat})[p^{-1}]  = \on{colim}_{R \in \mathfrak{D}} \prod_{\mathbb{N}}R. 
 \end{equation}
 By the last paragraph, every compact subset of $ \prod_{\mathbb{N}} W(K^{+, \flat})[p^{-1}]$
lies in 
$\prod_{i \in \mathbb{N}} R_i$ 
for some $R_i \in \mathfrak{D}$. Since $\mathfrak{D}$ is $\omega_1$-filtered, we can find $R\in \mathfrak{D}$ that contains all $R_i$, hence every compact subset of  $ \prod_{\mathbb{N}} W(K^{+, \flat})[p^{-1}]$ 
lies in $ \prod_{\mathbb{N}}R$ for some $R \in \mathfrak{D}$, proving the claim. 

 For every $C \in\mathfrak{C}$,  we have that $C = \on{lim}_n C_n$ with $C_n $  the image of $C$ in $W(K^{+, \flat})/(p, [\omega])^n$ since $C$  is closed, so $C$ is a profinite abelian group. By \cite[lemma 3.26]{Andreychev2021}, we have 
 \[
 C \otimes_{\mathbb{Z}_{\blacksquare}}  \prod_{\mathbb{N}}  \mathbb{Z}
 \cong 
  \prod_{\mathbb{N}} C,
 \]
 and it follows that 
 \begin{equation} \label{2hfoewhfo}
  C [[T]] \otimes_{(C, \mathbb{Z})_{\blacksquare}} C[[U]] 
  \cong 
  C \otimes_{\mathbb{Z}_{\blacksquare}} \mathbb{Z}[[T]] \otimes_{\mathbb{Z}_{\blacksquare}} \mathbb{Z}[[U]]
  \cong 
    C \otimes_{\mathbb{Z}_{\blacksquare} } \mathbb{Z}[[T,U]]
    \cong 
    C[[T,U]].  
 \end{equation}
Now the standard computation as in the proof of \cite[proposition A.66]{Bosco2021Drinfeld} gives rise to
\begin{align}\label{acowur4389}
C[p^{-1}] [[T]] \otimes_{(C[p^{-1}], \mathbb{Z})_{\blacksquare}} C[p^{-1}][[U]] 
 &
 \cong 
 (\on{colim}_{f:\mathbb{N} \rightarrow \mathbb{N}} \prod_{\mathbb{N}} p^{-f(i)} C)
 \otimes_{(C[p^{-1}], \mathbb{Z})_{\blacksquare}}
 \on{colim}_{g:\mathbb{N} \rightarrow \mathbb{N}} \prod_{\mathbb{N}} p^{-g(j)} C
 \\
 & 
 \cong 
 \on{colim}_{f,g:\mathbb{N} \rightarrow \mathbb{N}} \prod_{\mathbb{N}\times\mathbb{N}} p^{-f(i)-g(j)} C
 \\
  &
  \cong 
   C[p^{-1}] [[T,U]] 
\end{align} 
where the first isomorphism is the usual compact cover of the infinite product, and the function $f: \mathbb{N} \rightarrow \mathbb{N}$ are ordered by pointwise order, i.e. $f \leq g$ if and only if $f(n) \leq g(n)$ for every $n$. 
The second isomorphism uses (\ref{2hfoewhfo}), and the last isomorphism holds because for any function $\mathbb{N} \times \mathbb{N} \rightarrow \mathbb{N}$ we can find a bigger one of the form $(x,y) \rightarrow f(x)+g(y)$ for some $f, g: \mathbb{N} \rightarrow \mathbb{N}$, see the proof of \cite[proposition A.66]{Bosco2021Drinfeld}  for example. 

 We now compute
 \begin{align*}
 &
      W(K^{+, \flat})[p^{-1}][[T]]
 \otimes_{  (W(K^{+, \flat})[p^{-1}], \mathbb{Z}_p)_{\blacksquare}}
 W(K^{+, \flat})[p^{-1}][[U]] 
 \\
 \cong 
 &
 \on{colim}_{[n] \in \Delta} 
  W(K^{+, \flat})[p^{-1}][[T]] \otimes_{\mathbb{Z}_{\blacksquare}}  
   W(K^{+, \flat})[p^{-1}]^{\otimes_{\mathbb{Z}_{\blacksquare}}^n}
  \otimes_{\mathbb{Z}_{\blacksquare}}
  W(K^{+, \flat})[p^{-1}][[U]]
 \\
 \cong 
 &
 \on{colim}_{[n] \in \Delta} 
 (\on{colim}_{R \in \mathfrak{D}} R[[T]] )
 \otimes_{\mathbb{Z}_{\blacksquare}} 
  (\on{colim}_{R_1 \in \mathfrak{D}} R_1) 
  \otimes_{\mathbb{Z}_{\blacksquare}} \cdots  \otimes_{\mathbb{Z}_{\blacksquare}}  (\on{colim}_{R_n \in \mathfrak{D}} R_n)
   \otimes_{\mathbb{Z}_{\blacksquare}}
  (\on{colim}_{S \in \mathfrak{D}} S[[T]] )
 \\
 \cong 
 &
  \on{colim}_{[n] \in \Delta} 
  \on{colim}_{R \in \mathfrak{D}} 
 R[[T]] \otimes_{\mathbb{Z}_{\blacksquare}}  
   R^{\otimes_{\mathbb{Z}_{\blacksquare}}^n}
  \otimes_{\mathbb{Z}_{\blacksquare}}
  R[[U]]
  \\
  \cong 
  & 
    \on{colim}_{R \in \mathfrak{D}} 
    R[[T]] \otimes_{(R,\mathbb{Z})_{\blacksquare}} R[[U]]
    \\
    \cong
    & 
       \on{colim}_{R \in \mathfrak{D}} R[[T,U]]
       \\
       \cong 
       & 
         W(K^{+, \flat})[p^{-1}][[T,U]],
 \end{align*} 
 where the first isomorphism is the usual bar complex on relative tensor product, and the second is (\ref{2374035325fg}) and (\ref{zcveeuoro9}). In the third isomorphism we use that $\mathfrak{D}$ is filtered, namely it is sifted, so the diagonal embedding is cofinal, and we can replace the colimit indexed by $\mathfrak{D}^{n+2}$ by the diagonal $\mathfrak{D}$. The second to last isomorphism is (\ref{acowur4389}), and the last one is (\ref{zcveeuoro9}).
 \end{proof}

We also record the following lemma for later use. 

\begin{lemma} \label{vcoerionvcd}
    Let $(K, K^+)$ be a perfectoid ring over $(\mathbb{Q}_p, \mathbb{Z}_p)$, and $\mathcal{R}$ analytic ring  over $(K, K^+)_{\blacksquare}$,  then 
    \[
    \on{lim}_{e} \widehat{\on{dR}}_{\mathcal{R}/(\mathbb{B}_{\on{dR}}^+(K)/\xi^e, \mathbb{Z}_p)_{\blacksquare}} 
    \cong 
    \widehat{\on{dR}}_{\mathcal{R}/(\mathbb{B}_{\on{dR}}^+(K), \mathbb{Z}_p)_{\blacksquare}}
    \]
\end{lemma}

\begin{proof}
    There exists a canonical (filtered) morphism  from the right to the left, to prove it is an equivalence, it is enough to prove that it is an isomoprhism on the graded pieces. By (\ref{zdfhweuoofd}), it is enough to prove that 
    \[ \on{Sym}^i_{\mathcal{R}} \Sigma \mathbb{L}_{ \mathcal{R}/(\mathbb{B}_{\on{dR}}^+(K), \mathbb{Z}_p)_{\blacksquare} } 
    \cong 
    \on{lim}_e 
    \on{Sym}^i_{\mathcal{R}}
    \Sigma \mathbb{L}_{\mathcal{R}/(\mathbb{B}_{\on{dR}}^+(K)/\xi^e, \mathbb{Z}_p)_{\blacksquare} }.
    \]
    When $i=0$, there is nothing to prove as both sides are $\mathcal{R}$. For $i=1$, we use the cofiber sequence
    \begin{equation} \label{23049dsfhoi}
    \mathcal{R} \otimes_{(\mathbb{B}_{\on{dR}}^+(K)/\xi^e, \mathbb{Z}_p)_{\blacksquare} }
    \mathbb{L}_{(\mathbb{B}_{\on{dR}}^+(K)/\xi^e, \mathbb{Z}_p)_{\blacksquare} /(\mathbb{B}_{\on{dR}}^+(K), \mathbb{Z}_p)_{\blacksquare} }
    \longrightarrow
    \mathbb{L}_{ \mathcal{R}/(\mathbb{B}_{\on{dR}}^+(K), \mathbb{Z}_p)_{\blacksquare} }
    \longrightarrow \mathbb{L}_{\mathcal{R}/(\mathbb{B}_{\on{dR}}^+(K)/\xi^e, \mathbb{Z}_p)_{\blacksquare} }
    \end{equation}
    and we compute that 
    \[
      \mathcal{R} \otimes_{(\mathbb{B}_{\on{dR}}^+(K)/\xi^e, \mathbb{Z}_p)_{\blacksquare} }
    \mathbb{L}_{(\mathbb{B}_{\on{dR}}^+(K)/\xi^e, \mathbb{Z}_p)_{\blacksquare} /(\mathbb{B}_{\on{dR}}^+(K), \mathbb{Z}_p)_{\blacksquare} }
    \cong 
    (\xi^e)/(\xi^{2e}) 
     \otimes_{(\mathbb{B}_{\on{dR}}^+(K)/\xi^e, \mathbb{Z}_p)_{\blacksquare} }
     \mathcal{R}. 
    \]
    We have the commutative diagram
    \[
    \begin{tikzcd}
      (\xi^{e+i})/(\xi^{2e+2i}) 
     \otimes_{(\mathbb{B}_{\on{dR}}^+(K)/\xi^{e+1}, \mathbb{Z}_p)_{\blacksquare} }
     \mathcal{R}   
     \arrow[d]
     \arrow[r]
     &
      \mathbb{L}_{ \mathcal{R}/(\mathbb{B}_{\on{dR}}^+(K), \mathbb{Z}_p)_{\blacksquare} }
      \arrow[d,equal]
      \arrow[r]
      &
      \mathbb{L}_{\mathcal{R}/(\mathbb{B}_{\on{dR}}^+(K)/\xi^{e+i}, \mathbb{Z}_p)_{\blacksquare} }
      \arrow[d]
        \\
     (\xi^e)/(\xi^{2e}) 
     \otimes_{(\mathbb{B}_{\on{dR}}^+(K)/\xi^e, \mathbb{Z}_p)_{\blacksquare} }
     \mathcal{R}   
     \arrow[r]
     &
      \mathbb{L}_{ \mathcal{R}/(\mathbb{B}_{\on{dR}}^+(K), \mathbb{Z}_p)_{\blacksquare} }
      \arrow[r]
      &
      \mathbb{L}_{\mathcal{R}/(\mathbb{B}_{\on{dR}}^+(K)/\xi^{e}, \mathbb{Z}_p)_{\blacksquare} }
    \end{tikzcd}
    \]
    where the left vertical map is induced by natural inclusion $(\xi^{e+i}) \subset (\xi^e)$, and it follows that for $i \geq e$, the left vertical map is $0$. Then the limit of the left vertical map is $0$, and taking limits along the vertical maps of the commutative diagram, we obtain that 
    \[
    \mathbb{L}_{ \mathcal{R}/(\mathbb{B}_{\on{dR}}^+(K), \mathbb{Z}_p)_{\blacksquare} } 
    \cong 
    \on{lim}_e 
    \mathbb{L}_{\mathcal{R}/(\mathbb{B}_{\on{dR}}^+(K)/\xi^e, \mathbb{Z}_p)_{\blacksquare} },
    \]
    since the horizontal maps are cofiber sequences. 

    Taking $\on{Sym}^i_{\mathcal{R}}$ of the cofiber sequence (\ref{23049dsfhoi}), we obtain a natural filtration on $\on{Sym}^i_{\mathcal{R}} \mathbb{L}_{ \mathcal{R}/(\mathbb{B}_{\on{dR}}^+(K), \mathbb{Z}_p)_{\blacksquare} }  $
    with graded pieces 
    \[
     \on{Sym}^{j}_{\mathcal{R}}
     (
     \Sigma  (\xi^e)/(\xi^{2e}) 
     \otimes_{(\mathbb{B}_{\on{dR}}^+(K)/\xi^e, \mathbb{Z}_p)_{\blacksquare} }
     \mathcal{R} ) 
    \otimes_{\mathcal{R}}
    \on{Sym}^{i-j}_{\mathcal{R}} \Sigma \mathbb{L}_{\mathcal{R}/(\mathbb{B}_{\on{dR}}^+(K)/\xi^e, \mathbb{Z}_p)_{\blacksquare} } 
    \]
    taking limits along $e$, we see that
    $\on{Sym}^{j}_{\mathcal{R}}
     (
     \Sigma  (\xi^e)/(\xi^{2e}) 
     \otimes_{(\mathbb{B}_{\on{dR}}^+(K)/\xi^e, \mathbb{Z}_p)_{\blacksquare} }
     \mathcal{R} )$
     has limit $0$ for all $j \geq 1$ since the transition maps are eventually $0$ as we see above, so the only remaining term in the limit  is $ \on{lim}_e 
    \on{Sym}^i_{\mathcal{R}}
    \Sigma \mathbb{L}_{\mathcal{R}/(\mathbb{B}_{\on{dR}}^+(K)/\xi^e, \mathbb{Z}_p)_{\blacksquare} }$, hence
      \[ \on{Sym}^i_{\mathcal{R}} \Sigma \mathbb{L}_{ \mathcal{R}/(\mathbb{B}_{\on{dR}}^+(K), \mathbb{Z}_p)_{\blacksquare} } 
    \cong 
    \on{lim}_e 
    \on{Sym}^i_{\mathcal{R}}
    \Sigma \mathbb{L}_{\mathcal{R}/(\mathbb{B}_{\on{dR}}^+(K)/\xi^e, \mathbb{Z}_p)_{\blacksquare} }
    \]
    as desired.
\end{proof}
 
\begin{lemma} \label{f;zdfoweiuf}
    Let $f: (A, A^+) \rightarrow (B,B^+)$ be a map of complete Huber pairs over $(\mathbb{Q}_p, \mathbb{Z}_p)$, and $\widehat{\on{dR}}^{\on{an}}_{B/A}$ the analytic derived de Rham complex defined in \cite[construction 4.3.]{GuoLi2021}, then we have a canonical filtered equivalence
    \[
    \widehat{\on{dR}}_{(B,B^+)_{\blacksquare}/(A,A^+)_{\blacksquare}} (*)
    \cong 
    \widehat{\on{dR}}^{\on{an}}_{B/A}
    \]
\end{lemma}

\begin{proof}
    Recall that $ \widehat{\on{dR}}^{\on{an}}_{B/A} $
    is defined to be 
    \[
     \widehat{\on{dR}}^{\on{an}}_{B/A}= \on{colim}_{(B^0, A^0)} \widehat{\on{dR}}_{B^0/ A^0}^{\wedge} [\frac{1}{p}]
    \]
    where $A^0 $ and $B^0$ ranges through rings of definition of $(A,A^+)$ and $(B,B^+)$ such that $f (A^0) \subset B^0$, and $\on{colim}_{(B^0, A^0)} \widehat{\on{dR}}_{B^0/ A^0}^{\wedge}$ 
    is the derived $p$-adic completion of the Hodge filtered derived de Rham complex of $B^0$ over $A^0$, and the colimit is taken in the completely filtered category. Note that this is a filtered colimit. 
    
    We claim that there is a natural filtered equivalence
    \[
    \widehat{\on{dR}}_{(B,B^+)_{\blacksquare}/(A,A^+)_{\blacksquare}} 
    \cong 
    \on{colim}_{(B^0, A^0)}  \widehat{\on{dR}}_{(B,B^0)_{\blacksquare}/(A,A^0)_{\blacksquare}}
    \cong 
    \on{colim}_{(B^0, A^0)} \widehat{\on{dR}}_{\underline{B^0}_{\blacksquare}/ \underline{A^0}_{\blacksquare}} [\frac{1}{p}]
    \]
    where $\underline{B^0}$ and $\underline{A^0}$
    are the condensed ring associated to the topological ring $A^0$ and $B^0$ with $p$-adic topology. By definition of $(A, A^+)_{\blacksquare}$, we have 
    \[
    (A, A^+)_{\blacksquare} = \on{colim}_{A^0} (A, A^0)_{\blacksquare}
    \]
    as analytic rings, 
    where $A^0 \subset A^+$ ranges through ring of definitions of $(A, A^+)$, and similarly for $ (B, B^+)_{\blacksquare} $. 
   We observe that $\widehat{\on{dR}}_{(B,B^+)_{\blacksquare}/(A,A^+)_{\blacksquare}}$ 
    commutes with colimits on the first  variable $(B,B^+)_{\blacksquare}$. Indeed,  $\widehat{\on{dR}}_{(B,B^+)_{\blacksquare}/(A,A^+)_{\blacksquare}}$ is defined by taking the underlying condensed graded ring of $\widehat{\mathbb{DR}}_{(B,B^+)_{\blacksquare}/(A,A^+)_{\blacksquare}} $ as in definition \ref{zasdlfowei}, and then applying the functor $|-|^{\geq *}$ in (\ref{zcbksueife}), and  we know that  $\widehat{\mathbb{DR}}_{(B,B^+)_{\blacksquare}/(A,A^+)_{\blacksquare}} $ commutes with colimits in the variable $B$, which follows from   proposition \ref{423087403284}, and the description of colimits in $\on{GrAnDAG}$ in lemma \ref{weuryeiwyrue} and proposition \ref{wyrw3874}. Taking underlying rings commutes with sifted colimits in this setting, with the same proof as in \cite[corollary 4.1.12.]{Notesonsolidgeometry}. Since our colimit  is filtered and $|-|^{\geq *}$ clearly commutes with colimits, we are done.  We can now prove the claim using the change of base formula in lemma \ref{2830948302}. Indeed, we have filtered equivalences
    \begin{align*}
    \widehat{\on{dR}}_{(B,B^+)_{\blacksquare}/(A,A^+)_{\blacksquare}} 
     &
     \cong 
    \widehat{\on{dR}}_{(B,B^+)_{\blacksquare}/(\mathbb{Q}_p, \mathbb{Z}_p)_{\blacksquare}} 
    \coprod_{ \widehat{\on{dR}}_{(A,A^+)_{\blacksquare}/(\mathbb{Q}_p, \mathbb{Z}_p)_{\blacksquare}} }
    \underline{A}
    \\
  &   \cong 
  \on{colim}_{B^0}  \widehat{\on{dR}}_{(B,B^0)_{\blacksquare}/(\mathbb{Q}_p, \mathbb{Z}_p)_{\blacksquare}}
  \coprod_{  \on{colim}_{A^0}  \widehat{\on{dR}}_{(A,A^0)_{\blacksquare}/(\mathbb{Q}_p, \mathbb{Z}_p)_{\blacksquare}} }
    \underline{A}
    \\
   & \cong
    \on{colim}_{(B^0, A^0)} (  \widehat{\on{dR}}_{(B,B^0)_{\blacksquare}/(\mathbb{Q}_p, \mathbb{Z}_p)_{\blacksquare}}
  \coprod_{  \widehat{\on{dR}}_{(A,A^0)_{\blacksquare}/(\mathbb{Q}_p, \mathbb{Z}_p)_{\blacksquare}} }
    \underline{A} )
    \\
    & \cong  \on{colim}_{(B^0, A^0)} \widehat{\on{dR}}_{(B,B^0)_{\blacksquare}/(A,A^0)_{\blacksquare}}.
    \end{align*}

    We next observe that 
    $\widehat{\on{dR}}_{\underline{B^0}_{\blacksquare}/ \underline{A^0}_{\blacksquare}} $
    is already $p$-complete, this follows from looking at the graded pieces which are symmetric powers of $\Sigma \mathbb{L}_{\underline{B^0}_{\blacksquare}/ \underline{A^0}_{\blacksquare}}$
    by (\ref{zdfhweuoofd}). $\mathbb{L}_{\underline{B^0}_{\blacksquare}/ \underline{A^0}_{\blacksquare}}$ 
    is $p$-complete because $B^0$ is $p$-complete, and the symmetric powers are also $p$-complete by \cite[lemma 7.4.10.]{Notesonsolidgeometry}. Then we have a canonical map 
    $\widehat{\on{dR}}_{\underline{B^0}/ \underline{A^0}}^{\wedge}
    \longrightarrow \widehat{\on{dR}}_{\underline{B^0}_{\blacksquare}/ \underline{A^0}_{\blacksquare}} $, where we use that $\widehat{\on{dR}}_{\underline{B^0}_{\blacksquare}/ \underline{A^0}_{\blacksquare}}$ 
    is the solid completion of $\widehat{\on{dR}}_{\underline{B^0}/ \underline{A^0}}$. It is an equivalence because we can check it modulo $p$, which reduces to the discrete case and is clear. The lemma then follows by evaluating at the one point profinite set. 
\end{proof}

\begin{theorem} \label{faqper32089}
    Let $ (R, R^+) $ be a smooth complete Huber ring over $(K, K^+)$, which is a perfectoid field over $(\mathbb{Q}_p, \mathbb{Z}_p)$, and $(R,R^+)_{\blacksquare}$ the corresponding analytic animated rings over $\mathbb{Z}_{\blacksquare}$ as defined in \cite[theorem 3.2.8.]{Andreychev2021}, then there exists a canonical equivalence
    \[
    \widehat{\on{dR}}_{(R,R^+)_{\blacksquare}/(\mathbb{Q}_p,\mathbb{Z}_p)_{\blacksquare}} (*)
    \cong 
    R \Gamma_{\on{inf}}( \on{Spa}(R, R^+)/\Sigma )
    \]
    of $\on{B}_{\on{dR}}^+$-algebras,
    where $(*)$ denotes evaluation at the one point profinite set, and the right hand side is the infinitesmial cohomology constructed in \cite[theorem 1.2.7.]{Guo2025Crystalline}. The $\on{B}_{\on{dR}}^+$-algebra structure on the left hand side is through its $\widehat{\on{dR}}_{(K,K^+)_{\blacksquare}/(\mathbb{Q}_p,\mathbb{Z}_p)_{\blacksquare}}  $-algebra structure, which is equivalent to $\mathbb{B}_{\on{dR}}^+(K)$ by proposition \ref{afhsjfheosdf}. 
    
    Moreover, it is a filtered equivalence with the Hodge filtration on the left hand side and the infinitesmial filtration on the right hand side. The Katz-Oda filtration with respect to the $(K,K^+)$-algebra structure of $(R,R^+)$ on the left hand side, as defined in lemma \ref{2374032493}, corresponds to the $\xi$-adic filtration on the right hand side.  
\end{theorem}

\begin{proof}
It follows from 
lemma \ref{weruoincoe35} and lemma \ref{vcoerionvcd}  that 
\[
\widehat{\on{dR}}_{(R,R^+)_{\blacksquare}/(\mathbb{Q}_p,\mathbb{Z}_p)_{\blacksquare}} 
\cong 
\widehat{\on{dR}}_{(R,R^+)_{\blacksquare}/(\mathbb{B}_{\on{dR}}^+(K),\mathbb{Z}_p)_{\blacksquare}}
\cong 
\on{lim}_e \widehat{\on{dR}}_{(R,R^+)_{\blacksquare}/(\mathbb{B}_{\on{dR}}^+(K)/\xi^e,\mathbb{Z}_p)_{\blacksquare}}, 
\]
where the structure map $(\mathbb{B}_{\on{dR}}^+(K),\mathbb{Z}_p)_{\blacksquare} \longrightarrow (R,R^+)_{\blacksquare}$
is the composition of  the quotient map 
$(\mathbb{B}_{\on{dR}}^+(K),\mathbb{Z}_p)_{\blacksquare} \longrightarrow (K,K^+)_{\blacksquare} $ and the algebra map of $(R,R^+)$ over $(K,K^+)$. 
We know from \cite[corollary 7.2.6.]{Guo2025Crystalline} that 
\[
R \Gamma_{\on{inf}}( \on{Spa}(R, R^+)/\Sigma ) \cong \on{lim}_{e} \widehat{\on{dR}}^{\on{an}}_{R/\on{B}_{\on{dR}}^+/\xi^e}  
\]
 with $ \widehat{\on{dR}}^{\on{an}}_{R/\on{B}_{\on{dR}}^+/\xi^e}$ defined in \cite[construction 4.3.]{GuoLi2021}. It is also a filtered isomorphism carrying the  infinitesmial filtration to Hodge filtration using \cite[corollary 5.5.2.]{Guo2025Crystalline}, and we note that this is the only place where we use the smoothness assumption. 
The theorem then follows from 
    lemma \ref{f;zdfoweiuf}. The claim on Katz-Oda filtration is obvious. 
\end{proof}

\begin{theorem}
 Let $ (R, R^+) $ be a smooth complete Huber ring over $(K, K^+)$, which is a perfectoid field over $(\mathbb{Q}_p, \mathbb{Z}_p)$,  then there exists filtration on $\on{HH}((R,R^+)_{\blacksquare}/(\mathbb{Q}_p,\mathbb{Z}_p)_{\blacksquare})^{t\mathbb{S}^1}$
and $\on{HH}((R,R^+)_{\blacksquare}/(\mathbb{Q}_p,\mathbb{Z}_p)_{\blacksquare})^{h\mathbb{S}^1}$
with graded pieces
\[
\on{gr}^{i}
\on{HH}((R,R^+)_{\blacksquare}/(\mathbb{Q}_p,\mathbb{Z}_p)_{\blacksquare})^{t\mathbb{S}^1} (*)
= 
\Sigma^{2i}  R \Gamma_{\on{inf}}( \on{Spa}(R, R^+)/\Sigma )
\]
\[
\on{gr}^{i}
\on{HH}((R,R^+)_{\blacksquare}/(\mathbb{Q}_p,\mathbb{Z}_p)_{\blacksquare})^{h\mathbb{S}^1}(*)
= \Sigma^{2i} R \Gamma_{\on{inf}}( \on{Spa}(R, R^+)/\Sigma )^{\geq i}
\]
where $R \Gamma_{\on{inf}}( \on{Spa}(R, R^+)/\Sigma )^{\geq i}$
denote the $i$-th infinitesmial filtration of $R \Gamma_{\on{inf}}( \on{Spa}(R, R^+)/\Sigma )$
as defined in \cite[section 7.2.]{Guo2025Crystalline}. 

Moreover, they  are algebras over
\[
\on{gr}^{\bullet}\on{HH}((K,K^+)_{\blacksquare}/(\mathbb{Q}_p,\mathbb{Z}_p)_{\blacksquare})^{t\mathbb{S}^1} = \mathbb{B}_{\on{dR}}^+(K) [\sigma, \sigma^{-1}] 
\]
and 
\[
\on{gr}^{\bullet}\on{HH}((K,K^+)_{\blacksquare}/(\mathbb{Q}_p,\mathbb{Z}_p)_{\blacksquare})^{h\mathbb{S}^1} 
\cong \mathbb{B}_{\on{dR}}^+(K)[\Tilde{u}, v]/(\Tilde{u}v- \xi),
\]
respectively, 
and
the multiplication by $v$ 
induces the filtration map $R \Gamma_{\on{inf}}( \on{Spa}(R, R^+)/\Sigma )^{\geq i+1}\longrightarrow R \Gamma_{\on{inf}}( \on{Spa}(R, R^+)/\Sigma )^{\geq i}$,
and the multiplication map by $\Tilde{u}$ induces the map $R \Gamma_{\on{inf}}( \on{Spa}(R, R^+)/\Sigma )^{\geq i} \longrightarrow R \Gamma_{\on{inf}}( \on{Spa}(R, R^+)/\Sigma )^{\geq i+1} $ 
which is multiplication by $\xi$. 

Lastly, the de Rham comparison 
\[
R \Gamma_{\on{inf}}( \on{Spa}(R, R^+)/\Sigma )
\otimes_{\on{B}_{\on{dR}}^+} K
\cong
\widehat{\on{dR}}_{\on{Spa}(R, R^+)/\on{Spa}(K,K^+) }^{\on{an}}
\]
as established in \cite[theorem 1.2.7. (1)]{Guo2025Crystalline}
is induced by the graded pieces of the canonical isomorphism
\[
\on{HH}((R,R^+)_{\blacksquare}/(\mathbb{Q}_p,\mathbb{Z}_p)_{\blacksquare})^{t\mathbb{S}^1} 
\otimes_{\on{HH}((K,K^+)_{\blacksquare}/(\mathbb{Q}_p,\mathbb{Z}_p)_{\blacksquare})^{t\mathbb{S}^1}} 
\underline{K}^{ t\mathbb{S}^1}
\cong 
\on{HH}((R,R^+)_{\blacksquare}/(K,K^+)_{\blacksquare})^{t\mathbb{S}^1} .
\]
On the other hand, we have  
\[
\on{cofib}  \Bigl\{
R \Gamma_{\on{inf}}( \on{Spa}(R, R^+)/\Sigma )^{\geq i-1}
\overset{\xi}{\longrightarrow}
R \Gamma_{\on{inf}}( \on{Spa}(R, R^+)/\Sigma )^{\geq i}
 \Bigl\}
\cong
\widehat{\on{dR}}_{\on{Spa}(R, R^+)/\on{Spa}(K,K^+) }^{\on{an}, \geq i}
\]
i.e. the Hodge filtration of the de Rham cohomology is equivalent to the graded pieces of the "reversed" infinitesmial filtration,  
which is induced by the graded pieces of the natural equivalence
\[
\on{HH}((R,R^+)_{\blacksquare}/(\mathbb{Q}_p,\mathbb{Z}_p)_{\blacksquare})^{h\mathbb{S}^1} 
\otimes_{\on{HH}((K,K^+)_{\blacksquare}/(\mathbb{Q}_p,\mathbb{Z}_p)_{\blacksquare})^{h\mathbb{S}^1}} 
\underline{K}^{ h\mathbb{S}^1}
\cong 
\on{HH}((R,R^+)_{\blacksquare}/(K,K^+)_{\blacksquare})^{h\mathbb{S}^1}. 
\]    
\end{theorem}

\begin{proof}
    This follows directly from lemma \ref{qrehjglsg} 
    and theorem \ref{faqper32089}.
\end{proof}

\section{Comparison with Guo's prismatic $\on{B}_{\on{dR}}^+$-cohomology}

In this section, we prove comparison with the prismatic $\on{B}_{\on{dR}}^+$-cohomology constructed in \cite{Guo2026Prismatic} again in the smooth case. 
This cohomology theory coincides with the pro-etale $\on{B}_{\on{dR}}^+$-cohomology constructed in \cite{SCHOLZE2013} in the smooth case. 

Let us first recall the functor $\Psi$ constructed in \cite[section 4.1.]{Guo2026Prismatic}.  

\begin{definition}[{\cite[construction 4.1.1.]{Guo2026Prismatic}}]
Let $\on{M^{\bullet}} \in \on{Fil}^{\geq 0} D(\on{B}_{\on{dR}}^+) $, then we define 
\[
\Psi(M) := 
\on{colim}
\left\{
\begin{tikzcd}
M^0 \arrow[d, "\xi"]
\\
M^0
  & M^1 \arrow[l] \arrow[d, "\xi"]
\\
  & M^1
  & M^2 \arrow[l] \arrow[d, "\xi"]
\\
  & & \ddots
\end{tikzcd}
\right\} \in  D(\on{B}_{\on{dR}}^+)
\]
where the horizontal maps are the filtration maps. 
\end{definition}

By \cite[remark 4.1.3.]{Guo2026Prismatic}, $\Psi$
is lax symmetric monoidal, so induces a functor on filtered commutative algebras. 

On the other hand,  we construct a closely related functor as follows.  

\begin{definition}
    Let $N^{\bullet} \in \on{Mod}_{ (\on{B}_{\on{dR}}^+, (\xi)^{\bullet} ) } \on{Fil}^{\geq 0} D(\on{B}_{\on{dR}}^+) $, i.e. $N^{\bullet}$ is a filtered module over $\on{B}^+_{\on{dR}}$ with $\xi$-adic filtration, then we define 
    \[
    \Phi( N^{\bullet} ) := 
    \on{colim} 
\left\{ 
N^0 \overset{\xi}{\longrightarrow} N^1 \overset{\xi}{\longrightarrow} N^2 \overset{\xi}{\longrightarrow} \cdots
    \right\} 
    \in  D(\on{B}_{\on{dR}}^+). 
    \]
\end{definition}

\begin{lemma} \label{cvnosiruew}
    Let $\on{M^{\bullet}} \in \on{Fil}^{\geq 0} D(\on{B}_{\on{dR}}^+) $, then we have a canonical equivalence
    \[
    \Phi(
    M^{\bullet} \otimes_{\on{B}_{\on{dR}}^+} (\on{B}_{\on{dR}}^+, (\xi)^{\bullet}) )
    \cong 
    \Psi ( M^{\bullet}), 
    \]
    where the tensor product is the Day convolution tensor product of filtered objects, and the 
    $(\on{B}_{\on{dR}}^+, (\xi)^{\bullet})$-module structure on the tensor product is the one on the second factor. 
\end{lemma}

\begin{proof}
   Let $C_{i, j} = M^i\otimes_{\otimes_{\on{B}_{\on{dR}}^+}} (\xi^j)  $ be the diagram 
    \[
    C_{\bullet,\bullet}:= \left\{
    \begin{tikzcd}
  \vdots \arrow[d]  
  &
    \vdots \arrow[d] \arrow[l]
  &
   \vdots \arrow[d] \arrow[l]
   &
    \vdots \arrow[d] \arrow[l]
    \\
     M^0 \otimes_{\otimes_{\on{B}_{\on{dR}}^+}} (\xi^2)   
     \arrow[d]
     &
     M^1\otimes_{\on{B}_{\on{dR}}^+} (\xi^2)
     \arrow[l]
      \arrow[d]
     &
      M^2\otimes_{\on{B}_{\on{dR}}^+} (\xi^2)
     \arrow[l]
      \arrow[d]
     & 
     \cdots \arrow[l]
      \arrow[d]
    \\
        M^0 \otimes_{\otimes_{\on{B}_{\on{dR}}^+}} (\xi)   
         \arrow[d]
     &
     M^1\otimes_{\on{B}_{\on{dR}}^+} (\xi)
     \arrow[l]
      \arrow[d]
     &
      M^2\otimes_{\on{B}_{\on{dR}}^+} (\xi)
     \arrow[l]
      \arrow[d]
     & 
     \cdots \arrow[l]
      \arrow[d]
        \\
     M^0 \otimes_{\otimes_{\on{B}_{\on{dR}}^+}} (\xi^0)   
     &
     M^1\otimes_{\on{B}_{\on{dR}}^+} (\xi^0)
     \arrow[l]
     &
      M^2\otimes_{\on{B}_{\on{dR}}^+} (\xi^0)
     \arrow[l]
     & 
     \cdots \arrow[l]
    \end{tikzcd}
    \right\}
    \]
    where the horizontal maps are given by the filtration maps of $M^{\bullet}$, and the vertical maps are the inclusion of ideals. Then  by definition of Day convolution, the $n$-th filtration of 
    $ M^{\bullet} \otimes_{\on{B}_{\on{dR}}^+} (\on{B}_{\on{dR}}^+, (\xi)^{\bullet})$
    is the colimit of the sub-diagram $C_{i,j}|_{i+j \geq n}$, and it is easy to see that it is the colimit of  the cofinal sub-diagram 
    \[
    \begin{tikzcd}
        M^0\otimes_{\otimes_{\on{B}_{\on{dR}}^+}} (\xi^{n+1}) \arrow[d, ""] & & &
\\
M^0\otimes_{\otimes_{\on{B}_{\on{dR}}^+}} (\xi^{n})
  & M^1\otimes_{\otimes_{\on{B}_{\on{dR}}^+}} (\xi^{n}) \arrow[l] \arrow[d, ""]
  &
  &
\\
  &  \ddots & \ddots \arrow[d] &
  \\
  & &
  M^{n-1}\otimes_{\otimes_{\on{B}_{\on{dR}}^+}} (\xi^{})
  & 
  |[red]| M^n\otimes_{\otimes_{\on{B}_{\on{dR}}^+}} (\xi^{})
  \arrow[d, red] \arrow[l, red]
  \\
  & & & 
  |[red]|M^n\otimes_{\otimes_{\on{B}_{\on{dR}}^+}} (\xi^{0}).
    \end{tikzcd}
    \]
    The $ (\on{B}_{\on{dR}}^+, (\xi)^{\bullet})$-module  map  
    \begin{equation} \label{zvweorue9}
    \xi:  (M^{\bullet} \otimes_{\on{B}_{\on{dR}}^+} (\on{B}_{\on{dR}}^+, (\xi)^{\bullet}))^{n-1} 
    \longrightarrow
    (M^{\bullet} \otimes_{\on{B}_{\on{dR}}^+} (\on{B}_{\on{dR}}^+, (\xi)^{\bullet}))^{n}
    \end{equation}
    is induced by the map of diagrams $\xi: C_{i, j} \longrightarrow C_{i,j+1}$ which is multiplication by $\xi$ on the factor  $(\xi^j)$, so it 
    is  (after multiplication by $\xi$) the inclusion of the black sub-diagram above. If we fix trivializations of the ideals $(\xi^i) \cong \on{B}_{\on{dR}}^+$ by the generators $\xi^i$, the diagram becomes 
     \[
    \begin{tikzcd}
        M^0 \arrow[d, "\xi"] & & &
\\
M^0
  & M^1 \arrow[l] \arrow[d, "\xi"]
  &
  &
\\
  &  \ddots & \ddots \arrow[d,
  "\xi"] &
  \\
  & &
  M^{n-1}
  & 
  |[red]| M^n
  \arrow[d, red,"\xi"] \arrow[l, red]
  \\
  & & & 
  |[red]|M^n.
    \end{tikzcd}
    \]
    with the vertical maps multiplication by $\xi$, and the map (\ref{zvweorue9}) is simply the inclusion of the black sub-diagram. 
    It follows that after these trivializations of ideals, we have
    \[
    \Phi (M^{\bullet} \otimes_{\on{B}_{\on{dR}}^+} (\on{B}_{\on{dR}}^+, (\xi)^{\bullet}))
    = \on{colim}
\left\{
\begin{tikzcd}
M^0 \arrow[d, "\xi"]
\\
M^0
  & M^1 \arrow[l] \arrow[d, "\xi"]
\\
  & M^1
  & M^2 \arrow[l] \arrow[d, "\xi"]
\\
  & & \ddots
\end{tikzcd}
\right\} 
    \]
which is exactly $\Psi(M^{\bullet})$. 
\end{proof}

\begin{remark}
    The argument is formal, and the proof (and definitions) works verbatim for filtered condensed modules over $\mathbb{B}_{\on{dR}}^+(K)$ as well. We will use freely the condensed version below as well. 
\end{remark}

\begin{theorem} \label{cvoseufoe}
    Let $ (R, R_+) $ be a smooth complete Huber ring over  $(K, K^+)$, which is a perfectoid field over $(\mathbb{Q}_p, \mathbb{Z}_p)$, then there exists a canonical equivalence
    \begin{equation} \label{2370df}
    \widehat{\prism}_{(R,R^+)_{\blacksquare}/\mathbb{B}_{\on{dR}}^+(K)} (*)
    \cong 
    R \Gamma_{\prism}( \on{Spa}(R, R^+)/\on{B}_{\on{dR}}^+ )
    \end{equation}
    of $\on{B}_{\on{dR}}^+$-algebras,
    where the left hand side is defined as in lemma  \ref{92385093283fsd}, using $\mathbb{B}_{\on{dR}}^+(K) =\widehat{\on{dR}}_{(K,K^+)_{\blacksquare}/(\mathbb{Q}_p,\mathbb{Z}_p)_{\blacksquare}} $ 
    proved in proposition \ref{afhsjfheosdf}, and
    $(*)$ denotes evaluation at the one point profinite set. The right hand side is the prismatic cohomology constructed in \cite[theorem 1.1.2.]{Guo2026Prismatic}. 

    Moreover, if we modulo out by $\xi$ on both sides, then the Hodge-Tate filtration on the left hand side as defined in lemma \ref{2f937r03fs}
    corresponds to the Hodge-Tate filtration of the right hand side as defined in \cite[section 5]{Guo2026Prismatic}.
\end{theorem}

\begin{proof}
  By  \cite[corollary 4.3.15]{Guo2026Prismatic}, we have an  equivalence
  \begin{equation} \label{fw9eur382}
   R \Gamma_{\prism}( \on{Spa}(R, R^+)/\on{B}_{\on{dR}}^+ )
   \cong 
   \Psi(   \on{lim}_e  \widehat{\on{dR}}_{(P/(I,\xi^e))/\on{B}_{\on{dR}}^+/\xi^e}^{\on{an}} )
  \end{equation}
  where $P =\on{B}_{\on{dR}}^+\langle T_i \rangle$ is a completed polynomial ring together with a Koszul regular ideal $I$ such that $P/(I, \xi) \cong R$, and $P/I$ is flat over $\on{B}_{\on{dR}}^+$. Here $\Psi$ is defined in \cite[section 4.1.]{Guo2026Prismatic} and $\widehat{\on{dR}}_{(P/(I,\xi^e))/\on{B}_{\on{dR}}^+/\xi^e}^{\on{an}}$ is
  defined in \cite[construction 4.3.]{GuoLi2021}.  By \cite[example 4.3.2.]{Guo2026Prismatic}, Zariski locally on $ \on{Spa}(R, R^+) $
  there exists such a pair $(P,I)$, and we fix such a pair from now on.  Since both sides of (\ref{2370df}) satisfies Zariski descent, it is enough to prove the equivalence in this case.

We claim that  
\[
 \on{lim}_e  \widehat{\on{dR}}_{(P/(I,\xi^e), \mathbb{Z}_p)_{\blacksquare}/(\on{B}_{\on{dR}}^+/\xi^e, \mathbb{Z}_p)_{\blacksquare}}
 \cong 
 \widehat{\on{dR}}_{(P/I, \mathbb{Z}_p)_{\blacksquare}/(\on{B}_{\on{dR}}^+, \mathbb{Z}_p)_{\blacksquare}}.
\]
We have
\[
 \on{lim}_e  \mathbb{L}_{(P/(I,\xi^e), \mathbb{Z}_p)_{\blacksquare}/(\on{B}_{\on{dR}}^+/\xi^e, \mathbb{Z}_p)_{\blacksquare}}
 \cong 
 \mathbb{L}_{(P/I, \mathbb{Z}_p)_{\blacksquare}/(\on{B}_{\on{dR}}^+, \mathbb{Z}_p)_{\blacksquare}}
\]
because 
the right hand side is already $\xi$-adically complete which is a miracle of solidness. Now \cite[lemma 7.4.10.]{Notesonsolidgeometry}
shows that  $\on{Sym}^i_{(P/I, \mathbb{Z}_p)_{\blacksquare}}\Sigma \mathbb{L}_{(P/I, \mathbb{Z}_p)_{\blacksquare}/(\on{B}_{\on{dR}}^+, \mathbb{Z}_p)_{\blacksquare}}$
 is $\xi$-adic complete, and we have the corresponding limit expression for ever $i$, so the claim follows.

Moreover, we have
\[
\widehat{\on{dR}}^{\geq \bullet}_{(R, \mathbb{Z}_p)_{\blacksquare}/(\on{B}_{\on{dR}}^+, \mathbb{Z}_p)_{\blacksquare}}
\cong
\widehat{\on{dR}}^{\geq \bullet}_{(P/I, \mathbb{Z}_p)_{\blacksquare}/(\on{B}_{\on{dR}}^+, \mathbb{Z}_p)_{\blacksquare}} 
\otimes_{(\on{B}_{\on{dR}}^+, \mathbb{Z}_p)_{\blacksquare} }
\widehat{\on{dR}}^{\geq \bullet}_{(\on{B}_{\on{dR}}^+/\xi, \mathbb{Z}_p)_{\blacksquare}/(\on{B}_{\on{dR}}^+, \mathbb{Z}_p)_{\blacksquare}} 
\]
as filtered objects
because 
$
(R, \mathbb{Z}_p)_{\blacksquare} \cong (P/I, \mathbb{Z}_p)_{\blacksquare}  
\otimes_{(\on{B}_{\on{dR}}^+, \mathbb{Z}_p)_{\blacksquare} }
(\on{B}_{\on{dR}}^+/\xi, \mathbb{Z}_p)_{\blacksquare}
$
by our assumption on $(P,I)$
and Hodge completed  derived de Rham preserves colimits. 
We know 
\[
\widehat{\on{dR}}^{\geq \bullet}_{(\on{B}_{\on{dR}}^+/\xi, \mathbb{Z}_p)_{\blacksquare}/(\on{B}_{\on{dR}}^+, \mathbb{Z}_p)_{\blacksquare}}  \cong \mathbb{B}_{\on{dR}}^+(K)
\]
with $\xi$-adic filtration by lemma \ref{weruoincoe35}
and proposition \ref{afhsjfheosdf}, and it follows from 
  lemma \ref{cvnosiruew}
  that 
  \[
  \Phi(\widehat{\on{dR}}^{\geq \bullet}_{(R, \mathbb{Z}_p)_{\blacksquare}/(\on{B}_{\on{dR}}^+, \mathbb{Z}_p)_{\blacksquare}}  )
  \cong 
  \Psi(\widehat{\on{dR}}^{\geq \bullet}_{(P/I, \mathbb{Z}_p)_{\blacksquare}/(\on{B}_{\on{dR}}^+, \mathbb{Z}_p)_{\blacksquare}}).
  \]
  We know
  $
  \widehat{\prism}_{(R,R^+)_{\blacksquare}/\mathbb{B}_{\on{dR}}^+(K)}
  =
  \Phi(\widehat{\on{dR}}^{\geq \bullet}_{(R, \mathbb{Z}_p)_{\blacksquare}/(\on{B}_{\on{dR}}^+, \mathbb{Z}_p)_{\blacksquare}}  )$
  by its very definition in lemma  \ref{92385093283fsd},\footnote{ We use that $\widehat{\on{dR}}^{\geq \bullet}_{(R, R^+)_{\blacksquare}/(\on{B}_{\on{dR}}^+, \mathbb{Z}_p)_{\blacksquare}}  = \widehat{\on{dR}}^{\geq \bullet}_{(R, \mathbb{Z}_p)_{\blacksquare}/(\on{B}_{\on{dR}}^+, \mathbb{Z}_p)_{\blacksquare}}  $
here, it holds because 
$\mathbb{L}_{(R, R^+)_{\blacksquare}/(\on{B}_{\on{dR}}^+, \mathbb{Z}_p)_{\blacksquare} }
= \mathbb{L}_{(R, \mathbb{Z}_p)_{\blacksquare}/(\on{B}_{\on{dR}}^+, \mathbb{Z}_p)_{\blacksquare} }$
using our assumption of smoothness of $(R, R^+)$ over $(K,K^+)$. }
so we have proved that
\begin{equation} \label{cwoeruvdsf}
\widehat{\prism}_{(R,R^+)_{\blacksquare}/\mathbb{B}_{\on{dR}}^+(K)}
\cong
\Psi(\on{lim}_e  \widehat{\on{dR}}_{(P/(I,\xi^e), \mathbb{Z}_p)_{\blacksquare}/(\on{B}_{\on{dR}}^+/\xi^e, \mathbb{Z}_p)_{\blacksquare}}).
\end{equation}
The theorem then follows from (\ref{fw9eur382}) and lemma \ref{f;zdfoweiuf}. 

To prove the compatibility with Hodge-Tate filtration, recall that the Hodge-Tate filtration on 
$ R \Gamma_{\prism}( \on{Spa}(R, R^+)/\on{B}_{\on{dR}}^+ ) \otimes_{\on{B}_{\on{dR}}^+ } K$ 
is its Postnikov filtration by \cite[theorem 5.1.1.]{Guo2026Prismatic}. On the other hand, we know from lemma \ref{2f937r03fs} that the Hodge-Tate filtration on
\[
\widehat{\prism}_{(R,R^+)_{\blacksquare}/\mathbb{B}_{\on{dR}}^+(K)} \otimes_{(\mathbb{B}_{\on{dR}}^+(K),\mathbb{Z}_p)_{\blacksquare}} K 
\cong
\overline{\prism}_{(R,R^+)_{\blacksquare}/\mathbb{B}_{\on{dR}}^+(K)} 
\]
is complete and
has graded pieces $\Sigma^{-2*}\on{Sym}_{(R, \mathbb{Z}_p)_{\blacksquare}}^{\bullet} (\Sigma \mathbb{L}_{(R, \mathbb{Z}_p)_{\blacksquare}/(K, \mathbb{Z}_p)_{\blacksquare}}[1]) $, with each piece  concentrated in a different homological degree since $R$ is smooth over $K$, so it is also the Postnikov filtration.  
\end{proof}

\begin{theorem}
    Let $ (R, R_+) $ be a smooth complete Huber ring over  $(K, K^+)$, which is a perfectoid field over $(\mathbb{Q}_p, \mathbb{Z}_p)$, then there exists a canonical equivalence 
    \[
\on{gr}^{\bullet}\on{HH}((R,R^+)_{\blacksquare}/(\mathbb{Q}_p, \mathbb{Z}_p)_{\blacksquare})^{h\mathbb{S}^1}[\Tilde{u}^{-1}] (*) = 
 R \Gamma_{\prism}( \on{Spa}(R, R^+)/\on{B}_{\on{dR}}^+ )[\Tilde{u},\Tilde{u}^{-1}],  
\]
of graded $\on{B}_{\on{dR}}^+$-algebras,
where the graded $\on{B}_{\on{dR}}^+$-algebras structure on
the left hand side is given by the 
$\on{gr}^{\bullet}\on{HH}((K,K^+)_{\blacksquare}/(\mathbb{Q}_p, \mathbb{Z}_p)_{\blacksquare})^{h\mathbb{S}^1}   $-algebra structure, which is equivalent to $\mathbb{B}_{\on{dR}}^+(K)[\Tilde{u},v]/(\Tilde{u}v-\xi)$ 
by lemma \ref{qrehjglsg} and proposition \ref{afhsjfheosdf}. 

Moreover, we have
\[
\on{gr}_{\on{HKR}}^{\bullet}\on{HH}((R,R^+)_{\blacksquare}/(\mathbb{Q}_p, \mathbb{Z}_p)_{\blacksquare})[u^{-1}] (*)
\cong
(
R \Gamma_{\prism}( \on{Spa}(R, R^+)/\on{B}_{\on{dR}}^+ ) \otimes_{\on{B}_{\on{dR}}^+ } K) [u, u^{-1}]
\]
with $u$ of graded degree 1 and homological degree 2, 
and 
\[
\on{gr}_{\on{HKR}}^{i}\on{HH}((R,R^+)_{\blacksquare}/(\mathbb{Q}_p, \mathbb{Z}_p)_{\blacksquare})(*)
\cong
\Sigma^{2i}
F^{i}_{\on{HT}}
(
R \Gamma_{\prism}( \on{Spa}(R, R^+)/\on{B}_{\on{dR}}^+ ) \otimes_{\on{B}_{\on{dR}}^+ } K) 
\]
where $F^{i}_{\on{HT}}$ is the Hodge-Tate filtration as defined in \cite[theorem 5.1.1.]{Guo2026Prismatic}. 
\end{theorem}

\begin{proof}
    The combination of theorem \ref{cvoseufoe}, lemma \ref{92385093283fsd} and lemma \ref{2f937r03fs}. 
\end{proof}

\section{Breuil-Kisin cohomology}

Let $k$ be a perfect field of characteristic $p$, $K$ be a finite ramified extension of $W(k)[p^{-1}]$ with uniformizer $\omega$,  and $E$ be the minimal polynomial of $\omega$ over $ W(k)[p^{-1}]$. Let 
\[
\mathfrak{S} := \underline{W(k)[p^{-1}]}[T]^{\wedge}_{E} \in  \on{CAlg}D(\mathbb{Z}_{\blacksquare})
\]
where we give $W(k)[p^{-1}]$ the nonarchimedean topology. Note that we have 
$ \mathfrak{S} \cong 
K[[T-\omega]]$. 

\begin{proposition}
    We have a canonical equivalence
    \[
    \mathfrak{S} \cong 
    \widehat{\on{dR}}_{(K, \mathbb{Z}_p)_{\blacksquare}/(\mathbb{Q}_p[T], \mathbb{Z}_p)_{\blacksquare}}  
    \]
    where $ \mathbb{Q}_p[T] \longrightarrow K$ 
    sends $T$ to $\omega$. Moreover, this is a filtered isomorphism with the $E$-adic filtration on the left, and the Hodge filtration on the right. 
\end{proposition}

\begin{proof}
This is identical to the proof of proposition \ref{afhsjfheosdf}. 
We observe first that 
\begin{equation} \label{cv3ruvfoaiufe}
\mathbb{L}_{ (K, \mathbb{Z}_p)_{\blacksquare}/(\mathbb{Q}_p[T], \mathbb{Z}_p)_{\blacksquare}} 
\cong 
\mathbb{L}_{ (K, \mathbb{Z}_p)_{\blacksquare}/(W(k)[p^{-1}][T], \mathbb{Z}_p)_{\blacksquare}} 
\cong 
\Sigma
(E)/(E^2)
\cong 
\Sigma K
\end{equation}
    where the first equivalence follows from 
    \[ \mathbb{L}_{ (W(k)[p^{-1}][T], \mathbb{Z}_p)_{\blacksquare}/ (\mathbb{Q}_p[T], \mathbb{Z}_p)_{\blacksquare}
    } 
    \cong \mathbb{L}_{ (W(k)[T], \mathbb{Z}_p)_{\blacksquare}/ (\mathbb{Z}_p[T], \mathbb{Z}_p)_{\blacksquare}
    } [p^{-1}] 
    \cong
    0
    \]
which holds because $\mathbb{L}_{ (W(k)[T], \mathbb{Z}_p)_{\blacksquare}/ (\mathbb{Z}_p[T], \mathbb{Z}_p)_{\blacksquare}
    } $
    is already $p$-complete, and modulo $p$ we have 
    $\mathbb{L}_{k[T]/\mathbb{F}_p[T] } \cong \mathbb{L}_{k/\mathbb{F}_p} \otimes_{\mathbb{F}_p} \mathbb{F}_p[T] \cong 0$
    since $k$ is perfect. The rest follows directly from $K \cong  W(k)[p^{-1}][T]/E$, and $E$ is a non-zero divisor. Moreover, the computation also shows that the canonical differential 
    \[
    d: K \longrightarrow 
    \mathbb{L}_{ (K, \mathbb{Z}_p)_{\blacksquare}/(\mathbb{Q}_p[T], \mathbb{Z}_p)_{\blacksquare}} 
    \cong \Sigma (E)/(E^2)
    \]
    corresponds to the extension 
    \begin{equation} \label{85023kfdhj}
    0
    \longrightarrow
    (E)/(E^2)
    \longrightarrow
     W(k)[p^{-1}][T]/E^2 \longrightarrow   W(k)[p^{-1}][T]/E = K \longrightarrow 0. 
    \end{equation}

  We view $\mathfrak{S}$ with $E$-adic filtration as an element of $\on{Fil}^{\wedge}D((\mathbb{Q}_p[T],\mathbb{Z}_p)_{\blacksquare}) $, which defines an element of $ \on{cMod}_{\mathbb{D}_+^{\vee}} \on{Gr} D((\mathbb{Q}_p[T],\mathbb{Z}_p)_{\blacksquare}) $
    under the equivalence (\ref{zcbksueife}). It is also a commutative algebra object since the equivalence (\ref{zcbksueife}) is symmetric monoidal, so it defines an element of $\on{cMod}_{\mathbb{D}_+^{\vee}} (\on{DAlg} \on{Cond} \on{Gr} D(\mathbb{Z})_{\mathbb{Q}_p[T]/})  $
    since derived algebra objects coincides with commutative algebra objects in characteristic 0. The universal property in proposition \ref{423087403284} induces a map 
    \[
    L\Omega^{+,\bullet}_{\underline{K}/\underline{ \mathbb{Q}_p }[T]} \longrightarrow 
    \mathfrak{S}
    \]
    of $\on{cMod}_{\mathbb{D}_+^{\vee}} (\on{DAlg} \on{Cond} \on{Gr} D(\mathbb{Z})_{\mathbb{Q}_p[T]/})$
    which on degree 0 is the identity map of $K$. Give $\mathfrak{S}$ the analytic graded structure induced from its degree 0 part $K$ and its analytic structure $(K, \mathcal{O}_K)_{\blacksquare}$ using proposition \ref{wury9387r9}, the functoriality as in  
   lemma \ref{wityuoiewuroe} gives rise to a map 
     \[
    L\Omega^{+,\bullet}_{(K,\mathbb{Z}_p)_{\blacksquare}/(\mathbb{Q}_p[T], \mathbb{Z}_p)_{\blacksquare}} \longrightarrow 
   \mathfrak{S}.
    \]
Indeed,  lemma \ref{wityuoiewuroe} induces a map $\widehat{\mathbb{DR}}^{\on{un}}_{(K,\mathbb{Z}_p)_{\blacksquare}/(\mathbb{Q}_p[T], \mathbb{Z}_p)_{\blacksquare}} \longrightarrow \mathfrak{S}$
in $\on{GrAnDAG}^{\on{un}}$, passing to completion gives rise to 
$\widehat{\mathbb{DR}}_{(K,\mathbb{Z}_p)_{\blacksquare}/(\mathbb{Q}_p[T], \mathbb{Z}_p)_{\blacksquare}} \longrightarrow \mathfrak{S}$, where we use that 
$\mathfrak{S}$ is already complete as analytic graded derived ring since its underlying graded ring is $\on{Sym}^{\bullet}_{K} (\Sigma^2 K[1])$ that lies in $\on{Gr} D((K,\mathbb{Z}_p)_{\blacksquare})$. Taking the underlying condensed ring gives the desired map. 
    
    Since the underlying graded ring of both sides coincides with $\on{Sym}^{\bullet}_{K} (\Sigma^2 K[1])$ by lemma \ref{4702374903} and the definition of the equivalence (\ref{zcbksueife}), to prove it is an isomorphism it is enough to prove that it induces an isomorphism of degree 1, this is the same as in the proof proposition \ref{3074375372} where we observe that the $\mathbb{D}_+$-module structures on 1-truncated parts of both sides give rise to a commutative diagram 
    \[
    \begin{tikzcd}
     \mathbb{L}_{(K,\mathbb{Z}_p)_{\blacksquare}/(\mathbb{Q}_p[T], \mathbb{Z}_p)_{\blacksquare}} 
     \arrow[r]
     &
     \Sigma (E)/(E^2)
        \\
     K \arrow[u,"d"] \arrow[r,equal]  & K \arrow[u]
    \end{tikzcd}
    \]
    we have observed that the left hand side corresponds to the extension (\ref{85023kfdhj}) and it follows directly from the definition that the right hand side also corresponds to the same extension, so the horizontal arrow is an isomorphism by the universal property of cotangent complex. 

    Lastly, translating everything from comodules over $\mathbb{D}_{+}^{\vee}$ to completely filtered objects as in definition \ref{sfdsfjekl} gives the equivalence we want to prove.
    
\end{proof}

As in the case of $\on{B}_{\on{dR}}^+$-cohomology , we now define two $\mathfrak{S}$-cohomology theories of any analytic animated rings over perfectoid solid ring $(K,\mathbb{Z}_p)_{\blacksquare}$
with expected properties. 

\begin{definition}
    Let $(K,\mathbb{Z}_p)_{\blacksquare} \longrightarrow \mathcal{R}$
    be
    a map of analytic animated rings, then we define
    \[
    R\Gamma_{\on{inf}}(\mathcal{R}/\mathfrak{S} ):= \widehat{\on{dR}}_{\mathcal{R}/ (\mathbb{Q}_p[T], \mathbb{Z}_p)_{\blacksquare}}
    \]
    where the structure map is the composition
    $(\mathbb{Q}_p[T], \mathbb{Z}_p)_{\blacksquare}  \longrightarrow (K,\mathbb{Z}_p)_{\blacksquare} \longrightarrow  \mathcal{R}$
    with $T$ goes to $\omega$, and is naturally a $\mathfrak{S} $-algebra, with the algebra structure from algebra structure over $\widehat{\on{dR}}_{(K,\mathbb{Z}_p)_{\blacksquare}/ (\mathbb{Q}_p[T], \mathbb{Z}_p)_{\blacksquare}} \cong \mathfrak{S}$.  
    
   We equip it with the Hodge filtration
    \[
    R\Gamma_{\on{inf}}(\mathcal{R}/\mathfrak{S} )^{\geq \bullet}:= \widehat{\on{dR}}_{\mathcal{R}/ (\mathbb{Q}_p[T], \mathbb{Z}_p)_{\blacksquare}}^{\geq \bullet},
    \]
    and it is clear that it is compatible with the $E$-adic filtration on $\mathfrak{S}$. 
\end{definition}

\begin{proposition}
   Let $(K,\mathbb{Z}_p)_{\blacksquare} \longrightarrow \mathcal{R}$
   be
    a map of analytic animated rings,   there exists filtration on $\on{HH}((K, \mathbb{Z}_p)_{\blacksquare}/(\mathbb{Q}_p[T], \mathbb{Z}_p)_{\blacksquare})^{t\mathbb{S}^1}$
and $\on{HH}((K, \mathbb{Z}_p)_{\blacksquare}/(\mathbb{Q}_p[T], \mathbb{Z}_p)_{\blacksquare})^{h\mathbb{S}^1}$
with graded pieces
    \[
\on{gr}^{\bullet}
\on{HH}((K, \mathbb{Z}_p)_{\blacksquare}/(\mathbb{Q}_p[T], \mathbb{Z}_p)_{\blacksquare})^{t\mathbb{S}^1}
\cong
\mathfrak{S} [\sigma,\sigma^{-1}]
\]
with $\sigma$ of graded degree 1 and homological degree 2, and 
\[
\on{gr}^{\bullet}
\on{HH}((K, \mathbb{Z}_p)_{\blacksquare}/(\mathbb{Q}_p[T], \mathbb{Z}_p)_{\blacksquare})^{h\mathbb{S}^1}
\cong \mathfrak{S}[\Tilde{u},v]/(\Tilde{u}v-E)
\]
with $\Tilde{u}$ of homological degree 2 and graded degree 1. 
Moreover, the limit map 
\[
\mathfrak{S}[\Tilde{u}, v]/(\Tilde{u}v- \xi)
\cong 
\on{gr}^{\bullet}
\on{HH}((K, \mathbb{Z}_p)_{\blacksquare}/(\mathbb{Q}_p[T], \mathbb{Z}_p)_{\blacksquare})^{h\mathbb{S}^1}
\longrightarrow 
\on{gr}^{\bullet}\on{HH}((K, \mathbb{Z}_p)_{\blacksquare}/(\mathbb{Q}_p[T], \mathbb{Z}_p)_{\blacksquare})
\cong K [u]
\]
with $u$ the image of $\Tilde{u}$,
is simply quotient out by $v$. 
\end{proposition}

\begin{proof}
     Lemma \ref{242039htwo} and (\ref{cv3ruvfoaiufe}). 
\end{proof}

\begin{theorem}
    Let  $(K,\mathbb{Z}_p)_{\blacksquare} \longrightarrow \mathcal{R}$
  be  a map of analytic animated rings, then we have the de Rham comparison
    \[
     R\Gamma_{\on{inf}}(\mathcal{R}/\mathfrak{S} ) \otimes_{(\mathfrak{S},\mathbb{Z}_p)_{\blacksquare}  } K 
     \cong 
     \widehat{\on{dR}}_{\mathcal{R}/ (K, \mathbb{Z}_p)_{\blacksquare}}
    \]

    Moreover,  there exists filtration on $\on{HH}(\mathcal{R}/(\mathbb{Q}_p[T], \mathbb{Z}_p)_{\blacksquare})^{t\mathbb{S}^1}$
and $\on{HH}(\mathcal{R}/(\mathbb{Q}_p[T], \mathbb{Z}_p)_{\blacksquare})^{h\mathbb{S}^1}$
with graded pieces
\[
\on{gr}^{i}
\on{HH}(\mathcal{R}/(\mathbb{Q}_p[T], \mathbb{Z}_p)_{\blacksquare})^{t\mathbb{S}^1}
= 
\Sigma^{2i}  R\Gamma_{\on{inf}}(\mathcal{R}/\mathfrak{S} ) 
\]
\[
\on{gr}^{i}
\on{HH}(\mathcal{R}/(\mathbb{Q}_p[T], \mathbb{Z}_p)_{\blacksquare})^{h\mathbb{S}^1}
= \Sigma^{2i} R\Gamma_{\on{inf}}(\mathcal{R}/\mathfrak{S})^{\geq i}
\]
such that the canonical map from homotopy fixed points to Tate construction induces the filtration maps on graded pieces. 
\end{theorem}

\begin{proof}
The de Rham comparison is lemma \ref{2830948302}, and the rest is lemma \ref{qrehjglsg}.
\end{proof}

\begin{definition}
 Let $(K,K^+)_{\blacksquare} \longrightarrow \mathcal{R}$ be
    a map of analytic animated rings, then we define 
    \[
    \widehat{\prism}_{\mathcal{R}/\mathfrak{S} }: = \on{colim}( \widehat{\on{dR}}_{\mathcal{R}/(\mathbb{Q}_p[T], \mathbb{Z}_p)_{\blacksquare}}^{\geq 0} \overset{E}{\longrightarrow} \widehat{\on{dR}}_{\mathcal{R}/(\mathbb{Q}_p[T], \mathbb{Z}_p)_{\blacksquare}}^{\geq 1} 
    \overset{E}{\longrightarrow} \widehat{\on{dR}}_{\mathcal{R}/(\mathbb{Q}_p[T], \mathbb{Z}_p)_{\blacksquare}}^{\geq 2}
    \overset{E}{\longrightarrow}
    \cdots ),
    \]
    as in lemma \ref{92385093283fsd}, which is a $\mathfrak{S} $-algebra in solid modules,  and 
    \[
    \overline{\prism}_{\mathcal{R}/\mathfrak{S}  } 
    := \widehat{\prism}_{\mathcal{R}/\mathfrak{S} } \otimes_{ (\mathfrak{S}, \mathbb{Z}_p )_{\blacksquare}} K. 
    \]
\end{definition}

\begin{theorem}[Hodge-Tate filtration]
     Let $(K,\mathbb{Z}_p)_{\blacksquare} \longrightarrow \mathcal{R}$
   be  a map of analytic animated rings, then there exists a natural Hodge-Tate filtration
    \[
F_{\on{HT}}^{i}\overline{\prism}_{\mathcal{R}/\mathfrak{S}   }
    \cong 
     \Sigma^{-2i}
    \on{Sym}^i_{\mathcal{R}} \Sigma \mathbb{L}_{\mathcal{R}/(\mathbb{Q}_p[T], \mathbb{Z}_p)_{\blacksquare} }
    \]
    on $\overline{\prism}_{\mathcal{R}/\mathfrak{S}  } $
    with graded pieces
    \[
    \on{gr}_{\on{HT}}^{\bullet} \overline{\prism}_{\mathcal{R}/\mathfrak{S} }  
    \cong 
    \Sigma^{-2*}\on{Sym}_{\mathcal{R}}^{\bullet} (\Sigma \mathbb{L}_{\mathcal{R}/(K,\mathbb{Z}_p)_{\blacksquare}}[1])
    \]

  Moreover, we have
\[
\on{gr}_{\on{HKR}}^{\bullet}\on{HH}(\mathcal{R}/(\mathbb{Q}_p[T], \mathbb{Z}_p)_{\blacksquare})[u^{-1}]
\cong
\overline{\prism}_{\mathcal{R}/\mathfrak{S} } [u, u^{-1}]
\]
with $u$ in homological degree 2 and graded degree 1,
and 
\begin{equation} \label{30r538fsd}
\on{gr}_{\on{HKR}}^{i}\on{HH}(\mathcal{R}/(\mathbb{Q}_p[T], \mathbb{Z}_p)_{\blacksquare})
\cong
\Sigma^{2i}
F^{i}_{\on{HT}}\overline{\prism}_{\mathcal{R}/\mathfrak{S}  },
\end{equation}
and $u$ acts by the filtration map of Hodge-Tate filtration. Lastly, the limit map 
\[
\on{HH}(\mathcal{R}/(\mathbb{Q}_p[T], \mathbb{Z}_p)_{\blacksquare})^{h\mathbb{S}^1}
\longrightarrow
\on{HH}(\mathcal{R}/(\mathbb{Q}_p[T], \mathbb{Z}_p)_{\blacksquare})
\]
induces on graded objects
\begin{equation}  \label{cnwerwidjfi}
\on{gr}^i R\Gamma_{\on{inf}}(\mathcal{R}/\mathfrak{S} )
\cong 
F^{i}_{\on{HT}}\overline{\prism}_{\mathcal{R}/\mathfrak{S}  }
\end{equation}
\end{theorem}

\begin{proof}
    This is lemma \ref{2f937r03fs} and lemma \ref{92385093283fsd}. 
\end{proof}

\begin{theorem}
Let $(K,\mathbb{Z}_p)_{\blacksquare} \longrightarrow \mathcal{R}$
  be  a map of analytic animated rings, then there exists a canonical equivalence
    \[
    \on{gr}^{\bullet}
\on{HH}(\mathcal{R}/(\mathbb{Q}_p[T], \mathbb{Z}_p)_{\blacksquare})^{h\mathbb{S}^1}[\Tilde{u}^{-1}]
\cong 
\widehat{\prism}_{\mathcal{R}/\mathfrak{S} } [ \Tilde{u}, \Tilde{u}^{-1}]
    \]
with $\Tilde{u}$ of homological degree 2 and graded degree 1, and the natural limit map 
\[
\on{HH}
(\mathcal{R}/(\mathbb{Q}_p[T], \mathbb{Z}_p)_{\blacksquare})^{h\mathbb{S}^1} 
\longrightarrow 
\on{HH}(\mathcal{R}/(\mathbb{Q}_p[T], \mathbb{Z}_p)_{\blacksquare})
\]
induces after inverting $\Tilde{u}$
the quotient map 
$\widehat{\prism}_{\mathcal{R}/\mathfrak{S}  } 
\longrightarrow \overline{\prism}_{\mathcal{R}/\mathfrak{S}  } $
on  0-th graded. 
  
   Moreover, the localization map 
   \[
   \on{HH}(\mathcal{R}/(\mathbb{Q}_p[T], \mathbb{Z}_p)_{\blacksquare})^{h\mathbb{S}^1}
   \longrightarrow
   \on{HH}(\mathcal{R}/(\mathbb{Q}_p[T], \mathbb{Z}_p)_{\blacksquare})^{h\mathbb{S}^1}[\Tilde{u}^{-1}]
   \]
   induces a comparison map 
   \[
     R\Gamma_{\on{inf}}(\mathcal{R}/\mathfrak{S} ) 
     \longrightarrow 
     \widehat{\prism}_{\mathcal{R}/\mathfrak{S}  }
   \]
   and it becomes an equivalence after inverting $E$, i.e. 
   \[
     R\Gamma_{\on{inf}}(\mathcal{R}/ \mathfrak{S} ) [E^{-1}]
     \cong 
     \widehat{\prism}_{\mathcal{R}/\mathfrak{S}  }[E^{-1}].
   \]
\end{theorem}

\begin{proof}
  This is   lemma \ref{92385093283fsd}. 
\end{proof}

\begin{remark}
Throughout the discussion above,
     we have implicitly fixed  generator $E$ and $\Tilde{u} $
    to trivialize the Breuil-Kisin twists.
\end{remark}

\bibliographystyle{alpha} 
\bibliography{ref}
\end{document}